\PassOptionsToPackage{unicode}{hyperref}
\PassOptionsToPackage{naturalnames}{hyperref}
\documentclass[svgnames]{amsart}
\usepackage{latexsym,bm}
\usepackage[svgnames]{xcolor}
\usepackage{xparse}
\usepackage{longtable}
\usepackage{tikz}
\usetikzlibrary{matrix,decorations.markings,decorations.text,calc}
\usetikzlibrary{decorations.pathreplacing,calligraphy}
\usepackage[utf8]{inputenc}
\usepackage[english]{babel}
\usepackage[T1]{fontenc}
\usepackage{amsmath, amsfonts, amssymb, amsthm, mathrsfs, pb-diagram}
\usepackage{lmodern}
\usepackage{fullpage}
\usepackage{microtype}
\usepackage{mathtools}
\usepackage{enumitem}
\setlist[enumerate]{
  label=\alph*\upshape),
  ref=\alph*,
  nolistsep
}

\usepackage{todonotes}

\NewDocumentCommand\Cmd{ sm }{\textsf{\textbackslash #2}\IfBooleanT{#1}{$\{\ldots\}$}}

\usepackage{etoolbox}
\usepackage{aliascnt}
\usepackage{cite}

\usepackage[pagebackref=false]{hyperref}
\hypersetup{%
  pdftitle={Skew cellularity of the Hecke algebras of type $G(\ell,p,n)$},
  pdfauthor={Jun Hu, Andrew Mathas and Salim Rostam},%
  colorlinks = true,%
  linkcolor =blue,%
  anchorcolor = red,%
  citecolor = blue,%
  urlcolor = blue%
}

\DeclareMathOperator\cha{\textrm{char}}
\renewcommand\bot{\mathop{\textrm{bot}}}
\renewcommand\top{\mathop{\textrm{top}}}

\newcommand\Std{\mathop{\rm Std}\nolimits}
\newcommand{\DeclareMyOperator}[1]{%
  \expandafter\DeclareMathOperator\csname #1\endcsname{#1}
}
\forcsvlist{\DeclareMyOperator}{rad, cross, id, res, Hom, Shape, Left, Right,reg, sing, Rep}

\newcommand\1{\mathbf{1}}

\DeclarePairedDelimiterX{\set}[1]{\{}{\}}{\setargs{#1}}
\NewDocumentCommand{\setargs}{>{\SplitArgument{1}{|}}m}{\setargsaux#1}
\NewDocumentCommand{\setargsaux}{mm}
{\IfNoValueTF{#2}{#1} {#1\,\delimsize|\,\mathopen{}#2}}

\def\map#1#2{\,{:}\,#1\!\longrightarrow\!#2}
\newcommand\bijection{\xrightarrow{\ \sim\ }}

\let\gedom=\trianglerighteq
\let\gdom=\vartriangleright
\let\ledom=\trianglelefteq
\let\ldom=\vartriangleleft

\newcommand\tdom{\mathop{\vartriangleright_{\dcharge}}}
\newcommand\tedom{\mathop{\trianglerighteq_{\dcharge}}}

\newcommand\restr[1][A^\sigma]{\!\downarrow}
\newcommand\ind[1][A]{\!\uparrow}

\newcommand\Uglov[1][n]{\mathscr{U}_{\brho,#1}^{\bLam}}
\newcommand\Nscr{{\mathscr{N}}}
\newcommand\Nodes[1][n]{\Nscr^\ell_{#1}}
\newcommand\Parts[1][n]{\mathscr{P}^\ell_{#1}}
\newcommand\Partss{\mathscr{P}^{\ell}_{\sigma,n}}
\newcommand\Partsp[1][n]{\mathscr{P}^\ell_{\sigma,p,#1}}
\newcommand\Partspz[1][n]{\mathscr{P}^\ell_{\sigma,p,#1,0}}
\newcommand\Comp[1][n]{\mathscr{C}^\ell_{#1}}
\newcommand\Ccal{\mathcal{C}}

\newcommand\Fcal{\mathcal{F}}
\newcommand\Ocal{\mathcal{O}}
\newcommand\Qcal{\mathcal{Q}}
\renewcommand\P{\mathcal{P}}
\newcommand\Psig{{\P_{\sigma,\sp}}}

\newcommand\blam{{\boldsymbol\lambda}}
\newcommand\bLam{{\boldsymbol\Lambda}}
\newcommand\brho{{\boldsymbol\rho}}
\newcommand\bmu{{\boldsymbol\mu}}
\newcommand\bnu{{\boldsymbol\nu}}
\newcommand\bom{{\boldsymbol\omega}}

\renewcommand\sp{\mathsf{p}}

\newcommand\plam{\sp_\lambda}
\newcommand\so{\mathsf{o}}

\newcommand\sigmatab[1][\lambda]{\sigma_{#1}}

\newcommand\sigmaclam[1][\lambda]{\sigma_{C,#1}}
\newcommand\compgamma[1][j]{\gamma_{\lambda,#1}}
\newcommand\sigmavg{\overline{\sigma}_A} 
\newcommand\Prep{\P_\sigma}

\newcommand\Tsig{{T_{\sigma}}}
\newcommand\lesig{\ldom_\sigma}
\newcommand\gesig{\gdom_\sigma}
\newcommand\leqsig{\ledom_\sigma}

\newcommand\olambdatab[1][\lambda]{\so_T(#1)}
\newcommand\region[2]{\mathsf{X}_{\dcharge}^{#1,#2}}
\newcommand\lex{<_\mathrm{lex}}

\newcommand{\dia}{\delta} 
\newcommand\regionnodes[1][\dia]{\mathscr{N}^{a,#1}}
\newcommand\regionblamfull[2]{\blam^{#1,#2}}
\NewDocumentCommand\regionblam{ O{\dia} O{a} }{\regionblamfull{#2}{#1}}
\newcommand\Ckt[1][\t]{B^{[k], \reg}_{\blam, #1}}
\newcommand\Cslam[1][k]{{B_{\blam}^{[#1]}}}
\newcommand\Cslamreg[1][k]{B_{\blam}^{[#1], \reg}}
\newcommand\Cslamsing[1][k]{B_{\blam}^{[#1], \sing}}
\newcommand\blamr[1][k]{\blam_R^{[#1]}}
\newcommand\blaml[1][k]{\blam_L^{[#1]}}
\newcommand\Stds[1][\blam]{\mathrm{Std}_\sigma(#1)}
\newcommand\Btreg[1][\t]{B_{#1}^{\reg}}
\newcommand\Ctreg[1][\t]{C_{#1}^{\reg}}
\newcommand\Btsing[1][\t]{B_{#1}^{\sing}}
\newcommand\Ctsing[1][\t]{B_{#1}^{\sing}}

\newcommand\Dkl[1][]{D_{\sigma^{k#1}\s,\sigma^{l#1}\t}^{\text{reg}}}
\newcommand\gbconf{\widehat\bnu}

\newcommand\Dkt[1][\t]{D^{(k)}_{#1}}

\let\emph\textbf

\newcommand\sig[1][]{\sigma{#1}}

\newcommand\tsig{\tilde\sigma}

\NewDocumentCommand\Csmu{ s }{
  \IfBooleanTF{#1}{C_{\bmu,\sigma}}{C_{\sigma,\bmu}}
}

\newcommand\Eom{\mathsf{E}_{\bom}}
\newcommand\Eoma[1][\alpha]{\mathsf{E}_{\bom,#1}}
\newcommand\EomA[1][\alpha]{\mathsf{E}_{\bom,[#1]}}

\newcommand{\str}[1][s]{\mathsf{#1}}
\newcommand{\gstr}[1][g]{\mathsf{#1}}
\newcommand{\rstr}{\mathsf{r}}

\let\<=\langle
\let\>=\rangle

\def\N{\mathbb{N}}
\def\Q{\mathbb{Q}}
\def\R{\mathbb{R}}
\def\Z{\mathbb{Z}}

\def\a{\mathfrak{a}}
\def\b{\mathfrak{b}}
\def\b{\mathfrak{b}}

\def\s{\mathfrak{s}}
\def\t{\mathfrak{t}}
\def\u{\mathfrak{u}}
\def\v{\mathfrak{v}}

\newcommand\Add{{\mathcal{A}}}
\newcommand\Rem{{\mathcal{R}}}

\let\Sec\S

\def\lam{\lambda}

\def\Sym{\mathfrak{S}}

\let\eps\varepsilon
\newcommand\epsl[1][\lambda]{\eps_{#1}}
\newcommand\zl[1][\lambda]{z}
\newcommand\NN{\tfrac1N }

\newcommand{\bQ}{\mathbf{Q}}
\newcommand{\bvQ}{\bQ^{\vee\varepsilon}}

\newcommand\noedge{\mathop{\rlap{\space\slash}\rule[0.5ex]{1.1em}{0.1ex}}}

\newcommand\Gln{G(\ell,1,n)}
\newcommand\Glpn{G(\ell,p,n)}
\newcommand\RR{\mathscr{R}}
\NewDocumentCommand\RG{O{n} d()}{\RR_{#1}^{\bLam\IfNoValueF{#2}{,#2}}}
\newcommand\RGa[1][\alpha]{\RR_{[#1]}^{\bLam}}

\NewDocumentCommand\Rlpn{ O{n} d()}{\RR^{\bLam\IfNoValueF{#2}{,#2}}_{p,#1}}
\newcommand\HH{\mathscr{H}}
\NewDocumentCommand\Hln{D(){\q,\bvQ}}{\HH_{n}(#1)}
\NewDocumentCommand\Hlpn{O{p,n} D(){\q,\bvQ} D<>{} }{\HH^{#3}_{#1}(#2)}
\NewDocumentCommand\WA{ O{n} }{\mathbb{A}_{#1}^{\dcharge}}
\NewDocumentCommand\WAA{ O{\alpha} }{\mathbb{A}_{[#1]}^{\dcharge}}
\newcommand\q{q}
\newcommand\cyc{{\operatorname{\mathrm{cyc}}}}

\newcommand\Web{\mathscr{W}_{\dcharge}}

\newcommand\xcoord[1][\dcharge]{\mathsf{x}_{#1}}

\newcommand\ltx{<_{\mathsf{x}}}

\newcommand\charge{\boldsymbol\kappa}
\newcommand\dcharge{\boldsymbol\brho}

\newcommand\bi{\mathbf{i}}
\newcommand\bj{\mathbf{j}}

\newcommand\cI{\mathcal{I}}
\newcommand\cIa{\mathcal{I}^\alpha}
\newcommand\cIsig{\mathcal{I}^n_\sigma}

\newcommand\enumref[2]{\hyperref[#2]{\autoref*{#1}(\autoref*{#2})}}

\def\NewTheorem#1{%
  \newaliascnt{#1}{equation}%
  \newtheorem{#1}[#1]{#1}%
  \aliascntresetthe{#1}%
  \expandafter\def\csname #1autorefname\endcsname{#1}%
}

\newtheorem*{MainTheorem}{Main Theorem}
\numberwithin{equation}{section}

\NewTheorem{Proposition}
\NewTheorem{Theorem}
\NewTheorem{Corollary}
\NewTheorem{Lemma}
\NewTheorem{Definition}
\theoremstyle{definition}
\NewTheorem{Notation}
\NewTheorem{Convention}
\NewTheorem{Example}
\AtEndEnvironment{Example}{\null\hfill$\Diamond$}%
\NewTheorem{Examples}
\AtEndEnvironment{Examples}{\vskip-10mm\null\hfill$\Diamond$}%
\theoremstyle{remark}
\NewTheorem{Remark}
\NewTheorem{Assumption}

\newcounter{theoremcases}
\AtBeginEnvironment{proof}{\setcounter{theoremcases}{0}}
\newcommand\Case[1]{%
  \refstepcounter{theoremcases}%
  \noindent\textbf{Case \arabic{theoremcases}.}\space#1:
}

\ExplSyntaxOn
\prop_new:N \l_charge_prop
\prop_new:N \l_tab_prop 
\prop_new:N \l_bom_prop 

\seq_new:N  \l_tab_seq  
\clist_new:N  \l_word_clist 

\int_new:N \l_m_int     
\int_new:N \l_l_int     
\int_new:N \l_r_int     
\int_new:N \l_c_int     
\int_new:N \l_ell_int   
\int_new:N \l_n_int     

\fp_new:N \l_top_fp     
\fp_new:N \l_bot_fp     

\int_new:N \l_N_int     
\int_set:Nn \l_N_int {10}

\fp_new:N  \l_xcoord_fp 
\tl_new:N  \l_xcoord_tl 

\tl_new:N  \l_charge_tl 

\cs_generate_variant:Nn \prop_put:Nnn {NVx}
\cs_new:Npn \set_charge:nnn #1#2#3
{ 
  \prop_clear:N \l_charge_prop
  \clist_set:Nn \l_tmpa_clist {#3}
  \int_set:Nn \l_ell_int { \clist_count:N \l_tmpa_clist }
  \int_zero:N \l_l_int
  \clist_map_inline:Nn \l_tmpa_clist
  {
    \fp_set:Nn \l_xcoord_fp { ##1/#2 - \l_l_int/(#2*(\l_ell_int+1)) }
    \int_incr:N \l_l_int
    \prop_put:NVx \l_charge_prop \l_l_int { \fp_to_decimal:N \l_xcoord_fp }
    \draw[redstring] (\fp_to_decimal:N \l_xcoord_fp,0)--+(0,#1);
  }
}

\cs_new:Npn \set_xcoord:nnn #1#2#3
{
  \prop_get:NoN \l_charge_prop {#3} \l_charge_tl 
  \fp_set:Nn \l_xcoord_fp { \l_charge_tl+#2-#1 - (#1+#2)/\l_N_int }
  \tl_set:Nx \l_xcoord_tl {\fp_to_decimal:N\l_xcoord_fp}
}

\cs_generate_variant:Nn \prop_put:Nnn {NVo}
\cs_new:Npn \set_tab_prop:n #1
{
  \prop_clear:N \l_tab_prop
  \seq_clear:N \l_tab_seq
  \int_zero:N \l_l_int 
  \int_zero:N \l_n_int 
  \clist_set:Nn \l_tmpa_clist {#1}
  \clist_map_inline:Nn \l_tmpa_clist
  {
   \int_incr:N \l_l_int
   \prop_get:NVN \l_charge_prop \l_l_int \l_charge_tl 
   \int_zero:N \l_r_int 
   \clist_set:Nn \l_tmpb_clist {##1}
   \clist_map_inline:Nn \l_tmpb_clist
   {
     \int_incr:N \l_r_int
     \int_zero:N \l_c_int 
     \clist_set:Nn \l_tmpc_clist {####1}
     \clist_map_inline:Nn \l_tmpc_clist
     {
       \int_incr:N \l_c_int
       \set_xcoord:nnn {\l_r_int} {\l_c_int} {\int_use:N\l_l_int}
       \seq_put_right:NV \l_tab_seq \l_xcoord_tl
       \prop_put:NVn \l_tab_prop \l_xcoord_tl {########1}

       \int_incr:N \l_n_int
       \set_xcoord:nnn {\l_n_int} {1} {1}
       \prop_put:NVV \l_bom_prop \l_n_int \l_xcoord_tl
     }
   }
  }
  \seq_sort:Nn \l_tab_seq
  {
    \fp_compare:nNnTF { ##1 } < { ##2 }
    { \sort_return_swapped: }
    { \sort_return_same: }
  }
}

\NewDocumentCommand\DrawTableauDiagram{ s O{1} m D(){0} m }
{
  \IfBooleanTF{#1}
    { 
      \fp_set:Nn \l_top_fp {0}
      \fp_set:Nn \l_bot_fp {#2}
    }
    {
      \fp_set:Nn \l_top_fp {#2}
      \fp_set:Nn \l_bot_fp {0}
    }

  \set_charge:nnn{#2}{#3}{#4}
  \set_tab_prop:n{#5}

  \int_zero:N \l_m_int
  \seq_map_inline:Nn \l_tab_seq {
    \int_incr:N \l_m_int
    \prop_get:NVN \l_bom_prop \l_m_int \l_tmpa_tl
    \prop_get:NnN \l_tab_prop {##1} \l_tmpb_tl
    \prop_get:NVN \l_bom_prop \l_tmpb_tl \l_tmpc_tl
    \draw[solid](\tl_use:N##1,\fp_to_decimal:N\l_top_fp)
              --(\tl_use:N\l_tmpa_tl,0.5*#2)
              --(\tl_use:N\l_tmpc_tl,\fp_to_decimal:N\l_bot_fp);
    \draw[ghost](\tl_use:N##1+1,\fp_to_decimal:N\l_top_fp)
              --(\tl_use:N\l_tmpa_tl+1,0.5*#2)
              --(\tl_use:N\l_tmpc_tl+1,\fp_to_decimal:N\l_bot_fp);
  }
}

\NewDocumentCommand\DrawOmegaDiagram{ O{0.6} m D(){0} m m }
{
  \int_set:Nn \l_n_int {#4}
  \clist_set:Nn \l_word_clist {#5}

  \set_charge:nnn{#1*\clist_count:N\l_word_clist}{#2}{#3}

  \int_zero:N \l_m_int
  \prop_clear:N \l_bom_prop
  \int_do_while:nn { \l_m_int < \l_n_int }
  {
     \int_incr:N \l_m_int
     \set_xcoord:nnn {\l_m_int} {1} {1}
     \prop_put:NVV \l_bom_prop \l_m_int \l_xcoord_tl
  }

  \fp_set:Nn \l_top_fp {#1}
  \fp_zero:N \l_bot_fp

  \int_zero:N \l_m_int
  \clist_map_inline:Nn \l_word_clist {
    \int_zero:N \l_m_int
    \int_do_while:nn { \l_m_int < \l_n_int }
    { 
        \int_incr:N \l_m_int
        \prop_get:NVN \l_bom_prop \l_m_int \l_tmpa_tl
        \int_compare:nNnTF {\l_m_int} = {##1}
        { 
            \int_incr:N \l_m_int
            \prop_get:NVN \l_bom_prop \l_m_int \l_tmpb_tl
            \draw[solid](\tl_use:N\l_tmpa_tl,\fp_to_decimal:N\l_bot_fp)
                        --(\tl_use:N\l_tmpb_tl,\fp_to_decimal:N\l_top_fp);
            \draw[solid](\tl_use:N\l_tmpb_tl,\fp_to_decimal:N\l_bot_fp)
                        --(\tl_use:N\l_tmpa_tl,\fp_to_decimal:N\l_top_fp);
            \draw[ghost](\tl_use:N\l_tmpa_tl+1,\fp_to_decimal:N\l_bot_fp)
                        --(\tl_use:N\l_tmpb_tl+1,\fp_to_decimal:N\l_top_fp);
            \draw[ghost](\tl_use:N\l_tmpb_tl+1,\fp_to_decimal:N\l_bot_fp)
                        --(\tl_use:N\l_tmpa_tl+1,\fp_to_decimal:N\l_top_fp);
        }
        { 
            \draw[solid](\tl_use:N\l_tmpa_tl,\fp_to_decimal:N\l_bot_fp)
                        --(\tl_use:N\l_tmpa_tl,\fp_to_decimal:N\l_top_fp);
            \draw[ghost](\tl_use:N\l_tmpa_tl+1,\fp_to_decimal:N\l_bot_fp)
                        --(\tl_use:N\l_tmpa_tl+1,\fp_to_decimal:N\l_top_fp);
        }
    }
    \fp_add:Nn \l_top_fp {#1}
    \fp_add:Nn \l_bot_fp {#1}
  }
}

\NewDocumentCommand\DrawWebsterIdempotent{ O{1} m D(){0} m }
{
  \set_charge:nnn{#1}{#2}{#3}
  \int_zero:N \l_l_int 
  \clist_set:Nn \l_tmpa_clist {#4}
  \clist_map_inline:Nn \l_tmpa_clist
  {
    \int_incr:N \l_l_int
    \int_zero:N \l_r_int 
    \clist_set:Nn \l_tmpb_clist {##1}
    \clist_map_inline:Nn \l_tmpb_clist
    {
     \int_incr:N \l_r_int
     \int_step_inline:nnn 1 {####1}
     {
       \set_xcoord:nnn {\l_r_int} {########1} {\int_use:N\l_l_int}
       \draw[solid](\fp_to_decimal:N \l_xcoord_fp, 0)--++(0,#1);
       \draw[ghost](\fp_to_decimal:N \l_xcoord_fp+1, 0)--++(0,#1);
     }
    }
  }
}
\ExplSyntaxOff

\tikzset{
  baseline=(current bounding box.center),
  centered/.style = {
     baseline = {([yshift=#1]current bounding box.center)}
  },
  centered/.default=0ex,
  dots/.style={line width=1pt, line cap=round, gray,
               dash pattern=on 0pt off 2\pgflinewidth},
  frame/.style={fill=FloralWhite, draw=BurlyWood,very thin},
  circled/.style = {fill=Tan},
  redstring/.style = {line width=0.8mm, red},
  semisolid/.style = {draw=DarkBlue, thin},
  solid/.style = {draw=DarkBlue, line width=0.4mm},
  ghost/.style = {draw=LightSlateGrey, dashed, line width=0.4mm},
  dot/.style = {
    decoration={markings, mark=at position #1 with {\node[circle, inner sep=1.8pt, fill=DarkBlue]{};}},
    postaction={decorate}
  },
  dot/.default=0.5,
  ghostdot/.style = {
    decoration={markings, mark=at position #1 with {\node[circle, inner sep=1.8pt, fill=LightSlateGrey]{};}},
    postaction={decorate}
  },
  ghostdot/.default=0.5,
  pics/dot/.style = {
    code={\node[fill=LightSlateGrey,circle,inner sep=1.8pt,outer sep=0]at(0,0){};}
  },
  pics/TableauDiagram/.style={
    code={ \DrawTableauDiagram #1 }
  },
  pics/OmegaDiagram/.style={
    code={ \DrawOmegaDiagram #1 }
  },
  pics/WebsterIdempotent/.style={
    code={ \DrawWebsterIdempotent #1 }
  },
  wall/.style={black!80, thick},
  innerwall/.style={gray, thin, fill=white},
  russianTableau/.style={
    scale=0.2,
    draw/.append style={thick,black}
  },
  pics/RussianTableau/.style 2 args={
    code = { \begin{scope}[rotate=10, #1]
               \def\lastC{0}
               \foreach \row [count=\r] in {#2} {
                   \foreach \ent [count=\c] in \row {
                      \coordinate(-\r-\c) at (\c-\r,\c+\r-2) {};
                      \draw[innerwall] (\c-\r,\c+\r-2) --++(1,1)--++(-1,1)--++(-1,-1)--cycle;
                      \node at (\c-\r,\c+\r-1) {\ent};
                      \xdef\cc{\c}
                   }
                   \ifnum\r=1\draw[wall](0,0)--(\cc,\cc);\fi
                   \ifnum\lastC>\cc\draw[wall](\cc-\r+1,\cc+\r-1)--(\lastC-\r+1,\lastC+\r-1);\fi
                   \draw[wall](1-\r,\r-1)--++(-1,1);
                   \draw[wall](\cc-\r+1,\cc+\r-1)--++(-1,1);
                   \xdef\lastC{\cc}
                   \xdef\rr{\r}
               }
               \draw[wall](-\rr,\rr)--++(\lastC,\lastC);
            \end{scope}
    }
  }
}

\newcommand\RussianTableau[2][]{%
  \begin{tikzpicture}[russianTableau,#1]
    \draw (0,0) pic{RussianTableau={scale=0.4}{#2}};
  \end{tikzpicture}
}

\newcommand\RussianTableauscale[2][]{%
  \begin{tikzpicture}[russianTableau,#1]
    \draw (0,0) pic{RussianTableau={scale=0.3}{#2}};
  \end{tikzpicture}
}

\NewDocumentCommand\CtDiagram{ s O{1.5} m D(){0} m }
{
  \begin{tikzpicture}[centered=-1ex]
    \IfBooleanTF{#1}
    { \draw(0,0) pic{TableauDiagram={*[#2]{#2}(#4){#5}}}; }
    { \draw(0,0) pic{TableauDiagram={[#2]{#2}(#4){#5}}}; }
  \end{tikzpicture}
}

\NewDocumentCommand\CstDiagram{ O{1.5} m D(){0} m O{#3} }
{
  \begin{tikzpicture}[centered=-1ex]
    \draw(0,1.5) pic{TableauDiagram={*[#1]{#2}(#3){#4}}};
    \draw(0,0) pic{TableauDiagram={ [#1]{#2}(#3){#4}}};
  \end{tikzpicture}
}

\NewDocumentCommand\OmegaDiagram{ O{1.5} m D(){0} m m  }
{
  \begin{tikzpicture}[centered=-1ex]
    \draw(0,0) pic{OmegaDiagram={ [#1]{#2}(#3){#4}{#5}}};
  \end{tikzpicture}
}

\NewDocumentCommand\WebsterIdempotent{ O{1.5} m D(){0} m }
{
  \begin{tikzpicture}[centered=-1ex]
    \draw(0,0) pic{WebsterIdempotent={[#1]{#2}(#3){#4}}};
  \end{tikzpicture}
}

\title[Hecke algebras of type $G(\ell,p,n)$]
{Skew cellularity of the Hecke algebras of type $G(\ell,p,n)$}

\author{Jun Hu}\address{School of Mathematical and Statistics\\
  Beijing Institute of Technology\\
  Beijing, 100081, P.R. China}
\email{junhu404@bit.edu.cn}

\author{Andrew Mathas}
\address{School of Mathematics and Statistics, University of Sydney, NSW 2006, Australia}
\email{andrew.mathas@sydney.edu.au}

\author{Salim Rostam}
\address{Univ Rennes, CNRS, IRMAR - UMR 6625, F-35000 Rennes, France}
\email{salim.rostam@ens-rennes.fr}

\usepackage{multicol}

\begin{document}

\def\itemautorefname~#1\null{(#1)\null}
\def\subsectionautorefname{Section}
\def\sectionautorefname{Chapter}
\def\equationautorefname~#1\null{(#1)\null}

\def\Item(#1){\item\textbf{\upshape(#1\upshape)}}
\newcounter{relation}
\newcommand\Relation[1]{%
  \refstepcounter{relation}\label{#1}%
  \upshape(W$_{\arabic{relation}}$\upshape)\space%
}
\def\relationautorefname~#1\null{\upshape(W$_{#1}$\upshape)}
\let\ref\autoref
\let\eqref\autoref

\begin{abstract}
  This paper introduces (graded) skew cellular algebras, which generalise Graham and Lehrer's cellular algebras. We show that all of the main results from the theory of cellular algebras extend to skew cellular algebras and we develop a ``cellular algebra Clifford theory'' for the skew cellular algebras that arise as fixed point subalgebras of cellular algebras.

  As an application of this general theory, the main result of this paper proves that the Hecke algebras of type $G(\ell,p,n)$ are graded skew cellular algebras. In the special case when $p = 2$ this implies that the Hecke algebras of type $G(\ell,2,n)$ are graded cellular algebras. The proof of all of these results rely, in a crucial way, on the diagrammatic Cherednik algebras of Webster and Bowman.  Our main theorem extends Geck's result that the one parameter Iwahori-Hecke algebras are cellular algebras in two ways. First, our result applies to all cyclotomic Hecke algebras in the infinite series in the Shephard-Todd classification of complex reflection groups.  Secondly, we lift cellularity to the graded setting.

  As applications of our main theorem, we show that the graded decomposition matrices of the Hecke algebras of type $G(\ell,p,n)$ are unitriangular, we construct and classify their graded simple modules and we prove the existence of ``adjustment matrices'' in positive characteristic.
\end{abstract}

\maketitle
\NewDocumentCommand\notation{ mmO{}mO{}}{%
  \hyperref[#1]{$\csuse{#2}#3$} & #4 & \pageref{#1}\\
}

\section{Introduction}

  The Hecke algebras of complex reflection groups were introduced by
  Ariki and Koike \cite{AK,Ariki:Grpn} and Brou\'e and
  Malle~\cite{BM:cyc}, as generalisations of the Iwahori-Hecke algebras
  of Coxeter groups. Cyclotomic Hecke algebras have been studied
  extensively both because of their rich representation thery and
  because of their connections to reductive
  groups~\cite{Broue:conjectures}.  Interest in these algebras
  intensified with the introduction of the quiver Hecke algebras, or KLR
  algebras,  which  categorify the integrable
  highest weight representations of Kac-Moody
  algebras~\cite{KhovLaud:diagI,Rouquier:QuiverHecke2Lie}. In
  particular, Brundan--Kleshchev~\cite{BK:GradedKL} and
  Rouquier~\cite{Rouquier:QuiverHecke2Lie} proved that the Ariki-Koike
  algebras, which are the Hecke algebras associated to the complex
  reflection groups $G(\ell,1,n)$ in the classification of Shephard and
  Todd~\cite{BM:cyc,AK}, are isomorphic to the quiver Hecke
  algebras $\RG$ of type~$A$.

  The theory of cellular algebras, which was introduced by Graham and
  Lehrer~\cite{GL}, gives a framework for constructing all the
  irreducible modules of an algebra. In particular, Graham and Lehrer
  proved that the Ariki--Koike algebras are cellular. Later, Hu and
  Mathas~\cite{HuMathas:GradedCellular}, the first two named authors of this paper, extended this result to show
  that these algebras are \textit{graded} cellular algebras.

  Rostam~\cite{Rostam:Grpn}, the third named author of this paper,
  introduced the quiver Hecke algebras $\Rlpn$ of type $G(\ell,p,n)$ as
  fixed point subalgebras of $\RG$. Extending  Brundan and Kleshchev's
  graded isomorphism theorem, Rostam proved that $\Rlpn$ is isomorphic
  to a cyclotomic Hecke algebra of type $G(\ell,p,n)$. Under some strong
  assumptions on the parameters, Rostam
  used~\cite{HuMathas:GradedCellular} to prove that $\Rlpn$ is a
  cellular algebra. For example, he showed that $\Rlpn$ is a graded
  cellular algebra if~$p$ and $n$ are coprime. In general, he proved
  that a natural basis of~$\Rlpn$, arising from  a particular cellular basis
  of the Ariki--Koike algebra cannot be an ``adapted'' cellular basis; see~\cite[\textsection 5.2]{Rostam:Stuttering}.

  We can now state the main result of this paper, which shows that the
  Hecke algebras of type $G(\ell,p,n)$ are graded skew cellular
  algebras, with the important case $p = 2$.

  \begin{MainTheorem}\phantomsection\label{MainTheorem}
    Let $\Rlpn$ be the KLR algebra of type $G(\ell,p,n)$. Then
    $\Rlpn$ is a graded skew cellular algebra. Moreover, if $p = 2$ then $\Rlpn$
    is a graded cellular algebra.
  \end{MainTheorem}

  This result recovers and generalises known results from the (ungraded)
  representation theory of the Hecke algebra of $G(\ell,p,n)$ to the
  graded setting. In particular, this proves that the graded
  decomposition matrices are unitriangular, extending \cite{GenetJacon},
  and using~\cite[\textsection 10]{Bowman:ManyCellular}
  and~\cite{Kerschl:simples} we obtain a new construction and
  classification of the graded simple modules, extending
  \cite{GenetJacon,Hu:simpleGrpn}.

  As important special cases, our \hyperlink{MainTheorem}{Main Theorem}
  shows that the Iwahori-Hecke algebras of types $A_{n-1}$, $B_n=C_n$,
  $D_n$ and $I_2(n)$ are graded skew cellular algebras. These are the
  complex reflection groups of types $G(1,1,n), G(2,1,n), G(2,2,n)$ and
  $G(n,n,2)$, respectively, in the Shephard--Todd classification.  This
  result extends Geck's theorem~\cite{Geck:cellular}, which shows that
  the Iwahori-Hecke algebras of finite Coxeter groups are cellular
  algebras.  In particular, our main theorem gives a new graded cellular
  algebra structure on the Iwahori-Hecke algebras of type~$D_n$, which
  are the Hecke algebras of type $G(2,2,n)$. More generally, we show
  that the Hecke algebras of type $G(2d,2,n)$ are graded cellular
  algebras, for $d,n \geq 1$.  When $d>1$, this result is completely new, even in the ungraded
  setting. If $d = 1$ and $n \geq 2$ then we generalise Geck's result to the graded setting.  Geck's proof relies on Kazhdan--Lusztig theory, which does
  not exist for complex reflection groups. The proof of our main theorem
  relies in a crucial way on the diagram calculus introduced by
  Webster~\cite{Webster:RouquierConjecture} and
  Bowman~\cite{Bowman:ManyCellular}. In related work, LePage and
  Webster~\cite[\S4]{LePageWebster:RClpn} generalise Webster's
  diagrammatic algebras to give diagrammatic Cherednik algebras of type
  $G(\ell,p,n)$ but they do not consider questions relating to
  cellularity. Finally, note that since the first version of this paper appeared online, Lehrer and Lyu~\cite{LL} used our theory of skew cellular algebra to prove that the generalised Temperley--Lieb algebra of type $G(r,p,n)$ is graded cellular.

  To prove our main theorem, \autoref{D:SkewCellular}
  introduces \textit{(graded) skew cellular algebras}, which can be
  viewed as an analogue of Clifford theory for cellular algebras. More
  precisely, skew cellular algebras generalise the cellular algebra
  framework to certain fixed-point subalgebras of cellular algebras. We
  show that the main structural results of cellular algebras hold for
  skew cellular algebras. In particular, we show that:
  \begin{itemize}
    \item each (graded) skew cellular algebra has a family of (graded) skew cell modules
    \item the (graded) simple modules of a (graded) skew cellular
    algebra arise in a unique way as quotients of the (graded) skew cell
    modules
    \item the (graded) decomposition matrices of skew cellular algebras
    are unitriangular.
  \end{itemize}
  In contrast to cellular algebras, the simple modules of a skew cellular
  algebra are not necessarily self-dual; see \autoref{P:DualSimples} for a
  precise statement.

%
%
  The outline of this paper is as follows. \autoref{S:skew_cellularity}
  introduces and then develops the representation theory of skew
  cellular algebras, together with the closely related notion of a shift
  automorphism of a cellular algebra. \autoref{subsection:Clifford} develops Clifford theory in this setting.  \autoref{S:Hecke_algebras}
  recalls and extends the definitions and known results about the
  cyclotomic KLR algebras of type $G(\ell,p,n)$ and about the
  Webster-Bowman diagram calculus for the diagrammatic Cherednik
  algebras. \autoref{S:ShiftedRegularity} is the technical heart of the
  paper where we use the diagrammatic Cherednik algebras to define an explicit
  diagrammatic basis of $\RG$ (\autoref{D:Cslam}), which has the
  properties that we need to prove that $\Rlpn$ is a skew cellular algebra.
  \autoref{S:SkewCellularity} uses the diagram basis of $\RG$
  constructed in \autoref{S:ShiftedRegularity} to show that $\RG$ has a shift
  automorphism. Using the results from \autoref{S:skew_cellularity}, this
  implies that $\Rlpn$ is a skew cellular algebra, establishing our
  \hyperlink{MainTheorem}{Main Theorem}. Finally, as two applications of our main results,
  \autoref{subsection:adjustment} gives an ``adjustment matrix'' result for the Hecke
  algebras of type~$G(\ell,p,n)$ and \autoref{SS:Simples} gives the
  classification of the graded simple~$\Rlpn$-modules.

  An \hyperref[notation]{index of notation} can be found at the end of the paper.

%

  \subsection*{Acknowledgements}

  Jun Hu was supported by the National Natural Science Foundation of
  China (No. 12171029). Andrew Mathas was supported, in part, by the
  Australian Research Council. The authors are thankful to Chris Bowman
  for many discussions and to Lo\"ic Poulain d'Andecy for suggesting the
  term ``skew cellular''. We thank the referee for their comments and
  suggestions, which significantly improved our exposition.


\section{Skew cellular algebras}
\label{S:skew_cellularity}

  This chapter defines and then develops the representation theory of
  graded skew cellular algebras. The first section sets our
  notation for graded algebras.  The second section, which is the heart
  of the chapter, defines skew cellular algebras and shows how to extend
  the general theory of graded cellular algebras~\cite{GL,
  HuMathas:GradedCellular} to the skew setting. In the third
  section we study graded cellular algebras with shift automorphism that, like
  Clifford theory, provides a general tool for showing that fixed-point
  subalgebras of cellular algebras are skew cellular algebras. In the fourth section we study Clifford theory for the skew cellular algebras arising from a
  graded cellular algebra $A$ with a shift automorphism $\sigma$,
  especially when $\sigma_A$ is $\varepsilon$-splittable (in the sense of~\autoref{E:z_strong}).

  \subsection{Graded algebras}\label{S:Hecke}
  Throughout this paper we fix a commutative integral domain $R$ with
  one. In this paper a \textbf{graded $R$-module} is  a $\Z$-graded
  $R$-module $M=\bigoplus_{d\in\Z}M_d$. If $m\in M_d$ then $m$ is
  \textbf{homogeneous} of \textbf{degree}~$d$, for $d\in\Z$. If $M$ is a
  graded $R$-module and $s\in\Z$ let $M\<s\>$ be the graded $R$-module
  that is equal to $M$ as an (ungraded) $R$-module but where the grading
  is shifted so that the homogeneous component of $M\<s\>$ of degree
  $d$ is~$M\<s\>_d=M_{d-s}$, for $d\in\Z$.

  A \textbf{graded algebra} will always mean a $\Z$-graded algebra,
  which is a graded $R$-module such that $A_cA_d\subseteq A_{c+d}$, for
  $c,d\in\Z$. A (graded) $A$-module is a graded $R$-module
  $M=\bigoplus_{d\in\Z}M_d$  together with an $A$-action $A\times M\to M$
  such that $A_cM_d\subseteq M_{c+d}$, for $c,d\in\Z$.

  The category of graded $A$-modules has objects the graded $A$-modules
  and morphisms the homogeneous $A$-module maps of degree~$0$.  The \emph{graded dimension} of a graded
module $M=\bigoplus_{d\in\Z}M_d$ is the Laurent polynomial
\[
    \dim_t M = \sum_{d\in\Z} (\dim M_d)\,t^d\quad\in\Z[t,t^{-1}].
\] If
  the algebra $A$ comes equipped with an anti-involution then the
  \emph{dual} of a graded $A$-module~$M$ is the graded $A$-module
  \[
    M^*=\Hom_A(M,R) = \bigoplus_{d\in\Z}\Hom_A(M_d,R),
  \]
  where $R$ is in degree $0$ and $A$ acts on $M^*$ via its
  anti-involution.  A graded $A$-module $M$ is \textbf{self-dual} if
  $M\cong M^*$ as graded $A$-modules.

  Finally, if $M$ is a graded $R$-module let $\underline{M}$ be the
  ungraded $R$-module obtained by forgetting the grading. In particular,
  if $M$ is a graded $A$-module then $\underline{M}$ is an ungraded
  $\underline{A}$-module.

\subsection{Skew cellular algebras}
\label{S:skew_cellular_algebras}
  Like cellular algebras, skew cellular algebras are defined in terms of
  a skew cell datum. To describe these we first need some basic
  notation.

  Recall that a \textbf{poset}, or \textbf{partially ordered set}, is an
  ordered pair $(\P, \gedom)$, where $\P$ is a set and $\gedom$ is a
  reflexive, antisymmetric and transitive relation on~$\P$. If
  $x,y\in\P$ and $x\gedom y$ then we write $y\ledom x$. In addition, if $x\ne y$ then write
  $x\gdom y$ and $y\ldom x$.

  \begin{Definition} A \textbf{poset automorphism} of~$(\P,\gedom)$ is a permutation~$\sigma$ of $\P$ such that
    \[
      \lam\gedom\mu \text{ if and only if } \sigma(\lam)\gedom\sigma(\mu),
          \qquad\text{for all }\lam,\mu\in\P.
    \]
    If $\sigma = \iota$ is an involution we say that $\iota$ is a \emph{poset involution} of $(\P,\gedom)$.
   \end{Definition}

   Note that if $\sigma$ a poset automorphism of $(\P,\gedom)$ then $\lambda \gdom \mu \iff \sigma(\lambda) \gdom \sigma(\mu)$.
  Following \cite{GL,HuMathas:GradedCellular}, we can now define graded
  skew cellular algebras.

\begin{Definition}[$\Z$-graded skew cellular algebras]\label{D:SkewCellular}
Let $R$ be an integral domain and $A$ a $\Z$-graded $R$-algebra that is
free and of finite rank as $R$-module.

A \textbf{graded skew cell datum} for $A$ is an ordered quintuple
  $(\P,\iota, T,C,\deg)$ where $(\P,\unrhd)$ is a poset,
  $\iota$ is a poset involution of $\P$, for each
  $\lambda\in\P$ there is a finite set $T(\lambda)$ together with a bijection
  \[ \iota_\lam\map{T(\lam)} T(\iota(\lam)); \s\mapsto \iota(\s)=\iota_\lam(\s),\qquad
       \text{ for all }\s\in T(\lam),
  \]
  such that $\iota_{\iota(\lam)}\circ \iota_\lam=\id_{T(\lam)}$, and
  \[ C\map {\coprod_{\lambda\in\P}T(\lambda)\times T(\lambda)} A;
     (\s,\t)\mapsto c_{\s\t},
     \quad\text{and}\quad
     \deg: \coprod_{\lambda\in\P}T(\lambda)\rightarrow\Z
  \]
  are two functions such that $C$ is injective and
  \begin{enumerate}[label=\upshape(C$_{\arabic*}$\upshape), ref=C$_{\arabic*}$]
    \item \label{rel:cellular_homogeneous}
    Each element $c_{\s\t}$ is homogeneous of degree
    $\deg c_{\s\t}=\deg\s+\deg\t$, for $\lambda\in\P$ and
      $\s,\t\in T(\lambda)$.
    \item \label{rel:cellular_basis} The set
     $\set{c_{\s\t}|\s,\t\in T(\lambda), \lambda\in\P}$ is an $R$-basis of $A$.
    \item \label{rel:cellular_product}
    If $a\in A$ and $\s,\t\in T(\lambda)$, for $\lambda\in\P$, then
    there exist scalars $r_{\v\s}(a)$, which do not depend on $\t$, such that
    \[ a c_{\s\t}=\sum_{\v\in T(\lambda)}r_{\v\s}(a)
         c_{\v\t}\pmod {A^{\rhd\lambda}},
    \]
    where $A^{\rhd\lambda}$ is the $R$-submodule of $A$ spanned by
    $\set{c_{\a\b}|\a,\b\in T(\mu)\text{ for }\lambda\lhd\mu\in\P}$.
    \item \label{rel:cellular_star}
    There is a unique $R$-algebra anti-isomorphism $*\map AA$ such that
    $(c_{\s\t})^*=c_{\iota(\t)\iota(\s)}$, for all $\s,\t\in
    T(\lambda)$ and $\lambda\in\P$.
   \end{enumerate}
   A \textbf{$\Z$-graded skew cellular algebra} is a graded algebra
   that has a graded skew cell datum. The basis
   $\set{c_{\s\t}|\lambda\in\P\text{ and } \s,\t\in T(\lambda)}$
   is a \textbf{$\Z$-graded skew cellular basis} of~$A$.
 \end{Definition}

Applying the anti-isomorphism $\ast$ to
relation~\autoref{rel:cellular_product} and
using~\eqref{rel:cellular_star} together with the assumption that
$\iota$ is a poset involution, shows that if $a\in A$ and $\s,\t\in
T(\lambda)$, for $\lambda\in\P$, then
\[ c_{\iota(\s)\iota(\t)}a^*
    =\sum_{\u\in T(\lambda)}{r}_{\u\t}(a)c_{\iota(\s)\iota(\u)}\pmod
    {A^{\rhd\iota(\lambda)}}. \tag{C$'_3$}\label{rel:cell_star_prod}
\]
Therefore, after relabelling, if $a\in A$ and $\s,\t\in T(\lambda)$, for $\lambda\in\P$, then
\begin{equation} \label{equation:cst_astar}
    c_{\s\t}a = \sum_{\u \in T(\lam)} r_{\iota(\u),\iota(\t)}(a^*)c_{\s\u} \pmod{A^{\rhd\lam}},
\end{equation}
where the scalars $r_{\iota(\u),\iota(\t)}(a^*)$ are the same scalars
appearing in \autoref{rel:cellular_product}. In particular,
$r_{\iota(\u),\iota(\t)}(a^*)$ does not dependent on~$\s$.

\begin{Remark}\label{R:cellularAlgebras}
If $\iota=\id_\P$ is the identity map, and
$\iota_\lambda=\id_{T(\lambda)}$ for all~$\lambda\in\P$, then
\autoref{D:SkewCellular} recovers the definition of \textit{graded
cellular algebras} from \cite{HuMathas:GradedCellular}. If, in addition,
we forget the grading on~$A$ then \autoref{D:SkewCellular} reduces to
Graham and Lehrer's original definition of cellular algebras~\cite{GL}.
Thus, graded cellular algebras are given by a graded cellular datum
$(\P,T,C,\deg)$ and cellular algebras are given by a cell datum
$(\P,T,C)$. A \textbf{skew cellular algebra} is a graded skew cellular
algebra that is concentrated in degree $0$.  In particular, skew
cellular algebras are a generalisation of Graham and Lehrer's definition
of cellular algebras.
\end{Remark}

The reader might find it helpful to refer to the following example when
reading this section. More complicated examples of skew cellular algebras
are given in \autoref{Ex:SkewCellular} below.

\begin{Example}
  We give a ``toy example''. Let $R$ be any ring and let $x$ and $y$ be
  indeterminates over $R$. Fix an integer $m\ge1$ and set
  $A=R[x]/(x^m)\oplus R[y]/(y^m)$. Let $(\P,\gedom)$ be the poset
  $\P=\Z_2\times\set{0,1,\dots,m-1}$ with $(i,k)\gedom(i',k')$ only if $i=i'$
  and $k\ge k'$ (as integers). Define the poset involution
  $\iota\map\P\P$ by $\iota(i,k)=(i+1,k)$. For $\lambda=(i,k)\in\P$ set
  $T(\lambda)=\set{k}$ and $\deg(k) = k$. In particular, $\iota_\lambda(k)=k$.
   Then,  for $\lambda=(i,k)\in\P$ we have  $\s,\t\in T(\lambda)$ only if $\s=\t=k$, so
  define
  \[
  C_{\s\t} = C_{kk} = \begin{cases*}
             x^k,& if $i=0$,\\
             y^k,& if $i=1$.\\
    \end{cases*}
  \]
  Then $(\P,\iota,T,C,\deg)$ is a $\Z$-graded skew cell datum for $A$.
\end{Example}

For the rest of this section fix a graded skew cellular algebra $A$
with skew cell datum $(\P,\iota,T,C,\deg)$. We now study the graded
representation theory of $A$, generalising the results
of~\cite{GL,HuMathas:GradedCellular}.

\begin{Definition}\label{D:CellModule} Let $\lam\in\P$. The \textbf{(left) graded
  skew cell module} $C_\lambda$ is the left graded $A$-module with basis
  $\set{c_\s|\s\in T(\lam)}$ and with $A$-action determined by
\[
  ac_\s=\sum_{\t\in T(\lam)}r_{\t\s}(a)c_\t,\qquad
      \text{for $a\in A$ and $\s\in T(\lam)$,}
\]
where $r_{\t\s}(a)\in R$ is the scalar defined
in~\autoref{rel:cellular_product}.
\end{Definition}

\begin{Remark}
The name ``skew cell module'' also appears
in~\cite{BowmanEnyangGoodman:Diagram} but where the term ``skew''
 refers to skew Young diagrams. A priori, these are different objects.
\end{Remark}

Let $\lambda\in\P$. Then $C_\lambda=\bigoplus_{d\in\Z}(C_\lambda)_d$ is a graded $A$-module,
where $(C_\lambda)_d$ is the free $R$-module with basis
$\set{ c_\s | \s \in T(\lam)\text{ with } \deg \s = d}$.

If $\u,\s,\t,\v\in T(\lam)$ then, by \autoref{rel:cellular_product} and
\autoref{equation:cst_astar},
\[
c_{\u\s}c_{\t\v} \equiv\sum_{\a\in T(\lam)}r_{\a\t}(c_{\u\s})c_{\a\v}
    \equiv \sum_{\b\in T(\lam)}r_{\iota(\b)\iota(\s)}(c_{\t\v}^*)c_{\u\b}\pmod{A^{\rhd\lambda}}.
\]
It follows that
\begin{equation}
\label{E:r_at(c_us)}
r_{\a\t}(c_{\u\s})\neq 0 \text{ only if }\a=\u,
\end{equation}
and $r_{\iota(\b)\iota(\s)}(c_{\t\v}^*)\neq 0$ only if $\b=\v$ and,
consequently, that
\[
    r_{\u\t}(c_{\u\s})=r_{\iota(\v)\iota(\s)}(c_{\t\v}^*).
\]
By \autoref{rel:cellular_product} and \autoref{equation:cst_astar}, the
scalar $r_{\u\t}(c_{\u\t})$ depends only on $\s$ and~$\t$ and
\textit{not} on the choice of~$\u$ and~$\v$.

The next definition is motivated by \cite[Definition~2.3]{GL} and the
calculations above.
\begin{Definition}\label{D:cell_Form} Let $\lambda\in\P$. Let
  $\phi=\phi_\lambda\map{C_\lam\times C_\lam}R$ be the $R$-bilinear map
  determined by
  \[
  \phi(c_\s,c_{\t})\coloneqq r_{\u\t}(c_{\u\s})
            =r_{\iota(\v)\iota(\s)}(c_{\t\v}^*)\in R,
            \qquad\text{for all } \s,\t\in T(\lam).
  \]
\end{Definition}
Then, by the calculations above,
\begin{equation}\label{E:form}
   c_{\u\s}c_{\t\v} \equiv \phi_\lambda(c_\s,c_\t)c_{\u\v}\pmod{A^{\rhd\lambda}}.
\end{equation}


To better understand $\phi_\lambda$ we abuse notation and extend the map
$\iota_\lambda\map{T(\lambda)}{T(\iota(\lambda))}$ to an $R$-linear isomorphism
\[
    \iota_\lambda\map{C_\lambda}C_{\iota(\lambda)}; c_\s\mapsto c_{\iota(\s)},
        \qquad\text{for }s\in T(\lambda).
\]
In general, the $R$-linear map
$\iota_\lambda\map{C_\lambda}C_{\iota(\lambda)}$ is not an $A$-module
homomorphism. If $\lambda\in\P$ is fixed and $x\in C_\lambda$ then we
simplify our notation and write $\iota(x)=\iota_\lambda(x)\in
C_{\iota(\lambda)}$. In particular, $\iota(c_\s)=\iota_\lambda(c_\s)$,
for $\s\in T(\lambda)$.

The following lemma, which gives the main properties of $\phi_\lambda$,
is modelled on \cite[Proposition~2.9]{Mathas:ULect}.  However, note that
for skew cellular algebras the bilinear form $\phi_\lam$ \textit{not}
necessarily symmetric.

\begin{Lemma}\label{lemma:phi_symm} Let $\lam\in\P$ and $x, y \in C_\lam$.
  \begin{enumerate}
    \item \label{item:skew_sym}
    We have $ \phi_\lam(x, y)=\phi_{\iota(\lambda)}\bigl(\iota_\lam(y),\iota_\lam(x)\bigr)$.
    \item \label{item:rad_submodule}
      If $a\in A$ then
      $\phi_\lam(x,ay)=\phi_{\iota(\lam)}\bigl(\iota_{\lam}(y),a^*\iota_\lam(x)\bigr)
                 =\phi_\lam\bigl(\iota_{\iota(\lam)}(a^*\iota_\lam(x)),y\bigr)$.
    \item \label{item:cuv_x}
      If $\s,\u \in T(\lam)$ then $c_{\u\s} x = \phi_\lam(c_\s, x) c_\u$.
    \item \label{item:homogeneous_form}
    The form $\phi_\lambda$ is homogeneous of degree~$0$.
\end{enumerate}
\end{Lemma}

\begin{proof}
Since $\phi_\lam$ is bilinear, and $\iota_\lambda$ is linear, it suffices to consider the cases when
$x = c_\s$ and $y = c_\t$ for $\s,\t \in T(\lam)$.

By the discussion above \autoref{D:cell_Form}, for any $\u,\v \in T(\lam)$ we have
\[
    \phi_\lam(c_\s,c_\t) = r_{\u\t}(c_{\u\s})
                    = r_{\iota(\v),\iota(\s)}(c_{\t\v}^*)
                    = r_{\iota(\v),\iota(\s)}(c_{\iota(\v),\iota(\t)})
                    = \phi_{\iota(\lambda)}(c_{\iota(\t)},c_{\iota(\s)}),
\]
proving~\autoref{item:skew_sym}.
Now consider part~\autoref{item:rad_submodule}. Working modulo $A^{\rhd\lam}$, for any $\u,\v \in T(\lam)$
\begin{align*}
  \phi_\lambda(c_\s, a c_\t)c_{\u\v}
    &= \sum_{\t' \in T(\lam)} r_{\t'\t}(a) \phi_\lam(c_\s,c_{\t'}) c_{\u\v}
        & \text{by \autoref{D:CellModule}}, \\
    &\equiv \sum_{\t' \in T(\lam)} r_{\t'\t}(a) c_{\u\s} c_{\t'\v}
        & \text{by \autoref{E:form}}, \\
    &\equiv c_{\u\s} (a c_{\t\v})
        & \text{by \autoref{rel:cellular_product}}, \\
    &= (c_{\u\s}a)c_{\t\v} \\
    &\equiv \sum_{\s' \in T(\lam)} r_{\iota(\s'),\iota(\s)}(a^*) c_{\u\s'} c_{\t\v}
        & \text{by \autoref{rel:cell_star_prod}}, \\
    &\equiv \sum_{\s' \in T(\lam)} r_{\iota(\s'),\iota(\s)}(a^*) \phi_\lam(c_{\s'},c_\t)c_{\u\v}
        & \text{by \autoref{E:form}}, \\
    &\equiv \sum_{\s' \in T(\lam)} r_{\iota(\s'),\iota(\s)}(a^*)
        \phi_{\iota(\lambda)}(c_{\iota(\t)},c_{\iota(\s')})c_{\u\v}
        & \text{by part~(a)}, \\
    &= \phi_{\iota(\lambda)}(c_{\iota(\t)},a^* c_{\iota(\s)})c_{\u\v} \pmod{A^{\rhd\lam}}
        & \text{by \autoref{D:CellModule}}. \\
\end{align*}
Therefore, $\phi_\lambda(c_{\s},ac_{\t}) = \phi_{\iota(\lambda)}(c_{\iota(\t)},a^*c_{\iota(\s)})$.
Hence, $\phi_\lambda(x,ay)=\phi_{\iota(\lambda)}(\iota(y), a^*\iota(x))$,
proving the first identity in~\autoref{item:rad_submodule}. Applying part~\autoref{item:skew_sym}, we deduce that
\[
    \phi_\lambda(x,ay)=\phi_{\iota(\lambda)}(\iota_\lam(y), a^*\iota_\lam(x))
                      =\phi_\lambda\bigl(\iota_{\iota(\lam)}(a^*\iota_\lam(x)),y\bigr),
\]
completing the proof of~\autoref{item:rad_submodule}.

Finally, part~\autoref{item:cuv_x} follows immediately from \autoref{E:form} and part~\autoref{item:homogeneous_form}
follows by comparing the degrees on the left and right hand sides of
\autoref{E:form}.
\end{proof}


\begin{Remark}
  If $\s\in T(\lambda)$ and $a\in A$ then
$\iota_{\iota(\lambda)}\bigl(a\iota_\lam(x)\bigr)=ax$ if and only if
  $\iota_\lambda\map{C_\lambda}C_{\iota(\lambda)}$ is an $A$-module homomorphism.
  Hence, by \autoref{lemma:phi_symm}\autoref{item:rad_submodule}, the bilinear form
  $\phi_\lambda$ is associative if and only if $\iota_\lambda$ is an
  $A$-module homomorphism. In particular, $\phi_\lambda$ is symmetric and associative
  when $\iota_\lambda=\id_{C_\lambda}$, as is the case for (non-skew)
  cellular algebras.
\end{Remark}

For any $\lambda\in\P$, and for any ring $R$, the radical of $C_\lam$ is:
\begin{equation}\label{E:CellRadical}
  \rad(C_\lam) \coloneqq\set{y\in C_\lam |\phi_\lam(x, y)=0 \text{ for all }x \in C_\lam} .
\end{equation}

The next proposition is the skew cellular algebra analogue
of~\cite[Lemma 2.7]{HuMathas:GradedCellular}.


\begin{Proposition}\label{P:radical}
Let $\lambda\in\P$. Then the radical $\rad(C_\lam)$ is a graded $A$-submodule of $C_\lam$.
\end{Proposition}

\begin{proof}
 If $y\in\rad(C_\lam)$ and $a\in A$ then, by~\autoref{lemma:phi_symm}\autoref{item:rad_submodule},
 \[\phi(x, ay) = \phi\bigl(\iota_{\iota(\lam)}(a^*\iota_\lam(x)),y\bigr) = 0,\qquad
                \text{for any $x \in C_\lam$}.
 \]
 Therefore, $ay\in\rad(C_\lambda)$, showing that $\rad(C_\lambda)$ is an
 $A$-submodule of $C_\lambda$. By \autoref{lemma:phi_symm}(d), the form
 $\phi_\lambda$ is homogeneous of degree~$0$, so $\rad(C_\lambda)$ is a
 graded submodule of $C_\lambda$.
\end{proof}

\begin{Remark}
Note that because $\phi_\lambda$ is not symmetric this is the
\textit{right radical} of $\phi_\lambda$ and, \textit{a priori}, this is
different from the \textit{left radical} of $\phi_\lambda$.
It is not clear if the left radical of $\phi_\lambda$ is an $A$-submodule of $C_\lambda$ because there is no obvious left-handed analogue of \autoref{lemma:phi_symm}\autoref{item:rad_submodule}.
\end{Remark}


\begin{Definition}\label{D:CellularSimples}
  For $\lambda\in\P$ set $D_\lambda = C_\lambda/\rad(C_\lambda)$. Let
  $\P_0=\set{\lambda\in\P|D_\lambda\ne0}$.
\end{Definition}

By \autoref{P:radical}, $D_\lambda$ is a graded $A$-module. The next result extends the arguments
of~\cite[\textsection 2.2]{HuMathas:GradedCellular}, to
characterise the graded simple $A$-modules. Recall that the \emph{Jacobson radical} of an $A$-module $M$ is the intersection of its  maximal $A$-submodules.


\begin{Theorem}\label{T:SkewCellularSimples}
Suppose that $R$ is a field.
\begin{enumerate}
\item If $\lam \in \P_0$, $W$ is an $A$-submodule of $\underline{C}_\lam$ and $\theta\in\Hom_A(\underline{C}_\lam,\underline{C}_\lam/W)$, then there exists a unit $r\in R^\times$ such that $\theta(x)=rx+W$ for all $x\in \underline{C}_\lam$.
\item
\label{I:Dlam_irred_unique_simple_head}
If $\lam \in \P_0$ then $D_\lam$ is an absolutely irreducible graded $A$-module. Moreover,  $\rad C_\lambda$ is the Jacobson radical of $C_\lambda$ and, consequently, $D_\lambda$ is the unique simple head of $C_\lam$.
\item If $\lambda, \mu \in \P_0$ and $D_\lam \simeq D_\mu\langle k\rangle$, for some $k \in \mathbb{Z}$, then $\lambda = \mu$ and $k = 0$.
\item The set $\set{ D_\lam\langle k\rangle | \lam \in \P_0\text{ and }k \in \Z}$
is a complete set of pairwise non-isomorphic graded simple $A$-modules.
\end{enumerate}
\end{Theorem}


\begin{proof} Part (a) follows the same argument from the proof of
  \cite[Proposition (2.6)]{GL}. Similarly, the proof
  of~\cite[Proposition (3.2)(ii)]{GL} and~\cite[\textsection
  2.2]{HuMathas:GradedCellular} show that if $\lam \in \P_0$
  then $D_\lam$ is an absolutely irreducible graded $A$-module. Now
  using~(a) we deduce that $D_\lam$ is the unique simple head of
  $C_\lam$ and, hence, that $\rad(C_\lam)$ is the unique maximal submodule of
  $C_\lam$, so $\rad C_\lambda$ is the Jacobson radical of $C_\lambda$. The remaining parts of the theorem
  follows using the arguments
  from~\cite[\textsection 2.2]{HuMathas:GradedCellular}.
\end{proof}

\begin{Corollary}
Suppose that $R$ is a field. Then $\set{\underline{D}_\lam | \lam \in
\P_0}$ is a complete set of pairwise non-isomorphic ungraded simple
$A$-modules.
\end{Corollary}

%


The next result describes the duals of the simple modules of skew
cellular algebras.

\begin{Proposition}\label{P:DualSimples}
  Suppose that $R$ is a field. Then $D_\lambda^* \simeq D_{\iota(\lambda)}$ as graded $A$-modules.
In particular, $\lambda \in \P_0$ if and only if $\iota(\lambda) \in \P_0$.
\end{Proposition}


\begin{proof}
Let  $\lambda \in \P_0$. Then, by definition, there exist $x,y \in C_\lambda$ such that $\phi_\lambda(x,y) \neq 0$. By~\autoref{lemma:phi_symm}\autoref{item:skew_sym}, we have $\phi_{\iota(\lambda)}(\iota_\lambda(y),\iota_\lambda(x))=\phi_\lambda(x,y)\neq 0$ and thus $\iota(\lambda) \in \P_0$.

  Now for each $y\in C_{\iota(\lambda)}$
  there is a well-defined $R$-linear map from $\theta_y\map{C_\lambda}R$ given by
  \[
    \theta_y(x) = \phi_{\lambda}(\iota_{\iota(\lambda)}(y), x),\qquad
    \text{for }x\in C_{\lambda}.
  \]
  By definition, $\theta_y(x)=0$ if $x\in\rad C_{\lambda}$, so we
  can consider $\theta_y$ as a map from  $D_{\lambda}$ to $R$.
  Again by~\autoref{lemma:phi_symm}\autoref{item:skew_sym},
  $
    \phi_{\lambda}(\iota_{\iota(\lambda)}(y), x) = \phi_{\iota(\lambda)}(\iota_{\lambda}(x),y),
  $
  so $\theta_y=0$ if $y\in\rad C_{\iota(\lambda)}$. Hence, there is a well
  defined map $\Theta_\lambda\map{D_{\iota(\lambda)}}{D_\lambda^*}$ given by
  $\Theta_\lambda(y+\rad C_{\iota(\lambda)}) = \theta_y$, for $y\in C_{\iota(\lambda)}$.
  The map~$\Theta_\lambda$ is homogeneous of degree $0$ since
  $\phi_\lambda$ is homogeneous of degree~$0$ by
  \autoref{lemma:phi_symm}\autoref{item:homogeneous_form}.
  Moreover, by \autoref{lemma:phi_symm}\autoref{item:rad_submodule}, if $a\in A$ then
  \[  \theta_{ay}(x) = \phi_\lambda\bigl(\iota_{\iota(\lambda)}(ay),x\bigr)
                     =\phi_\lambda(\iota_{\iota(\lambda)}(y),a^*x)
                     =\theta_y(a^*x)=a\theta_y(x).
  \]
  So $\Theta_\lambda$ is a morphism of $A$-modules. Moreover, since $\iota_{\iota(\lam)}$ is a bijection, it follows from the equality $\phi_\lambda(\iota_{\iota(\lambda)}(y), x) = \phi_{\iota(\lambda)}(\iota_\lambda(x),y)$ that $\Theta_\lambda\map{D_{\iota(\lambda)}}{D_\lambda^*}$ is injective.

 We now assume that $R$ is a field. Since $\Theta_\lambda$ is injective we deduce that $\dim_R D_{\iota(\lambda)} \leq \dim_R D_\lambda$. By the same argument, the map $\Theta_{\iota(\lambda)} : D_\lambda \to D_{\iota(\lambda)}^*$ is also injective, which gives the reverse inequality. We deduce that $\Theta_\lambda$ is an isomorphism and this concludes the proof.
%
 \end{proof}
%
%
%


%

By \autoref{P:DualSimples}, if $R$ is a field and $\lambda \in \P_0$
then the graded $A$-module $D_\lambda$ is self-dual if and only if $\lambda =
\iota(\lambda)$. In the special case when $A$ is a cellular algebra,
this recovers the well-known result that the simple modules of cellular
algebras are self-dual since the involution $\iota$ is the identity map
in this case by \autoref{R:cellularAlgebras}.


Finally, as in~\cite[\textsection 2.3]{HuMathas:GradedCellular}, define the
\textbf{graded decomposition matrix} of~$A$ to be the matrix
\begin{equation}\label{E:DecompNum}
D_A(t) = \bigl(d_{\lambda\mu}(t)\bigr),\quad\text{where }\quad
d_{\lambda\mu}(t) \coloneqq \sum_{k \in \mathbb{Z}}
    [C_\lambda : D_\mu\langle k\rangle] t^k,
    \text{ for $\lambda\in\P$ and $\mu\in\P_0$.}
\end{equation}
We order the rows and columns of $D_A(t)$ by any total order $\ge$ that
extends $\gedom$; that is, if~$\lambda\gedom\mu$ then $\lambda\ge\mu$, for
$\lambda,\mu\in\P$. The arguments for cellular algebras now generalise
to prove the following.

\begin{Proposition}\label{P:unitriangular}
Let $\lambda \in \P$ and $\mu \in \P_0$. Then:
\begin{enumerate}
\item $d_{\lambda\mu}(t) \in \mathbb{Z}_{\geq 0}[t,t^{-1}]$;
\item $d_{\lambda\mu}(1) = [\underline{C}_\lambda : \underline{D}_\mu]$;
\item $d_{\mu\mu}(t) = 1$ and $d_{\lambda\mu}(t) \neq 0$ only if $\lambda \trianglerighteq \mu$.
\end{enumerate}
In particular, the graded decomposition matrix $D_A(t)$ is upper unitriangular.
\end{Proposition}

\subsection{Shift automorphisms of graded cellular algebras}

This section defines \textit{shift automorphisms} of graded cellular
algebras, which provides a general framework for constructing skew
cellular algebras from cellular algebras. This framework is used to
prove all of the main results in this paper.

As in the last section, let $R$ be an integral domain with one. Recall
from \autoref{R:cellularAlgebras} that a graded cellular algebra~$A$ is
determined by a graded cell datum $(\P,T,C,\deg)$.

\begin{Definition}\label{D:sigma_cellular}
Let $A$ be a $\Z$-graded cellular $R$-algebra with graded cell datum $(\P,T,C,\deg)$. A
\textbf{shift automorphism} of $A$ is a triple of automorphisms
$\sigma=(\sigma_A, \sigma_\P, \sigma_T)$ where~$\sigma_A$ is an $R$-algebra automorphism
of $A$, $\sigma_\P$ is a poset automorphism of~$\P$ and $\sigma_T$ is an automorphism of
the set $T = \coprod_\lambda T(\lambda)$ such that:
\begin{enumerate}
  \item If $\s\in T(\lambda)$ then $\sigma_T(\s)\in T(\sigma_P(\lambda))$
    and $\deg(\sigma_T(\s))=\deg(\s)$.
  \item If $\s,\t\in T(\lambda)$ then $\sigma_A(c_{\s\t})=c_{\sigma_T(s)\sigma_T(\t)}$.
  \item If $\s,\t\in T(\lambda)$, for $\lambda\in\P$, then
  $\sigma_T^k(\t)=\t$ if and only if $\sigma_T^k(\s)=\s$, for $k\in\Z$.
\end{enumerate}
\end{Definition}


Throughout this section we fix a graded cellular algebra $A$ with a
shift automorphism~$\sigma$.  The algebra of \emph{$\sigma_A$-fixed
points} in $A$ is
\[
A^\sigma = \set{ a \in A | \sigma_A(a) = a}.
\]
The aim of this section is to show that $A^\sigma$ is a skew cellular
algebra.

In practice, a shift automorphism
$\sigma=(\sigma_A,\sigma_P,\sigma_T)$ is completely determined by the
map~$\sigma_T$. Part~(c) of \autoref{D:sigma_cellular}, which is used in~\autoref{E:ckst}, is a non-trivial assumption that ensures that whenever
$\sigma_T^k$ restricts to give an automorphism of $T(\lambda)$ then all
of the $\sigma_T^k$-orbits in~$T(\lambda)$ have the same size.  As the
meaning will always be clear from context, we often abuse notation and
simply write $\sigma$ instead of $\sigma_A$, $\sigma_\P$ or $\sigma_T$.

%

A trivial example of a shift automorphism is given by taking the
identity maps $(\id_A, \id_\P, \id_T)$. Here is a less trivial example.

\begin{Example} \label{example:sigma_cellular}
For $n \in \mathbb{Z}_{\geq 1}$, the full matrix algebra $A =
\mathrm{Mat}_n(R)$ is a cellular algebra, with $\P = \set{n}$ and
$T(n) = \set{1, \dots, n}$, where the cellular basis elements is given
by the set of elementary matrices $E_{ij}$, where $1\leq i,j\leq n$.
Following \cite[Example 2.1.3]{Mathas:Singapore}, fix integers $d_1,
\dots, d_n \in \mathbb{Z}$ such that $d_i + d_{n+1-i} = 0$ and define a
$\mathbb{Z}$-grading on $A$ by setting $\deg E_{ij} = d_i +d_{n+1-j}$
for $1\leq i,j\leq n$. Then $A$ is a graded cellular algebra with graded
cellular basis $\set{F_{ij}}$, where  $F_{ij} = E_{i(n+1-j)}$ and
$\deg(i)=d_i$, for $1\le i,j\le n$. The condition $d_j+d_{n+1-j}=0$
ensures that the degrees add in the relation
$F_{ij}F_{(n+1-j)k}=F_{ik}$.

To define a shift automorphism of $A$, suppose that $w\in\Sym_n$
is any permutation that is a product of $p$ disjoint $\frac np$-cycles
such that $w(n+1-i)=n+1-w(i)$, for $1\le i\le n$.  For example, $w$
could be the permutation given by $w(i)=n+1-i$, for $1\le i\le n$, when $n$ is even and $p = \frac{n}{2}$.
Then~$A$ has a unique shift automorphism
$\sigma=(\sigma_A,\sigma_P,\sigma_T)$ such that $\sigma_T(i)=w(i)$, for
$1\le i\le n$.  Explicitly, $\sigma_T=w$, $\sigma_\P=\id_\P$ and
$\sigma_A(F_{ij})=F_{w(i)w(j)}$, for $1\le i,j\le n$. The assumption
that $w$ is the product of $p$ disjoint $\frac np$-cycles ensures that
all of the $\sigma_T$-orbits have size~$p$, in accordance with
\autoref{D:sigma_cellular}(c). The second condition on $w$
is forced by the requirement that~$\sigma_A$ respect the relations
$F_{ij}F_{(n+1-j)k}=F_{ik}$. If $n$ is odd then by looking at the
sum $\tfrac{n(n+1)}2=\sum_{i=1}^nw(i)$, it follows that
$w(\tfrac{n+1}2)=\tfrac{n+1}2$, which forces $w=1_{\Sym_n}$. So,
$w\ne1$ only if $n$ is even.

If we drop the requirement that $A$ be a \textit{graded} algebra then
we do not need to assume that $w(i)=n+1-i$, for $1\le i\le n$.
\end{Example}

The following properties of shift automorphisms are immediate from
\autoref{D:sigma_cellular}.

\begin{Lemma}\label{L:SkewAutoProperties}
  Suppose that $(\sigma_A,\sigma_\P, \sigma_T)$ is a shift automorphism
  of the graded cellular algebra~$A$. Then $\sigma_A$ is homogeneous
  automorphism of~$A$ of degree zero such that
  $\sigma_A(A^{\gdom\lambda})=A^{\gdom\sigma_\P(\lambda)}$. Moreover,
  $\sigma_A(a^*)=(\sigma_A(a))^*$, for all $a\in A$.
\end{Lemma}

For the rest of this section fix a graded cellular algebra $A$ with
graded cell datum $(\P,T,C,\deg)$ and a skew cellular algebra
automorphism $\sigma=(\sigma_A,\sigma_\P, \sigma_T)$.
 Let
\[
        \sp =|\sigma_A|
        \qquad\text{and}\qquad
        \sp_\P=|\sigma_\P|
\]
be the orders of the automorphisms $\sigma_A$ and  $\sigma_\P$,
respectively. Note that  $\sp$ is also the order of~$\sigma_T$, since by~\autoref{D:sigma_cellular}b) if $k \in \Z$ then
\begin{align*}
\sigma_A^k = \id_A
&\iff
\sigma_A^k(c_{\s\t}) = c_{\s\t} && \text{for all } \s,\t \in T(\lambda)\text{ and } \lambda \in \P
\\
&\iff
c_{\sigma_T^k(\s),\sigma_T^k(\t)} = c_{\s\t} && \text{for all } \s,\t \in T(\lambda)\text{ and }\lambda \in \P
\\
&\iff
\sigma_T^k (\s) = \s && \text{for all } \s \in T(\lambda)\text{ and }\lambda \in \P
\\
&\iff
\sigma_T^k = \id_T.
\end{align*}
In particular, both $\sp_\P$ and $\sp$ are finite since $\P$ and $T$ are finite sets.
Finally, by~\autoref{D:sigma_cellular}a) we have that $\sp_\P$ divides $\sp$.
For the rest of this section we assume that $R$ contains a primitive
$\sp$th root of unity~$\eps$ and $\sp\cdot 1_R$ is invertible in $R$.


The cyclic group $\Z_\P=\<\sigma_\P\>\cong\Z/\sp_\P\Z$ acts on $\P$, let $\Prep$ be a set of representatives for this action. For example, if $\leq$ is any total order refining $\ledom$ then one could take
\[
    \Prep = \set{\lambda\in\P|\sigma_\P^k(\lambda)\le\lambda \text{ for }0\le k<\sp_\P}
\]
For each
$\lambda\in\P$ let
\label{E:so_sp}
\[
        \so_\lambda = \min\set{k\ge1|\sigma_\P^k(\lambda)=\lambda}
        \qquad \text{and}\qquad
        \plam = \sp / \so_\lambda.
\]
Then $\so_\lambda$ is the size of the $\Z_\P$-orbit of $\lambda$, so
$\so_\lambda$ divides $\sp_\P$ and $\plam\in\N$.


\begin{Lemma}
\label{L:no_chains_in_orbit}
The elements in the same $\Z_\P$-orbit are not comparable under
$\ledom$. That is, if $\lambda \in \P$ and $k \in
\{1,\dots,\so_\lambda-1\}$, the elements $\lambda$ and $\sigma^k\lambda$
are not comparable under $\ledom$.
\end{Lemma}

\begin{proof}
  Let $\lambda \in \P$ and let
  $[\lambda]=\set{\sigma_\P^k(\lambda)|0\le k<\so_\lambda}$ be the
$\Z_\P$-orbit of $\lambda$. By way of
contradiction, suppose that there exist $k, l \in \Z$ such
that $\sigma_\P^k\lambda \ldom \sigma_\P^l \lambda$. Since $\sigma_\P$ is a poset
automorphism, for any $m \in \Z$
\[
\sigma_\P^m \lambda = \sigma_\P^{m-k}\left(\sigma_\P^k\lambda\right) \ldom \sigma_\P^{m-k}\left(\sigma_\P^l\lambda\right) = \sigma_\P^{m-k+l}\lambda.
\]
We have shown that for any $\mu \in [\lambda]$, there exists $\nu \in
[\lambda]$ such that $\mu \ldom \nu$. But this is absurd because this
implies that the finite poset $([\lambda],\ledom)$ has no maximal
element.
\end{proof}

Define the binary relations $\lesig$ and $\leqsig$ on $\Prep$ by
\begin{align*}
\lambda \lesig \mu &\iff \text{there exists } k \in \Z, \sigma_\P^k \lambda \ldom \mu,
\\
\intertext{and}
\lambda \leqsig \mu &\iff \lambda = \mu \text{ or } \lambda \lesig \mu,
\end{align*}
for any $\lambda,\mu \in \Prep$. Since $\sigma_\P$ is a poset automorphism,
\begin{align*}
\lambda \lesig \mu
&\iff \text{there exists } l \text{ such that } \lambda \ldom \sigma_\P^l \mu,
\\
&\iff \text{there exist } k,l \text{ such that } \sigma_\P^k \lambda \ldom \sigma_\P^l\mu
\\
&\iff \text{for all } l, \text{ there exists } k \text{ such that } \sigma_\P^k \lambda \ldom \sigma_\P^l\mu
\\
&\iff \text{for all } k, \text{ there exists } l \text{ such that } \sigma_\P^k \lambda \ldom \sigma_\P^l\mu.
\end{align*}

\begin{Proposition}
The binary relation $\leqsig$ is a partial order on $\Prep$.
\end{Proposition}

\begin{proof}
By definition, $\lambda \leqsig \lambda$ for all $\lambda \in
\Prep$. To show that $\leqsig$ is transitive suppose that
$\lambda \leqsig \mu \leqsig \nu$, for $\lambda,\mu,\nu \in \Prep$.
If either $\lambda = \mu$ or $\mu =\nu$ then
$\lambda \leqsig \nu$, so we can assume that $\lambda \neq \mu \neq \nu$.
Then there exist $k,l$ such that $\sigma_\P^k \lambda \ldom  \mu
\ldom \sigma_\P^l \nu$  thus $\sigma_\P^k \lambda \ldom \sigma_\P^l \nu$  (since
$\ldom$ is transitive) thus $\lambda \lesig \nu$. Finally, if $\lambda
\leqsig \mu \leqsig \lambda$ then if $\lambda \neq \mu$  there exist
$k,l$ such that $\lambda \ldom \sigma_\P^k\mu
\ldom\sigma_\P^l\lambda$. Using again the transitivity of $\ldom$ we obtain
$\lambda \ldom \sigma_\P^l\lambda$ which
contradicts~\autoref{L:no_chains_in_orbit}.
\end{proof}


 Let
$\sigma_\lambda=\sigma_T^{\so_\lambda}$. Then the cyclic group
$\Z_\lambda=\<\sigma_\lambda\>\cong\Z/\plam\Z$ acts on
$T(\lambda)$. Let $\Tsig(\lambda)$ be any set of representatives for
the $\Z_\lambda$-orbits of~$T(\lambda)$. By
\autoref{D:sigma_cellular}(c), all of the $\Z_\lambda$-orbits in
$T(\lambda)$ have the same size, so $|\Tsig(\lambda)|$ divides
$|T(\lambda)|$. Let $\olambdatab=|T(\lambda)|/|\Tsig(\lambda)|$ be
the size of any $\Z_\lambda$-orbit in $T(\lambda)$. Note that $\olambdatab$ divides $\lvert \Z_\lambda\rvert = \sp_\lambda$, in particular $\olambdatab$ divides $\sp$ since $\sp_\lambda$ divides $\sp$. If~$\lambda$
and~$\mu$ are in the same $\Z_\P$-orbit then it is easy to see that
$\so_\lambda=\so_\mu$ and $\olambdatab=\olambdatab[\mu]$.

Let $\Psig=\set{(\lambda,k)|\lambda\in\Prep\text{ and } k \in \Z/\olambdatab\Z}$, considered as a poset with ordering $\leqsig$ given
by
\[
(\lambda,k) \leqsig (\mu,l) \iff (\lambda,k) = (\mu,l) \text{ or } \lambda \lesig \mu,
\]
for all $(\lambda,k),(\mu,l) \in \Psig$. We write $(\lambda,k) \lesig (\mu,l)$ if $(\lambda,k) \neq (\mu,l)$ and $(\lambda,k) \leqsig (\mu,l)$, that is, if $\lambda \lesig \mu$.
For $(\lambda,k)\in\Psig$, define $\Tsig(\lambda,k) = \Tsig(\lambda)$. Finally, set
\begin{equation}\label{E:ckst}
  c^{(k)}_{\s\t} = \sum_{j=0}^{\olambdatab-1}
  \epsl^{kj}\sigmavg(c_{\s,\sigma_\lambda^j(\t)}),
  \qquad\text{for }\s,\t\in\Tsig(\lambda,k),
\end{equation}
where $\epsl \coloneqq \eps^{\sp/\olambdatab}$ and $\sigmavg\coloneqq\sum_{l=0}^{\sp-1}\sigma_A^l$. To complete the
definition of a skew cell datum for~$A^\sigma$, let $\iota_\sigma$ be the
poset involution of $\Psig$ given by
$\iota_\sigma(\lambda,k)=(\lambda,-k)$ and let
$(\iota_\sigma)_{(\lambda,k)} : \Tsig(\lambda,k) \to \Tsig(\lambda,-k)$ be given by the identity map of $\Tsig(\lambda)$, for
$(\lambda,k)\in\Psig$. Finally, if
$\s,\t\in \Tsig(\lambda,k)$ set $C_\sigma(\s,\t)=c^{(k)}_{\s\t}$
and $\deg_\sigma(\s)=\deg(\s)$.

We can now show that a (graded) cellular algebra with a shift
automorphism gives rise to a (graded) skew cellular algebra in the
sense of \autoref{D:SkewCellular}. This result can be viewed as a
cellular algebra analogue of Clifford theory.

\begin{Theorem}
\label{proposition:sigma_cellular_implies_extended}
Suppose that  $A$ is a graded cellular algebra
with graded cell datum $(\P,T,C, \deg)$ and shift automorphism
$\sigma=(\sigma_A,\sigma_\P,\sigma_T)$ over the integral domain $R$ containing a primitive $\sp$th root of unity $\eps$,  where $\sp$ is the order of
$\sigma_A$. Assume that $\sp \cdot 1_R \in R^\times$. Then $A^\sigma$ is a
graded skew cellular algebra with skew cellular datum
$(\Psig, \iota_\sigma, \Tsig, C_\sigma, \deg_\sigma)$.
\end{Theorem}

\begin{proof}
  By construction, the fixed point subalgebra $A^\sigma$ is an
  $R$-subalgebra of $A$, so it remains to check that the quintuple
  $(\Psig, \iota_\sigma, \Tsig, C_\sigma, \deg_\sigma)$
  satisfies the assumptions of \autoref{D:SkewCellular}.


  First note that $\iota_\sigma$ is a poset automorphism of $\Psig$
  since $(\lambda,k)\gesig(\mu,l)$ if and only if $\lambda\gesig\mu$,
  which is if and only if $(\lambda,-k)\gesig(\mu,-l)$. We now check
  \ref{rel:cellular_homogeneous},
  \ref{rel:cellular_basis},
  \ref{rel:cellular_product} and
  \ref{rel:cellular_star} from \autoref{D:SkewCellular}. The first of
  these properties is easy but the others require more work.

  First, if $\s,\t\in \Tsig(\lambda,k)$ then $\sigmavg$ is
  homogeneous of degree~$0$ by \autoref{L:SkewAutoProperties} and
  $\deg(\sigma_T(\s))=\deg(\s)$, for all $s\in T(\lambda)$ by
  \autoref{D:sigma_cellular}(a).  Therefore,
  $\deg(c^{(k)}_{\s\t})=\deg(c_{\s\t})=\deg(\s)+\deg(\t)$, for all
  $\s,\t\in \Tsig(\lambda,k)$. Hence, \ref{rel:cellular_homogeneous}
  holds.

  Next consider \autoref{rel:cellular_basis}. If $a\in A$ then
  $\sigmavg(a)\in A^\sigma$, so by~\autoref{E:ckst} we have $c^{(k)}_{\s\t}\in A^\sigma$,
  for all $\s,\t\in\Tsig(\lambda,k)$ and $(\lambda,k)\in\Psig$.
  To show that $\set{c^{(k)}_{\s\t}}$ is a basis of $A^\sigma$ first
  observe that because $\sp\cdot 1_R\in R^\times$, for any $b\in A^\sigma$, $
b=\sp^{-1}(\sp b)=\sp^{-1}\sum_{k=0}^{\sp-1}\sigma_A^k(b)=\sp^{-1}\sigmavg(b)\in \mathrm{span}_R\set{\sigmavg(a)|a\in A}$.
It follows that the fixed point subalgebra $A^\sigma$ is spanned by
  $\set{\sigmavg(c_{\s\t})}$. By definition,
  \begin{equation}\label{E:sigmaSum}
    \sigmavg(c_{\s\t})=\sum_{k=0}^{\sp-1}\sigma_A^k(c_{\s\t})
                      =\sum_{k=0}^{\sp-1}c_{\sigma_T^k(\s)\sigma_T^k(\t)}.
  \end{equation}
  In particular, $\sigmavg(c_{\s\t})=\sigmavg(c_{\s'\t'})$ whenever
  $\s'=\sigma^k_T(\s)$ and $\t'=\sigma_T^k(\t)$, for some $k\ge0$.
  It follows that the algebra $A^\sigma$ is spanned by the set
  \begin{align*}
    \Ccal &= \set[\big]{\sigmavg(c_{\s\t})|
                \s\in \Tsig(\lambda), \t\in T(\lambda)\text{ and }\lambda\in\Prep}\\
          &= \set[\big]{\sigmavg(c_{\s\sigma_\lambda^j(\t)})|
                \s,\t\in \Tsig(\lambda), 0\le j<\olambdatab\text{ and }\lambda\in\Prep},
  \end{align*}
  where the second equality follows because
  $\sigma_\lambda=\sigma_T^{\so_\lambda}\colon T(\lambda)\bijection T(\lambda)$ is a
  bijection. We claim that $\Ccal$ is linearly
  independent and hence of basis of~$A^\sigma$. By \autoref{E:sigmaSum},
  when we expand $\sigmavg(c_{\s\t})$ in the $c$-basis of $A$ the
  supports of the different elements of~$\Ccal$ are disjoint, so $\Ccal$
  is linearly independent because $\set{c_{\s\t}}$ is a basis of~$A$.
  Finally, observe that, by \autoref{E:ckst}, the transition matrix
  between the basis in $\Ccal$ and $\set{c^{(k)}_{\s\t}}$ is given by
  Vandermonde matrices in $\set{\epsl^{kj}|0\leq k,j<\olambdatab}$, which are invertible over~$R$ since $\olambdatab\cdot 1_R = \prod_{i = 1}^{\olambdatab-1}(1-\epsl^i)$ and $\olambdatab$ divides $\sp$ in $\Z$.
More precisely, for any $\lambda \in \P_\sigma$ and $\s,\t \in \Tsig(\lambda)$
\begin{equation}
\label{E:spans}
\mathrm{span}_R\set[\big]{c^{(k)}_{\s\t} | k \in \{1,\dots,\olambdatab\}} = \mathrm{span}_R\set[\big]{\sigmavg(c_{\s,\sigmatab^j(\t)}) : j \in \{1,\dots,\olambdatab\}}.
\end{equation}
  Hence,
  $\set{c^{(k)}_{\s\t}}$ is a basis of $A^\sigma$, so
  \autoref{rel:cellular_basis} holds.

We now verify~\ref{rel:cellular_product}. Fix $(\lambda,k)\in\Psig$ and
$\s,\t\in \Tsig(\lambda,k)$.
Let $a\in A^\sigma$. Using \autoref{rel:cellular_product} for the $c$-basis of $A$ and the fact that $\sigmavg$ is $A^\sigma$-linear
\begin{align*}
 ac^{(k)}_{\s\t}    &=\sum_{j = 0}^{\olambdatab - 1} \epsl^{kj}\sigmavg(a c_{\s,\sigma_\lambda^j(\t)}) \\
    &=\sum_{j = 0}^{\olambdatab - 1} \sum_{\v\in T(\lambda)} \epsl^{kj}r_{\v\s}(a)\sigmavg(
           c_{\v,\sigma_\lambda^j(\t)}) + \sigmavg(b),
\end{align*}
for some $b\in A^{\gdom\lambda}$.

By direct calculation,
\begin{align*}
ac^{(k)}_{\s\t} - \sigmavg(b)
&=
\sum_{j = 0}^{\olambdatab - 1} \sum_{\v\in T(\lambda)} \epsl^{kj}r_{\v\s}(a)\sigmavg(
         c_{\v,\sigma_\lambda^j(\t)} )
\\
&=
\sum_{j = 0}^{\olambdatab - 1} \sum_{\v\in T^{\sigma}(\lambda)}\sum_{l = 0}^{\olambdatab-1} \epsl^{kj}r_{\sigma_\lambda^l(\v),\s}(a)\sigmavg(
         c_{\sigma_\lambda^l(\v),\sigma_\lambda^j(\t)} )
\\
&=
\sum_{j = 0}^{\olambdatab - 1} \sum_{\v\in T^{\sigma}(\lambda)}\sum_{l = 0}^{\olambdatab-1} \epsl^{kj}r_{\sigma_\lambda^l(\v),\s}(a)\sigmavg(
         c_{\v,\sigmatab^{(j-l)}\t} )
\\
&=
\sum_{j = 0}^{\olambdatab - 1} \sum_{\v\in T^{\sigma}(\lambda)}\sum_{l = 0}^{\olambdatab-1} \epsl^{k(j+l)}r_{\sigma_\lambda^l(\v),\s}(a)\sigmavg(
         c_{\v,\sigma_\lambda^j(\t)} )
\\
&=
\sum_{\v\in T^{\sigma}(\lambda)}
\sum_{l = 0}^{\olambdatab-1} \epsl^{kl} r_{\sigma_\lambda^l(\v),\s}(a)
\sum_{j = 0}^{\olambdatab - 1} \epsl^{kj}\sigmavg(
         c_{\v,\sigma_\lambda^j(\t)} )
\\
&=
\sum_{\v \in T^{\sigma}(\lambda)} r'_{\v\s}(a) c^{(k)}_{\v\t},
\end{align*}
where
\begin{equation}
\label{E:r'}
r'_{\v\s}(a) = \sum_{l = 0}^{\olambdatab-1} \epsl^{kl} r_{\sigma_\lambda^l(\v),\s}(a),
\end{equation}
does not depend on $\t$. To complete the proof
of~\ref{rel:cellular_product} it suffices to prove that
$\sigmavg(b)\in (A^\sigma)^{\gesig(\lambda,k)}$, where $(A^\sigma)^{\gesig(\lambda,k)}$ is the
$R$-submodule of $A^\sigma$ spanned by
\[
     \set[\big]{c^{(l)}_{\u\v} | \u, \v \in \Tsig(\mu,l)\text{ for }(\mu, l) \in
            \Psig\text{ with } (\mu,l) \gesig (\lambda,k)}.
\]
Since $b\in A^{\gdom\lambda}$, it suffices to prove that
$\sigmavg(c_{\u\v})\in(A^\sigma)^{\gesig(\lambda,k)}$, whenever $\u,\v\in
T(\mu)$ and $\mu\gdom \lambda$. Let $\mu_0 \in \P_\sigma$ be the representative of $\mu$ under the action of $\Z_\P$ on $\P$. Then $\mu = \sigma^k \mu_0$ for some $k \in \Z$ and since $\mu \gdom \lambda$ we deduce that $\mu_0 \gesig \lambda$ (recalling that $\lambda \in \P_\sigma$). Finally, if $l \in \Z$ is such that $\u_0 \coloneqq \sigma^l \u \in \Tsig(\mu_0)$ then $\sigmavg(c_{\u\v}) = \sigmavg(c_{\u_0\v_0})$, where $\v_0=\sigma^l \v$. Now if $\v_0 = \sigmatab[\mu_0]^m \v_1$ with $\v_1 \in \Tsig(\mu_0)$, by~\eqref{E:spans} the element $\sigmavg(c_{\u\v}) = \sigmavg(c_{\u_0\v_0})$ is in
\[
\mathrm{span}_R\set[\big]{c^{(j)}_{\u_0\v_1} | j \in \{1,\dots,\olambdatab[\mu_0]\}}.
\]
Thus, $\sigmavg(c_{\u\v})\in (A^\sigma)^{\gesig (\lambda,k)}$ since $\mu_0 \gesig \lambda$.
This proves that $A^\sigma$ satisfies
\autoref{rel:cellular_product}.
%

Finally, we prove~\ref{rel:cellular_star}. By \autoref{L:SkewAutoProperties}, the anti-isomorphism $*$ of $A$ restricts to an anti-isomorphism of $A^\sigma$. Moreover, if
$\s,\t \in \Tsig(\lambda,k)$ and $(\lambda,k) \in \Psig$ then
\begin{align*}
  \bigl(c^{(k)}_{\s\t}\bigr)^*
    &= \sum_{j = 0}^{\olambdatab-1} \epsl^{kj}\sigmavg(c^*_{\s,\sigma_\lambda^j(\t)})
     = \sum_{j = 0}^{\olambdatab-1} \epsl^{kj}\sigmavg(c_{\sigma_\lambda^j(\t),\s}) \\
    &= \sum_{j = 0}^{\olambdatab-1} \epsl^{kj}\sigmavg(c_{\t,\sigmatab^{-j\so_\lambda}\s})
     = \sum_{j = 0}^{\olambdatab-1} \epsl^{-kj}\sigmavg(c_{\t,\sigma_\lambda^j(\s)}) \\
    &= c^{(-k)}_{\t\s}.
\end{align*}
This completes the proof that $A^\sigma$ is a graded skew cellular algebra.
\end{proof}

%

It is an interesting question whether every skew cellular algebra arises
in this way, that is, as the fixed point subalgebra of a cellular algebra.


\begin{Corollary}
\label{C:skew_cellular_cellular}
In the setting of~\autoref{proposition:sigma_cellular_implies_extended},
if $\olambdatab\leq 2$ for all $\lambda \in \P$ then  $A^\sigma$ is a
graded \emph{cellular} algebra with cellular datum $(\Psig, \Tsig,
C_\sigma, \deg_\sigma)$. In particular, if $2$ is invertible in $R$ and
$A$ is a graded cellular algebra with graded cell datum $(\P,T,C, \deg)$
and shift automorphism $\sigma=(\sigma_A,\sigma_\P,\sigma_T)$ such that
$\sigma_A$ has order $2$ then $A^\sigma$ is a graded cellular algebra
with cell datum $(\Psig, \Tsig, C_\sigma, \deg_\sigma)$.
\end{Corollary}

\begin{proof}
  By~\autoref{proposition:sigma_cellular_implies_extended}, the algebra
  $A^\sigma$ is graded skew cellular with skew cellular datum
  $(\Psig,\iota_\sigma, \Tsig,C_\sigma,\deg_\sigma)$. By construction,
  the involution $\iota_\sigma$ of $\Psig$ is given by
  $\iota_\sigma(\lambda,k) = (\lambda,-k)$. Since $\olambdatab \leq 2$
  we have $k = -k$ in $\Z/\olambdatab\Z$, so $\iota_\sigma =
  \id_{\Psig}$. Moreover, still by construction, the map
  $(\iota_\sigma)_{(\lambda,k)}\map{ \Tsig(\lambda,k)}\Tsig(\lambda,-k)$
  is the identity map of $\Tsig(\lambda)$. Hence,
  $(\iota_\sigma)_{(\lambda,k)}$ is the identity map of
  $\Tsig(\lambda,k)$. Recalling~\autoref{R:cellularAlgebras}, this
  proves the first statement. We deduce the second statement by noting
  that if $2 \in R^\times$ then $-1_R \neq 1_R$ in $R$, so $-1_R$ is a
  primitive square root of unity in~$R$.
\end{proof}

\begin{Example}\label{Ex:SkewCellular}
  Maintain the notation from~\autoref{example:sigma_cellular}. In
  particular, $A=\mathrm{Mat}_n(R)$ has graded cellular
  basis $\set{F_{ij}|1\le i,j\le n}$ and $\sigma$ is the shift
  automorphism of $A$ given by $\sigma(F_{ij})=F_{w(i)w(j)}$, where
  $w\in\Sym_n$ is a permutation such that $w$ is the product of $p$
  disjoint $\frac np$-cycles and $w(n+1-i)=n+1-w(i)$, for $1\le i\le n$.
  Hence, $\sp=p$ and we need to assume that $R$ contains a primitive
  $p$th root of unity and that $p\cdot 1_R \in R^\times$. The reader can check that the skew cellular
  subalgebra of $\sigma$-fixed points is
  \[
    A^\sigma = \set[\big]{M=(m_{ij})\in A| m_{ij}=m_{w(i)w(j)} \text{ for } 1\le i, j\le n}.
  \]
  Possible choices for $w$ include the permutations $(1,2)(3,4)$,
  $(1,3)(2,4)$ or $(1,4)(2,3)$ when $n=4$ and $(1,2,3)(6,5,4)$ when
  $n=6$. In these examples we have $p = 2$ thus $A^\sigma$ is in fact cellular by~\autoref{C:skew_cellular_cellular}. An example where $A^\sigma$ is skew cellular but where~\autoref{C:skew_cellular_cellular} does not apply is with $w = (1,3)(2,5)(4,6)$, where $p = 3$.  Of course, if $w$ is trivial (which happens when $p = n$) then $A^\sigma=A$ is a cellular algebra. For
  example, $w$ must be trivial if $n$ is odd by
  \autoref{example:sigma_cellular}.

  We note that Xi and
  Zhang~\cite[Theorem~4.5]{XiZhang:CentralizerAlgebras} have shown that
  the algebra $A^\sigma$ is cellular; their proof is based on Jordan reduction.
\end{Example}

\subsection{Clifford theory}
\label{subsection:Clifford}

In this section we explicitly describe how Clifford theory works for skew cellular algebras that are obtained using a shift automorphism, as in~\autoref{proposition:sigma_cellular_implies_extended}. The results in this section should be compared with~\cite[\textsection 3.7]{HuMathas:DecHrpn}.

If $M$ is any $A$-module, let $\prescript{\sigma}{}{M}$ be the $A$-module $M$ where the action of $A$ is twisted by $\sigma$. In other words, for any $m \in M$ and $a \in A$ we have
\[
a \cdot_{\prescript{\sigma}{}{M}} m = \sigma(a) \cdot_M m.
\]
Let $M\restr$ be the restriction of $M$ to an $A^\sigma$-module. If $N$ is an $A^\sigma$-module, we write $N\ind$ for the induced $A$-module.

Recall that $A$ is a graded cellular algebra with a shift
automorphism~$\sigma$.

\begin{Lemma}
\label{L:rsigmavsigmas}
 Using the notation of~\autoref{rel:cellular_product}, $r_{\sigma\v,\sigma\s}\bigl(\sigma(a)\bigr) = r_{\v\s}(a)$, for all $a \in A$, $\lambda \in \P$ and $\v,\s \in T(\lambda)$. Consequently, $\phi_\lambda(c_\v,c_\s) = \phi_{\sigma\lambda}(c_{\sigma(\v)},c_{\sigma(\s)})$.
\end{Lemma}

\begin{proof}
By~\autoref{rel:cellular_product} we have
\[
a c_{\s\t} = \sum_{\v \in T(\lambda)} r_{\v\s}(a) c_{\v\t} \pmod {A^{\gdom \lambda}}.
\]
Applying $\sigma$ using \autoref{D:sigma_cellular}(b) and~\autoref{L:SkewAutoProperties},
\[
\sigma(a) c_{\sigma\s,\sigma\t} = \sum_{\v \in T(\lambda)} r_{\v\s}(a) c_{\sigma\v,\sigma\t} \pmod{A^{\gdom \sigma\lambda}},
\]
Hence, the first equality follows because $\{c_{\u\v}\}$ is an $R$-basis of~$A$. In turn, this implies the second equality by \autoref{D:cell_Form}.
\end{proof}

\begin{Proposition}\label{sigmatwist}
  Let $\lambda \in \P$. The $R$-linear map $\gamma_\lambda\map{C_\lambda}{C_{\sigma\lambda}}$ defined by $c_\s \mapsto c_{\sigma\s}$, for $\s\in T(\lambda)$, induces isomorphisms of graded $A$-modules $C_\lambda \simeq \prescript{\sigma}{}{C}_{\sigma\lambda}$ and $D_\lambda \simeq \prescript{\sigma}{}{D}_{\sigma\lambda}$.
\end{Proposition}

\begin{proof}
By definition, $\gamma_\lambda$ is an isomorphism of  $R$-modules. To show that it is an $A$-module isomorphism, suppose that $a \in A$. Then, using~\autoref{L:rsigmavsigmas},
\begin{align*}
  \gamma_\lambda(ac_\s)
  = \sum_{\v \in T(\lambda)} r_{\v\s}(a) \gamma_\lambda(c_\v)
  = \sum_{\v \in T(\lambda)} r_{\sigma\v,\sigma\s}
      \bigl(\sigma(a)\bigr) c_{\sigma\v}
  = \sigma(a) c_{\sigma\s},
\end{align*}
which proves the first isomorphism. The fact that this is an isomorphism of graded modules comes from~\autoref{D:sigma_cellular}a).  To prove that $D_\lambda \simeq \prescript{\sigma}{}{D}_{\sigma\lambda}$ it suffices to show that $\gamma_\lambda \bigl(\rad(C_\lambda)\bigr) = \rad(C_{\sigma\lambda})$, which follows because $\phi_\lambda(c_\s,c_\t) = \phi_{\sigma\lambda}(c_{\sigma\s},c_{\sigma\t})$ by~\autoref{L:rsigmavsigmas}.
\end{proof}

We now assume that $\eps \in R$ is a primitive $\sp$-th root of unity and that $\sp \cdot 1_R \in R^\times$, so that~\autoref{proposition:sigma_cellular_implies_extended} applies. For any $(\lambda,k) \in \Psig$, let $C^{(k)}_\lambda$ be the associated skew cell module of~$A^\sigma$, with $R$-basis $\set{c_\s^{(k)} | \s \in \Tsig(\lambda,k)}$ and $R$-bilinear form $\phi_\lambda^{(k)}$,  which is not symmetric in general.


\begin{Lemma}
\label{L:phik}
Let $(\lambda,k) \in \Psig$ and $\s,\t \in \Tsig(\lambda,k)$. Then
\[
\phi_\lambda^{(k)}\bigl(c_\s^{(k)},c_\t^{(k)}\bigr) = \plam \sum_{l = 0}^{\olambdatab-1} \epsl^{kl} \phi_\lambda\bigl(c_{\sigmatab^l(\s)},c_\t\bigr).
\]
\end{Lemma}

\begin{proof}
Note that the maps $a \mapsto r_{\s\t}(a)$  and $b \mapsto r'_{\s\t}(b)$ for $a \in A$ and $b\in A^\sigma$ are $R$-linear. Unravelling the definitions, if $\u\in\Tsig(\lambda)$ then
\begin{align*}
\phi_\lambda^{(k)}\bigl(c_\s^{(k)},c_\t^{(k)}\bigr)
&=
r'_{\u\t}\bigl(c_{\u\s}^{(k)}\bigr)
=
\sum_{j = 0}^{\olambdatab-1} \epsl^{kj} r'_{\u\t} \bigl(\sigmavg(c_{\u,\sigmatab^j(\s)})\bigr).
\end{align*}
Therefore, by~\eqref{E:r'},
\begin{align*}
\phi_\lambda^{(k)}\bigl(c_\s^{(k)},c_\t^{(k)}\bigr)
&=
\sum_{j,l = 0}^{\olambdatab-1} \epsl^{kj} \epsl^{kl} r_{\sigmatab^l(\u),\t} \bigl(\sigmavg(c_{\u,\sigmatab^j(\s)})\bigr)
\\
&=
\sum_{j,l = 0}^{\olambdatab-1} \sum_{m=0}^{\sp-1} \epsl^{k(j+l)} r_{\sigmatab^l(\u),\t} \bigl(c_{\sigma^m(\u),\sigma^m\sigmatab^j(\s)}\bigr).
\end{align*}
Now by~\autoref{E:r_at(c_us)} we have $r_{\sigmatab^l(\u),\t} \bigl(c_{\sigma^m\u,\sigma^m\sigmatab^j(\s)}\bigr) \neq 0$ only if $\sigmatab^l(\u) = \sigma^m(\u)$. In particular, this implies that $\sigma^m(\u) \in T(\lambda)$, so $\so_\lambda \mid m$. Writing $m = \so_\lambda m'$ we obtain $\sigmatab^l(\u) = \sigmatab^{m'}(\u)$. Now write $m' = a \olambdatab + m''$, with $0 \leq m'' < \olambdatab$. Then $\sigmatab^{m'}(\u) = \sigmatab^{m''}(\u)$ and so $m'' = l$, since the $\Z_\lambda$-orbit of $\u$ has exactly size $\olambdatab$. Recalling that $\olambdatab$ divides $\plam$, we have $r_{\sigmatab^l(\u),\t} \bigl(c_{\sigma^m\u,\sigma^m\sigmatab^j(\s)}\bigr) \neq 0$, for $0 \leq m < \sp$ only if $m = \so_\lambda(a\olambdatab+l)$ with $0 \leq a < \frac{\plam}{\olambdatab}$, in which case $\sigma^m(\v) = \sigmatab^l(\v)$ for all $\v \in T(\lambda)$. We thus obtain
\begin{align*}
\phi_\lambda^{(k)}\bigl(c_\s^{(k)},c_\t^{(k)}\bigr)
&=
\frac{\plam}{\olambdatab}
\sum_{j,l = 0}^{\olambdatab-1} \epsl^{k(j+l)} r_{\sigmatab^l(\u),\t} \bigl(c_{\sigmatab^l(\u),\sigmatab^{l+j}(\s)}\bigr)
\\
&=
\frac{\plam}{\olambdatab}
\sum_{j,l = 0}^{\olambdatab-1} \epsl^{k(j+l)} \phi_\lambda\bigl(c_{\sigmatab^{l+j}(\s)},c_\t\bigr)
\\
&=
\frac{\plam}{\olambdatab}
\sum_{j = 0}^{\olambdatab-1}\sum_{l = 0}^{\olambdatab-1} \epsl^{kl} \phi_\lambda\bigl(c_{\sigmatab^l(\s)},c_\t\bigr)
\\
&=
\plam \sum_{l = 0}^{\olambdatab-1} \epsl^{kl} \phi_\lambda\bigl(c_{\sigmatab^l(\s)},c_\t\bigr),
\end{align*}
as desired.
\end{proof}

\begin{Proposition}
\label{P:Clambda_sum_Clambdak}
Let $\lambda \in \Prep$. The $R$-linear map
\begin{align*}
\gamma'_\lambda : \bigoplus_{k = 0}^{\olambdatab-1} C^{(k)}_\lambda &\longrightarrow C_\lambda;\quad
c^{(k)}_\s \longmapsto \sum_{j = 0}^{\olambdatab-1} \epsl^{-kj} c_{\sigmatab^j(\s)},
 \end{align*}
  induces isomorphisms of graded $A^\sigma$-modules
\begin{align*}
\bigoplus_{k = 0}^{\olambdatab-1} C^{(k)}_\lambda \simeq C_\lambda\restr,
\qquad\text{and}\qquad
\bigoplus_{k = 0}^{\olambdatab-1} D^{(k)}_\lambda \simeq D_\lambda\restr.
\end{align*}
\end{Proposition}

\begin{proof}
  First note that $\gamma'_\lambda$ is homogeneous since $\deg^\sigma(\s) = \deg\s=\deg\sigma^j(\s)$ for all $\s \in \Tsig(\lambda,k)=\Tsig(\lambda)$ and all $j \in \Z$ by~\autoref{D:sigma_cellular}a). By the Vandermonde determinant  argument that we used in the proof of~\autoref{proposition:sigma_cellular_implies_extended}, the map $\gamma'_\lambda$ sends a basis to a basis, so is an $R$-module isomorphism.

  We prove that $\bigoplus_k C^{(k)}_\lambda\cong C_\lambda\restr$ as $A^\sigma$-modules. Recall that if $\s \in T(\lambda)$ and $a \in A$ then
\[
a c_\s = \sum_{\v \in T(\lambda)} r_{\v\s}(a) c_\v = \sum_{\v \in \Tsig(\lambda)} \sum_{l = 0}^{\olambdatab-1} r_{\sigmatab^l(\v),\s}(a) c_{\sigmatab^l(\v)},
\]
in~$C_\lambda$. Similarly,
if $\s \in \Tsig(\lambda,k)$ and $a \in A^\sigma$ then
\[
a c_\s^{(k)} = \sum_{\v \in \Tsig(\lambda,k)} r'_{\v\s}(a) c_{\v}^{(k)},
\qquad\text{where}\qquad
r'_{\v\s}(a) = \sum_{j = 0}^{\olambdatab-1} \epsl^{kj} r_{\sigmatab^j(\v),\s}(a),
\]
in $C_\lambda^{(k)}$.
For any $k \in \Z/\olambdatab\Z$, $\s \in \Tsig(\lambda,k)$  and $a \in A^\sigma$ we have
\begin{align*}
a \gamma'_\lambda\bigl(c_\s^{(k)}\bigr)
=
\sum_{j = 0}^{\olambdatab-1} \epsl^{-kj} a c_{\sigmatab^j(\s)}
=
\sum_{j = 0}^{\olambdatab-1} \epsl^{-kj} \sum_{\v \in \Tsig(\lambda)} \sum_{l = 0}^{\olambdatab-1} r_{\sigmatab^l(\v),\sigmatab^j(\s)}(a) c_{\sigmatab^l(\v)}.
\end{align*}
Now by~\autoref{L:rsigmavsigmas} we have $r_{\sigmatab^l(\v),\sigmatab^j(\s)}(a) = r_{\sigmatab^{l-j}(\v),\s}(a)$ since $a \in A^\sigma$. Thus, we obtain, recalling that $\epsl$ is an $\olambdatab$-th root of unity,
\begin{align*}
 a \gamma'_\lambda\bigl(c_\s^{(k)}\bigr)
 &=
 \sum_{\v \in \Tsig(\lambda)} \sum_{l = 0}^{\olambdatab-1} \epsl^{-kl} \sum_{j = 0}^{\olambdatab-1} \epsl^{k(l-j)} r_{\sigmatab^{l-j}(\v),\s}(a) c_{\sigmatab^l(\v)}
 \\
 &=
 \sum_{\v \in \Tsig(\lambda)} \sum_{l = 0}^{\olambdatab-1} \epsl^{-kl} \sum_{j = 0}^{\olambdatab-1} \epsl^{kj} r_{\sigmatab^j(\v),\s}(a) c_{\sigmatab^l(\v)}
 \\
 &=
 \sum_{\v \in \Tsig(\lambda)} \sum_{l = 0}^{\olambdatab-1} \epsl^{-kl} r'_{\v\s}(a) c_{\sigmatab^l(\v)}
 \\
 &=
 \sum_{v \in \Tsig(\lambda)} r'_{\v\s}(a) \gamma'_\lambda\bigl(c_\v^{(k)}\bigr)
 \\
 &=
 \gamma'_\lambda\bigl(a c_\s^{(k)}\bigr).
\end{align*}
This proves that $\gamma'_\lambda$ is $A^\sigma$-linear, thus
establishing that
$\bigoplus_k C^{(k)}_\lambda\cong C_\lambda\restr$ as $A^\sigma$-modules.

To prove the second isomorphism, it suffices to prove that
\begin{equation}
\label{E:gamma'rad}
\bigoplus_{k = 0}^{\olambdatab-1} \gamma'_\lambda\left(\rad C_\lambda^{(k)}\right) = \rad C_\lambda.
\end{equation}
Let $x = \sum_{\s \in \Tsig(\lambda)} x_\s c_\s^{(k)} \in \rad C_\lambda^{(k)}$ with $x_\s \in R$.


By~\autoref{L:phik}, if $\t \in \Tsig(\lambda)$ then
\[
0=
\phi_\lambda^{(k)}\bigl(c_\t^{(k)},x\bigr)
=
\sum_{\s \in \Tsig(\lambda)}x_\s \phi_\lambda^{(k)}\bigl(c_\t^{(k)},c_\s^{(k)}\bigr)
=
\plam\sum_{\s \in \Tsig(\lambda)} \sum_{j = 0}^{\olambdatab-1} \epsl^{kj} x_\s \phi_\lambda \bigl(c_{\sigmatab^j(\t)},c_{\s}\bigr).
\]
Thus, if $\t \in \Tsig(\lambda)$ and $0 \leq l < \olambdatab$ then, using~\autoref{L:rsigmavsigmas},
\begin{align*}
\phi_\lambda\bigl(c_{\sigmatab^l\t}, \gamma'_\lambda(x)\bigr)
&=
\sum_{\s \in \Tsig(\lambda)} x_\s \phi_\lambda\left(c_{\sigmatab^l(\t)}, \gamma'_\lam(c_\s^{(k)})\right)
\\
&=
\sum_{\s \in \Tsig(\lambda)} \sum_{j = 0}^{\olambdatab-1} \epsl^{-kj} x_\s \phi_\lambda\left(c_{\sigmatab^l(\t)}, c_{\sigmatab^j(\s)}\right)
\\
&=
\epsl^{-kl}\sum_{\s \in \Tsig(\lambda)} \sum_{j = 0}^{\olambdatab-1} \epsl^{k(l-j)} x_\s \phi_\lambda\left(c_{\sigmatab^{l-j}(\t)}, c_{\s}\right)
\\
&=
\epsl^{-kl}\sum_{\s \in \Tsig(\lambda)} \sum_{j = 0}^{\olambdatab-1} \epsl^{kj} x_\s \phi_\lambda\left(c_{\sigmatab^j(\t)}, c_{\s}\right)
\\
&=
0,
\end{align*}
proving that $\gamma'_\lambda(x) \in \rad C_\lambda$. To prove that $\rad C_\lambda \subseteq \oplus_k \gamma'_\lambda\bigl(\rad C_\lambda^{(k)}\bigr)$, first note that if $\s \in \Tsig(\lambda)$ and $j \in \Z/\olambdatab\Z$ then
\[
{\gamma'_\lambda}^{-1}\bigl(c_{\sigmatab^j\s}\bigr)= \frac1\olambdatab \sum_{k = 0}^{\olambdatab-1} \epsl^{kj} c_\s^{(k)}.
\]
Let $x = \sum_{\s \in \Tsig(\lambda)} \sum_{j = 0}^{\olambdatab-1} x_{\s,j} c_{\sigmatab^j(\s)} \in \rad C_\lambda$.
 We have
\[
{\gamma'_\lambda}^{-1}(x)
=
\sum_{\s \in \Tsig(\lambda)} \sum_{j = 0}^{\olambdatab-1} x_{\s,j} {\gamma'_\lambda}^{-1}\bigl(c_{\sigmatab^j(\s)}\bigr)
=\frac1\olambdatab
\sum_{\s \in \Tsig(\lambda)} \sum_{j,k = 0}^{\olambdatab-1} \epsl^{kj} x_{\s,j} c_\s^{(k)}.
\]
So, to complete the proof it is enough to show that if $0 \leq k <
\olambdatab$ then
\[
\sum_{\s \in \Tsig(\lambda)} \sum_{j = 0}^{\olambdatab-1} \epsl^{kj} x_{\s,j} c_\s^{(k)} \in \rad C_\lambda^{(k)},
\]
Using~\autoref{L:phik} and~\autoref{L:rsigmavsigmas},
if $\t \in \Tsig(\lambda)$ then
\begin{align*}
\sum_{\s \in \Tsig(\lambda)} \sum_{j = 0}^{\olambdatab-1} \epsl^{kj} x_{\s,j} \phi_\lambda^{(k)}\bigl(c_\t^{(k)},c_\s^{(k)}\bigr)
&=
\plam \sum_{\s \in \Tsig(\lambda)} \sum_{j,l = 0}^{\olambdatab-1}  \epsl^{k(j+l)} x_{\s,j} \phi_\lambda\bigl(c_{\sigmatab^l(\t)},c_{\s}\bigr)
\\
&=
\plam \sum_{\s \in \Tsig(\lambda)} \sum_{j,l = 0}^{\olambdatab-1}  \epsl^{k(j+l)} x_{\s,j} \phi_\lambda\bigl(c_{\sigmatab^{l+j}(\t)}, c_{\sigmatab^j(\s)}\bigr)
\\
&=
\plam \sum_{\s \in \Tsig(\lambda)} \sum_{j,l = 0}^{\olambdatab-1}  \epsl^{kl} x_{\s,j} \phi_\lambda\bigl(c_{\sigmatab^l(\t)}, c_{\sigmatab^j(\s)}\bigr)
\\
&=
\plam\sum_{l = 0}^{\olambdatab-1} \epsl^{kl}\phi_\lambda\bigl(c_{\sigmatab^l(\t)}, x\bigr)
\\
&=
0,
\end{align*}
where the last equality comes from the fact that $x \in \rad C_\lambda$.
\end{proof}


\begin{Definition}
\label{E:z_strong}
  The automorphism $\sigma_A$ is \emph{$\varepsilon$-splittable} if there exists an
  invertible element~$z\in A^\times$ that is homogeneous of
  degree~$0$ such that $\sigma_A( z) = \eps  z$.
\end{Definition}

Fix $z \in A^\times$ as in \autoref{E:z_strong}. Then
$z^i A^\sigma = \ker(\sigma_A - \eps^i)$, for $i\ge0$. The terminology of
\autoref{E:z_strong} is justified because if $R$ is a field then we can decompose~$A$ into a direct sum of $\sigma_A$-eigenspaces
\begin{equation}
\label{E:A_free}
A = \bigoplus_{i = 0}^{\sp-1}  z^i A^\sigma.
\end{equation}
since
$\sigma_A$ has order $\sp$. In particular, if $\sigma_A$ is $\varepsilon$-splittable
then $A$ is free, and hence projective, as an $A^\sigma$-module.


Recall from~\autoref{sigmatwist} that $\gamma_\lambda \map{ C_\lambda}{ C_{\sigma\lambda}}$ is an $R$-linear isomorphism such that $\gamma_\lambda(ax) = \sigma(a)\gamma_\lambda(x)$, for all $a \in A$ and $x \in C_\lambda$. For $j \geq 0$ define the $R$-isomorphism $\compgamma \map{C_\lambda}{C_{\sigma^j\lambda}}$ by
\begin{equation}
\label{E:compgamma}
\compgamma = \gamma_{\sigma^{j-1}\lambda} \circ \dots \circ \gamma_\lambda,
\end{equation}
and set $\sigmaclam \coloneqq \compgamma[\so_\lambda]$, an $R$-automorphism of $C_\lambda$.
Then $\sigmaclam(c_\s) = c_{\sigmatab\s}$, for all $\s \in T(\lambda)$. In particular $\sigmaclam$ has order $\olambdatab$.  The $R$-linear isomorphism $\compgamma$ satisfies
\begin{align}
\label{E:compgamma(ax)}
\compgamma(ax) &= \sigma_A^j(a)\compgamma(x),
  &&\text{ for all $a \in A$ and $x \in C_\lambda$}.
\intertext{In particular $\sigmaclam$ satisfies}
\label{E:sigmaclam(ax)}
\sigmaclam(ax) &= \sigma_A^{\so_\lambda}(a) \sigmaclam(x),
  &&\text{for all $a \in A$ and $x \in C_\lambda$.}
\end{align}
%
%
%

\begin{Proposition}
\label{P:olambdatab=plam}
Suppose that $\sigma_A$ is $\varepsilon$-splittable. Then
$\olambdatab = \plam$ for all $\lambda \in \P$.
\end{Proposition}

\begin{proof}
Under the isomorphism $\gamma'_\lambda$ of~\autoref{P:Clambda_sum_Clambdak},
the cell module $C^{(k)}_\lambda$ is sent into the eigenspace $\ker\bigl( \sigmaclam - \epsl^k\bigr)$. Since $\epsl \in R^\times$, these eigenspaces are in direct sum and we conclude that
\begin{equation}
\label{E:C_lambda_sum_eigenspaces}
C_\lambda = \bigoplus_{k = 0}^{\olambdatab-1} \ker\bigl(\sigmaclam - \epsl^k\bigr).
\end{equation}
Now let $x \in C_\lambda$ be an eigenvector for $\sigmaclam$ with eigenvalue $\epsl^k$. By~\autoref{E:sigmaclam(ax)},  we have
\begin{align*}
\sigmaclam(zx)
&=
\sigma_A^{\so_\lambda}(z)\sigmaclam(x)
=
\eps^{\so_\lambda} z\epsl^k x
=
\eps^{\so_\lambda} \epsl^k zx.
\end{align*}
Since $z \in A^\times$, we have $zx \neq 0$ and thus $zx$ is an eigenvector for $\sigmaclam$ with eigenvalue $\eps^{\so_\lambda} \epsl^k$. By~\autoref{E:C_lambda_sum_eigenspaces}, this implies that $\eps^{\so_\lambda} \in \langle \epsl\rangle$, thus $\eps^{\so_\lambda}$ is an $\olambdatab$-th root of unity. But $\eps^{\so_\lambda}$ has order $\frac{\sp}{\so_\lambda} = \plam$ thus $\plam \mid \olambdatab$ thus $\plam = \olambdatab$ since $\olambdatab \mid \plam$.
\end{proof}

In particular,~\autoref{P:olambdatab=plam} implies that $\epsl = \eps^{\so_\lambda}$. Now let $\tau$ be the homogeneous automorphism of~$A$ given by $\tau(a)=\zl a \zl^{-1}$, for $a\in A$. Note that $A^\sigma$ is stable under $\tau$. The next result is a complement to~\autoref{P:Clambda_sum_Clambdak}.

%

\begin{Proposition}
\label{P:tau}
Let $(\lambda,k) \in \Psig$. The map
\[\gamma_\lambda^{(k)} : C_\lambda^{(k)} \to C_\lambda^{(k+1)};\quad
  x \mapsto {\gamma'_\lambda}^{-1}\bigl( \zl \gamma_\lambda'(x)\bigr)
\]
induces graded $A^\sigma$-module isomorphisms
\[
{C_\lambda^{(k)}} \simeq \prescript{\tau}{}{C_\lambda^{(k+1)}}
\qquad\text{and}\qquad
{D_\lambda^{(k)}} \simeq \prescript{\tau}{}{D_\lambda^{(k+1)}}.
\]
\end{Proposition}


\begin{proof} During the proof of~\autoref{P:olambdatab=plam} we obtained $\gamma'_\lambda\bigl(C^{(k)}_\lambda\bigr) = \ker\bigl(\sigmaclam - \epsl^k\bigr)$ and
\begin{equation}
\label{E:zgamma'Clamk}
\zl \gamma'_\lambda\bigl(C_\lambda^{(k)}\bigr) = \gamma'_\lambda\bigl(C_\lambda^{(k+1)}\bigr).
\end{equation}
Thus, the map $\gamma^{(k)}_\lambda$ is well defined. Moreover, it is clearly bijective since $z \in A^\times$.
Equation~\autoref{E:zgamma'Clamk} implies that there is an $A^\sigma$-module isomorphism ${C_\lambda^{(k)}} \simeq
\prescript{\tau}{}{C_\lambda^{(k+1)}}$ because if
 $x \in C_\lambda^{(k)}$ and $a \in A^\sigma$ then
\begin{align*}
\gamma_\lambda^{(k)}(ax)
&= {\gamma'_\lambda}^{-1}\bigl(\zl \gamma_\lambda'(ax)\bigr)
 = {\gamma'_\lambda}^{-1}\bigl(\zl a \gamma'_\lambda(x)\bigr)
 = {\gamma'_\lambda}^{-1}\bigl(\tau(a) \zl \gamma'_\lambda(x)\bigr)
\\
&= \tau(a) {\gamma'_\lambda}^{-1}\bigl( \zl \gamma'_\lambda(x)\bigr)
 = \tau(a) \gamma_\lambda^{(k)}(x).
\end{align*}
Moreover the map $\gamma_\lambda^{(k)}$ is homogeneous of degree $0$ since $\zl$ has degree zero (and $\gamma'_\lambda$ is homogeneous).


Finally, by~\autoref{T:SkewCellularSimples}\autoref{I:Dlam_irred_unique_simple_head}, $\rad C^{(l)}_\lambda$ is the Jacobson radical of $C^{(l)}_\lambda$ for all $l$ so, since $\gamma_\lambda^{(k)}$ is an $A^\sigma$-module isomorphism,  $\gamma_\lambda^{(k)}\bigl(\rad C_\lambda^{(k)}\bigr)= \rad  C_\lambda^{(k+1)}$. Hence, $\gamma_\lambda^{(k)}$ induces an isomorphism ${D_\lambda^{(k)}} \simeq \prescript{\tau}{}{D_\lambda^{(k+1)}}$ of $A^\sigma$-modules.
\end{proof}

Since $A^\sigma$ is a skew cellular algebra,
\autoref{T:SkewCellularSimples} gives a classification of the graded
simple $A^\sigma$-modules. Combining the results above we obtain the
following classification of the simple $A^\sigma$-modules in terms of
the simple $A$-modules.


\begin{Theorem}
\label{T:ShiftedSimples}
Let $R$ be a field containing a primitive $\sp$th root of unity $\eps$.
Suppose that $A$ has graded cell datum $(\P,T,C, \deg)$ and a shift
automorphism $\sigma=(\sigma_A,\sigma_\P,\sigma_T)$ such that $\sigma_A$
is $\varepsilon$-splittable and has order $\sp$. Then
\[
  \set{D^{(k)}_\lambda\<s\>|D^\lambda\ne0 \text{ for }(\lambda,k)\in\Psig,
      \text{ and }s\in\Z}
\]
is a complete set of pairwise non-isomorphic graded simple
$A^\sigma$-modules.
\end{Theorem}

\begin{proof}
  First note that because $R$ contains a primitive $\sp$th root of unity
  the characteristic of~$R$ cannot divide $\sp$, so $\sp\cdot1_R$
  is invertible in~$R$. Therefore, $A^\sigma$ is a skew cellular algebra
  with skew cellular datum $(\Psig, \iota_\sigma, \Tsig, C_\sigma, \deg_\sigma)$
  by \autoref{proposition:sigma_cellular_implies_extended}. Therefore, by
  \autoref{T:SkewCellularSimples}, a complete of pairwise non-isomorphic
  graded simple $A^\sigma$-modules is given by the non-zero modules in
  the set $\set{D^{(k)}_\lambda\<s\>|(\lambda,k)\in\Psig, s\in\Z}$.
  By \autoref{P:Clambda_sum_Clambdak}, \autoref{P:olambdatab=plam} and
  \autoref{P:tau}, if $\lambda\in\Prep$ then
  \[
      D^\lambda\ne0 \Longleftrightarrow
      D^{(k)}_\lambda\ne0\text{ for some }0\le k<\plam
      \Longleftrightarrow
      D^{(k)}_\lambda\ne0\text{ for all }0\le k<\plam.
  \]
  Hence, the result follows.
\end{proof}

By \autoref{E:A_free}, if $M$ is an $A^\sigma$-module then the
induced $A$-module $M\ind$ is given by
\[
M\ind \;\simeq \bigoplus_{i = 0}^{\sp - 1}  z^i M,
\]
where the action of $a \in A$ on $ z^iM$ is given by
\[
ax = \sum_{j = 0}^{\sp-1}  z^j a_j x, \qquad\text{for }x\in M
\]
where $a z^i = \sum_j  z^j a_j$ with $a_j \in A^\sigma$.

Recall the $R$-linear maps $\compgamma \map{C_\lambda}{C_{\sigma^j\lambda}}$ from~\autoref{E:compgamma}.

%

\begin{Proposition}
\label{P:induction_Clambdak}
Assume that $R$ is a field and that $\sigma_A$ is $\varepsilon$-splittable.
For $(\lambda,k) \in \Psig$, let
\[
\hat\gamma''_\lambda : C_\lambda^{(k)} \to  \bigoplus_{j = 0}^{\so_\lambda-1} C_{\sigma^j\lambda}
\]
be the $A^\sigma$-linear map whose $j$-th component is given by $\compgamma \circ \gamma'_\lambda$ for $0 \leq j < \so_\lambda$. The unique corresponding $A$-linear map $\gamma''_\lambda \map{ C_\lambda^{(k)}\bigr\uparrow }{\bigoplus_{j = 0}^{\so_\lambda-1} C_{\sigma^j\lambda}}$ induces isomorphisms of graded $A$-modules
\[
\left. C_\lambda^{(k)}\right\uparrow \simeq \bigoplus_{j = 0}^{\so_\lambda-1} C_{\sigma^j\lambda}
\qquad\text{and}\qquad
\left. D_\lambda^{(k)}\right\uparrow \simeq \bigoplus_{j = 0}^{\so_\lambda-1} D_{\sigma^j\lambda}.
\]
\end{Proposition}


\begin{proof}
Since $\left.C_\lambda^{(k)}\right\uparrow \simeq  \oplus_{i = 0}^{\sp} z^i C_\lambda^{(k)}$ as an $A$-module, if $x_i \in C_\lambda^{(k)}$, where $0 \leq i < \sp$, by Frobenius reciprocity we have
\[
\gamma''_\lambda\left(\sum_{i=0}^{\sp-1} z^i x_i\right) = \sum_{i=0}^{\sp-1} z^i \hat \gamma''_\lambda(x_i).
\]
The map $\gamma''_\lambda$ is $A$-linear by construction. Recalling from~\autoref{P:olambdatab=plam} that $\olambdatab = \plam$, thus
\[
\sp \lvert \Tsig(\lambda)\rvert = \so_\lambda \plam \lvert\Tsig(\lambda)\rvert=\so_\lambda \lvert T(\lambda)\rvert.
\]
It follows that the starting and ending $R$-vector spaces have the same dimension. To prove that $\gamma''_\lambda$ is bijective it suffices to prove that it is injective.

Let $(x_i)_{0 \leq i < \sp}$ be as above and assume that $\gamma''_\lambda\left(\sum_{i = 0}^{\sp-1} z^i x_i\right) = 0$. Since $\hat\gamma''_\lambda(x) = \sum_{j = 0}^{\so_\lambda-1} \compgamma\bigl(\gamma'_\lambda(x)\bigr)$ for all $x \in C_\lambda^{(k)}$, we deduce that for all $0 \leq j < \so_\lambda$ we have
\[
\sum_{i = 0}^{\sp - 1} z^i \compgamma\bigl(\gamma'_\lambda(x_i)\bigr) = 0.
\]
Now using~\autoref{E:z_strong} and~\autoref{E:compgamma(ax)}, we deduce that
\[
\sum_{i = 0}^{\sp - 1} \eps^{-ij} \compgamma\bigl(z^i\gamma'_\lambda(x_i)\bigr) = 0.
\]
Since $\compgamma$ is an $R$-isomorphism, we deduce that for all $0 \leq j < \so_\lambda$ we have
\[
\sum_{i = 0}^{\sp-1} \eps^{-ij} z^i \gamma'_\lambda(x_i) = 0.
\]
Writing $i = a\plam + b$ for $0 \leq a < \so_\lambda$ and $0 \leq b < \plam$,
\[
\sum_{a = 0}^{\so_\lambda-1} \sum_{b = 0}^{\plam-1} \eps^{-(a\plam + b)j} z^{a\plam +b} \gamma'_\lambda(x_{a\plam + b}) = 0,
\qquad\text{for all }0 \leq j < \so_\lambda.
\]
By~\autoref{E:zgamma'Clamk},  $z^{a\plam + b} \gamma'_\lambda(C_\lambda^{(k)}) = \gamma'_\lambda(C_\lambda^{(k+b)})$. Thus, using~\autoref{P:Clambda_sum_Clambdak},
\[
\sum_{a = 0}^{\so_\lambda-1} \eps^{-aj\plam} z^{a\plam} \gamma'_\lambda(x_{a\plam+b}) = 0,
\qquad\text{for $0 \leq j < \so_\lambda$ and $0 \leq b < \plam$.}
\]
Since $\eps^{\plam}$ is a primitive $\so_\lambda$-th root of unity,
for a fixed $0 \leq b < \plam$ we obtain an invertible linear system, so $z^{a\plam} \gamma'_\lambda(x_{a\plam+b}) = 0$ in $C_\lambda$, for $0 \leq j,a< \so_\lambda$ and $0 \leq b < \plam$. Since $z \in A^\times$ we deduce that $\gamma'_\lambda(x_i) = 0$, for $0 \leq i < \sp$. Hence, $x_i = 0$ in $C_\lambda^{(k)}$ since $\gamma'_\lambda$ is injective. We conclude that $\gamma''_\lambda$ is injective, proving the first isomorphism of $A$-modules. Note that this isomorphism is homogeneous of degree $0$ since $\deg(z)=0$ and $\hat\gamma''_\lambda$ is homogeneous.

To prove the second isomorphism, by \autoref{sigmatwist} and~\autoref{P:Clambda_sum_Clambdak} we have
\[
\gamma''_\lambda\left(\bigoplus_{i = 0}^{\sp - 1} z^i \rad C_\lambda^{(k)}\right) \subseteq \bigoplus_{j  = 0}^{\so_\lambda - 1} \rad C_{\sigma^j \lambda}.
\]
Moreover we also obtain that
\[
\dim_R \rad C_{\sigma^j \lambda} = \dim_R \rad C_\lambda = \plam \dim_R \rad C_\lambda^{(k)},
\]
for all $0 \leq j < \so_\lambda$. Thus, the above inclusion is an equality and the proof is complete.
\end{proof}


\begin{Remark}
  With a little more care it is possible to prove
  \autoref{P:induction_Clambdak} over an integral domain that contains
  $\eps$.  As in~\cite[\textsection 3.7]{HuMathas:DecHrpn},
  the existence of the isomorphism of~\autoref{P:induction_Clambdak} can
  be deduced from more general results such as~\cite[Proposition
  2.2]{Genet:graded} (and~\cite[Appendix]{Hu:simpleGrpn}). The point of
  \autoref{P:induction_Clambdak} is to give an explicit isomorphism.
\end{Remark}

\section{Hecke algebras and diagrammatic Cherednik algebras}
\label{S:Hecke_algebras}

  Having set up the machinery of skew cellular algebras we are now
  ready to tackle the main results of this paper, which show that the
  Hecke algebras of type $G(\ell,p,n)$ are graded skew
  cellular algebras. To do this we use the cyclotomic KLR algebras of
  type~$A$, together with the diagrammatic Cherednik algebras, to
  construct a shift automorphism of these algebras.

  \subsection{Hecke algebras}\label{SS:CycHecke}
This section recalls the definitions and results from the literature
  that we need about the Hecke algebras of type $G(\ell,p,n)$.
  Throughout this paper we fix positive integers $n$, $p$ and~$d$, with
  $p\ge2$. Recall that $R$ is a commutative integral ring with~$1$. Let $K = \mathrm{Frac}(R)$ be the field of fractions of~$R$. We assume that $K$ contains a primitive $p$th root of unity $\eps$. 
  Set $\ell=pd$ and fix \textbf{cyclotomic parameters}
  $Q_1,\dots,Q_d\in K$.  Set $\bQ=(Q_1,\dots,Q_d)$ and
  \[
    \bvQ=(\eps Q_1,\eps^2Q_1,\dots,\eps^p Q_1,\eps
    Q_2,\dots,\eps^pQ_2,\dots,\eps Q_d,\dots,\eps^pQ_d).
  \]
  Finally, fix an invertible \textbf{Hecke parameter} $\q\in K$.

  \begin{Definition}[Ariki and Koike~\cite{AK}, Brou\'e and Malle~\cite{BM:cyc}]\label{SS:CycHecke1}
  The \textbf{Hecke algebra} of type $G(\ell,1,n)$ with Hecke
  parameter $\q$ and cyclotomic parameters $\bvQ$ is the unital
  associative $K$-algebra $\Hln$ with generators $T_0,T_1,\dots,T_{n-1}$
  and relations:
  { \setlength{\abovedisplayskip}{2pt}
    \setlength{\belowdisplayskip}{1pt}
    \begin{align*}
          \prod_{k=1}^d(T_0^p-Q_k^p) &=0, & T_0T_1T_0T_1 &=T_1T_0T_1T_0,\\
          (T_r-\q)(T_r+1)&=0,         & T_kT_{k+1}T_k&=T_{k+1}T_kT_{k+1},
    \end{align*}
    \[ T_iT_j = T_jT_i \quad\text{if}\quad |i-j|>1,\]
  }
  where $1\leq r<n, 1\leq k<n-1$ and $1\le i,j<n$.
  \end{Definition}

  \begin{Remark}
  The algebra $\Hln$ is in fact a special case of a Hecke
  algebra of type $G(\ell,1,n)$, which can have $\ell$ arbitrary cyclotomic
  parameters $Q_1,\dots,Q_\ell\in K$.
  \end{Remark}

  Inspecting the relations, $\Hln$ has a unique automorphism $\sigma$ of order $p$ such that
  \begin{equation}\label{E:sigma}
      \sigma(T_0)=\eps T_0\quad\text{and}\quad
      \sigma(T_i)=T_i\text{ for }1\le i<n.
  \end{equation}

  We can now define the main (ungraded) algebras of interest in this paper.

  \begin{Definition}[Ariki~\cite{Ariki:Grpn}, Brou\'e and Malle~\cite{BM:cyc}]\label{D:Hlpn}
    The \textbf{Hecke algebra} of type $\Glpn$ with
    parameters $\q\in K^\times$ and $\bvQ\in K^\ell$ is the fixed-point subalgebra
    \[ \Hlpn=\bigl(\Hln\bigr)^\sigma=\set{h\in\Hln|\sigma(h)=h}.\]
  \end{Definition}

  Equivalently, $\Hlpn$ is the subalgebra of $\Hln$ that is generated by
  $T_0^p$, $T_0^{-1}T_1T_0$ and $T_1,\dots,T_{n-1}$. Notice that
  $\Hlpn$ is an Iwahori-Hecke algebra of type~$D$ when
  $\ell=p=2$.

  From the relations it is clear that if $c\in K^\times$ is any non-zero
  scalar then $\Hln(\q,c\bvQ)\cong\Hln(\q,\bvQ)$ and hence that
  $\Hlpn(\q,c\bvQ)\cong\Hlpn(\q,\bvQ)$, where $c\bvQ=(c\eps
  Q_1,\dots,c\eps^p Q_d)$.  Moreover, by \cite{HuMathas:MoritaGrpn}, we can
  assume that the cyclotomic parameters $Q_1,\cdots,Q_d$ are in a
  single $(\eps,\q)$-orbit. That is, $Q_i/Q_j\in\eps^\Z \q^\Z$ for any
  $1\leq i,j\leq d$.

  Let $e = \min\set{e>0|1+\q+\dots+\q^{e-1}=0}$ and set $e=\infty$
  if no such integer exists.  Using Clifford theory, as discussed on
  \cite[Page 3383]{Hu:simpleGrpn}, we can further assume that
  $Q_i = \q^{\rho_i}$,
  \begin{equation}
  \label{E:bvQ}
    \bvQ=(\eps\q^{\rho_1},\eps^2\q^{\rho_1},\dots,\eps^p\q^{\rho_1},
          \eps \q^{\rho_2},\dots,\eps^p\q^{\rho_2},\dots,
          \eps\q^{\rho_d},\dots,\eps^p\q^{\rho_d}).
  \end{equation}
  and that we are in one of the following two cases:

  \Case{$\q^\Z\cap\eps^\Z\neq \set{1}$} Equivalently, $e<\infty$ and
  $\gcd(e,p)>1$. Let $m=\gcd(e,p)$ and write $p=mp'$. We may assume that
  $e=me'$ and that
  $\eps^{p'}=\q^{e'}$ is a primitive $m$th root of unity in $K$. Note
  that $p'=\min\set{1\leq a\leq p|\eps^a\in \q^\Z}$.

  \Case{$\q^\Z\cap\eps^\Z=\set{1}$} Equivalently, either $e<\infty$ and
  $\gcd(e,p)=1$, or $e=\infty$. For consistency of notation with Case~1,
  we assume that $0=\rho_1\le \rho_2\le\dots\le \rho_d$ and set $p'=p$,
  $e'=e$ and $m=1$. In fact, as noted in
  \cite[Corollary~2.10]{HuMathas:SeminormalQuiver}, if $e=\infty$ then
  we can replace~$\q$ with an $\hat e$ root of unity for some
  sufficiently large $\hat e$ without changing the (graded) isomorphism
  type of~$\Hln$. Henceforth, we assume that~$e$ is finite.

  Permuting the integers $\rho_1,\dots,\rho_d$ does not affect $\Hln$ up to
  isomorphism and, similarly, we can replace $\rho_a$ with $\rho_a+e$ since
  $q^{\rho_a+e}=q^{\rho_a}$. In order to be able to construct the basis that
  we need to prove our main results we adopt the following convention.

  \begin{Definition}\label{D:dCharge}
    A \textbf{$d$-charge} is a $d$-tuple of integers
    $\dcharge=(\rho_1,\dots,\rho_d)\in\Z^d$ such that
    \[
          \rho_{a+1}-\rho_a\geq(2n+3)e,\quad\text{ for $1\le a<d$}.
    \]
  \end{Definition}

  We assume that we have a fixed $\dcharge$ for the rest of this paper.
  For convenience, we assume that $\rho_1=0$.  In particular, this
  implies that $0=\rho_1<\rho_2<\dots<\rho_d$.

\subsection{Quiver Hecke algebras}
\label{subsection:QHA}
The quiver Hecke algebras, or KLR algebras, are a remarkable family of
$\Z$-graded algebras that were introduced by Khovanov and
Lauda~\cite{KhovLaud:diagI} and
Rouquier~\cite{Rouquier:QuiverHecke2Lie}.  Following~\cite{Rostam:Grpn},
and to a lesser extent \cite{BK:GradedKL}, this section defines the
quiver Hecke algebras of type~$G(\ell,1,n)$ that we need to study
$\Hlpn$.

\begin{Definition}\label{D:ICompositions} Set
  $\cI=\set{\eps^j \q^i|0\leq j<p \text{ and } 0\leq i<e}$.
  An \textbf{$\cI$-composition} of $n$ is a finitely supported tuple
  $\alpha=(\alpha_i)_{i\in\cI}$ of non-negative integers that sum to~$n$.
  Let $\Comp$ be the set of $\cI$-compositions of~$n$.
  If $\alpha=(\alpha_i)_{i\in\cI}\in\Comp$ let
  \[
     \cIa=\set[\big]{\bi=(i_1,\dots,i_n)\in\cI^n|
     \alpha_i=\#\set{1\leq k\leq n|i_k=i} \text{ for all } i\in \cI}.
  \]
  A \textbf{residue sequence} is an element of $\cI^n$.
\end{Definition}

Let $\Sym_n$ be the \textbf{symmetric group} of degree $n$. As a Coxeter group,
$\Sym_n$ is generated by $s_1,\dots,s_{n-1}$, where $s_r=(r,r+1)$. If
$w\in\Sym_n$ then a \textbf{reduced expression} for $w$ is any word
$w=s_{r_1}\dots s_{r_k}$ with $k$ minimal.

The symmetric group $\Sym_n$ acts on $\cI^n$ by place permutations and
the sets $\cIa$, for $\alpha\in\Comp$, are the $\Sym_n$-orbits of $\cI^n$.
In particular, the sets $\cIa$ are finite. Set
\[
      J'=\set{1,2,\cdots,p'} \quad\text{and}\quad
      I= \Z/e\Z.
\]
As noted in \cite{Rostam:Grpn}, there is a natural bijection $I\times J'\bijection\cI$ given
by sending $(i,j)$ to $\eps^j\q^i$, for $(i,j)\in I\times J'$.
Henceforth, we identify $I\times J'$ and $\cI$ using this bijection.

\begin{Definition}[{Rouquier~\cite[\Sec3.2.5]{Rouquier:QuiverHecke2Lie}}]\label{D:Gamma}
  Let $\Gamma$ be the quiver with vertex set $\cI$ and edges
  $i\longrightarrow \q i$, for $i\in\cI$.
  Let $\Gamma_e$ be the full subquiver of $\Gamma$ with vertex set
  $\set{\q^i|i\in I}$.
\end{Definition}

Notice that there is an isomorphism of quivers
$\Gamma_e\cong\Gamma^{(j)}_e$, where $\Gamma_e^{(j)}$ has vertex set
$\set{\eps^j\q^i|i\in I}$, for $1\le j\le p'$. Moreover, there are
no edges between the vertices of $\Gamma^{(j)}_e$ and the vertices of
$\Gamma_e^{(k)}$ if $j\ne k\in J'$. Hence,
$\Gamma=\Gamma_e^{(1)}\sqcup\dots\sqcup\Gamma_e^{(p')}$ is the
disjoint union of~$p'$ copies of the quiver~$\Gamma_e$, which is the
affine quiver of type $A^{(1)}_e$.

Following Rouquier~\cite[\Sec 3.2.4]{Rouquier:QuiverHecke2Lie}, for
$i,j\in\cI$ define homogeneous polynomials $Q_{i,j}(u,v)\in
R[u,v]$, where $u$ and $v$ are indeterminates, by
\[
     Q_{i,j}(u,v) =   \begin{cases*}
           (u-v)(v-u), & if $i\leftrightarrows j$,\\
                (v-u), & if $i\rightarrow j$,\\
                (u-v), & if $i\leftarrow j$,\\
                    1, & if $i\noedge j$,\\
                    0, & if $i=j$.
     \end{cases*}
\]
where all edges are in the quiver $\Gamma$. The \textbf{degree} of
$Q_{i,j}(u,v)$ is its homogeneous degree.

For $\iota=(i,j)\in I\times J'=\cI$, define
$\Lambda_\iota$ to be the multiplicity of $\eps^j\q^i$ in $\bvQ$.
Set $\bLam=(\Lambda_\iota)_{\iota\in\cI}$.

\begin{Definition}[{Khovanov and
Lauda~\cite{KhovLaud:diagI},
Rouquier~\cite{Rouquier:QuiverHecke2Lie}}]
\label{D:QuiverHecke}
  Let $\alpha$ be an $\cI$-composition of~$n$.
  The \textbf{quiver Hecke algebra of type $G(\ell,1,n)$ and weight~$\bLam$} is
  the unital associative $R$-algebra
  $\RG[\alpha]$  with generators
  \[
       \set{e(\bi)|\bi\in\cIa}\cup\set{y_1,\dots,y_n}
            \cup\set{\psi_1,\dots,\psi_{n-1}}
  \]
  and relations
  { \setlength{\abovedisplayskip}{2pt}
    \setlength{\belowdisplayskip}{1pt}
    \begin{align*}
      e(\bi) e(\bj) &= \delta_{\bi\bj} e(\bi),
        &{\textstyle\sum_{\bi \in I^\alpha}} e(\bi)&= 1,
        & y_1^{\Lambda_{i_1}}e(\bi)&=0, \\
      y_r e(\bi) &= e(\bi) y_r,
      &\psi_r e(\bi)&= e(s_r{\cdot}\bi) \psi_r,
      &y_r y_s &= y_s y_r,
    \end{align*}
    \begin{align*}
      \psi_r \psi_s &= \psi_s \psi_r,&&\text{if }|r-s|>1,\\
      \psi_r y_s  &= y_s \psi_r,&&\text{if }s \neq r,r+1,
    \end{align*}
    \begin{align*}
        \psi_r y_{r+1} e(\bi)=(y_r\psi_r+\delta_{i_ri_{r+1}})e(\bi),
        \quad
        y_{r+1}\psi_re(\bi)=(\psi_r y_r+\delta_{i_ri_{r+1}})e(\bi),
    \end{align*}
    \begin{align*}
      \psi_r^2e(\bi) &= Q_{i_r,i_{r+1}}(y_r,y_{r+1})e(\bi)\\
    (\psi_{r+1}\psi_{r}\psi_{r+1}-\psi_{r}\psi_{r+1}\psi_{r})e(\bi)
    &=\delta_{i_ri_{r+2}}\frac{Q_{i_r,i_{r+1}}(y_r,y_{r+1})-Q_{i_{r},i_{r+1}}(y_{r+2},y_{r+1})}
           {y_r-y_{r+2}}e(\bi)
    \end{align*}
  }
  for all admissible $r,s$ and $\bi,\bj\in \cIa$. Set
  $\RG = \displaystyle\bigoplus_{\alpha\in\Comp} \RG[\alpha]$.
\end{Definition}

\begin{Remark}
  In the literature, $\RG$ is often called a cyclotomic Hecke algebra
  of type $A$ and weight~$\bLam$. Our naming convention reflects the
  close connections between the algebras $\Hln$ and $\RG$.
\end{Remark}

An important consequence of these relations is that $\RG$ is a $\Z$-graded
algebra with
\[
   \deg e(\bi)=0,\quad
   \deg y_r=2\quad\text{and}\quad
   \deg\psi_re(\bi) = \begin{cases} \deg Q_{i_r,i_{r+1}}(u,v), &\text{if $i_r\neq i_{r+1}$}\\
   -2, &\text{if $i_r=i_{r+1}$.}\end{cases}
\]
Following the reformulation in \cite{Rostam:Grpn}, we can now state the
main result of \cite{BK:GradedKL} that we need in order to apply this
result to the algebra~$\Hlpn$.

\begin{Theorem}[Brundan and Kleshchev's isomorphism
  theorem~\cite{BK:GradedKL, Rostam:Grpn}]\label{T:BKIsomorphism}
  Assume that $R = K$ is a field.
  Then there is an isomorphism of $K$-algebras
  $f\colon\Hln\bijection\RG$.
\end{Theorem}

Motivated in part by \autoref{D:Hlpn}, the third named author~\cite{Rostam:Grpn}
generalised this result to show that $\Hlpn$ is isomorphic to the
fixed-point subalgebra of~$\RG$ under a certain homogeneous automorphism of
order~$p$. Recall the automorphism $\sigma$ of $\Hln$ from
\autoref{E:sigma}. By definition,~$\sigma$ has order~$p$, so
$\<\sigma\>\cong\<\eps\>$, which is a cyclic group of order~$p$.

If $\alpha\in\Comp$ let $\sigma\cdot\alpha$ be the $\cI$-composition of
$n$ given by
\begin{equation}
\label{E:action_sigma_compositions}
     (\sigma\cdot\alpha)_i = \alpha_{\eps^{-1}i},\qquad
     \text{ for all }i\in \cI.
\end{equation}
Observe that left multiplication by $\eps$ gives a map
$\cIa\longrightarrow\cI^{\sigma\cdot\alpha};
\bi\mapsto\eps\bi=(\eps i_1,\dots,\eps i_n)$. Moreover, by~\autoref{E:bvQ}, $\Lambda_{\eps \iota} = \Lambda_{\iota}$ for all $\iota \in \cI$.

\begin{Theorem}[{Rostam~\cite{Rostam:Grpn},\cite[\Sec 1.4]{Rostam:PhD}}]
  \label{sigmaIso}
  Let $\alpha\in\Comp$. There is a unique homogeneous $R$-algebra
  isomorphism $\sigma^\bLam_\alpha\map{\RG[\alpha]}\RG[\sigma\cdot\alpha]$
  such that
  \[
       \sigma^\bLam_\alpha\bigl(e(\bi)\bigr) = e(\eps\bi),\quad
       \sigma^\bLam_\alpha(y_r)=y_r\quad\text{and}\quad
       \sigma^\bLam_\alpha(\psi_s)=\psi_s,
  \]
  for all $1\le r\le n$, $1\le s<n$ and $\bi\in \cI^\alpha$.
\end{Theorem}
%

  Set $\sigma^\bLam_n = \bigoplus_\alpha\sigma^\bLam_\alpha$, so that
  $\sigma^\bLam_n$ is an automorphism of
  $\RG=\bigoplus_\alpha\RG[\alpha]$.  To ease the notation, we normally
  write $\sigma=\sigma^\bLam_n$.  We are abusing notation here because
  the automorphism $\sigma$ of~$\RG$ is not equal to the
  automorphism~$\sigma$ of~$\Hln$ that was defined in \autoref{E:sigma}.
  This abuse is justified by \autoref{rpnIso} below.

If $\alpha$ is an $\cI$-composition let
$[\alpha]=\set{\sigma^k\cdot\alpha|1\le k\le p}$ be the  orbit of $\alpha$ under the action of $\langle \sigma\rangle \simeq \Z/p\Z$. Let
$\cIsig=\set{[\alpha]|\alpha\in\Comp}$ be the set of
$\sigma$-orbits of $\Comp$ and if  $[\alpha]\in\cIsig$
set $\RGa=\bigoplus_{\beta\in[\alpha]}\RG[\beta]$.
By definition,
\[
    \RG=\bigoplus_{\alpha\in\Comp}\RG[\alpha]
       =\bigoplus_{[\alpha]\in\cIsig}\RGa
\]
and the isomorphism $\sigma$ of \autoref{sigmaIso} restricts to an
automorphism $\sigma$ of $\RGa$.  Hence, we can consider $\sigma$
as both an automorphism of $\RGa$ and as an automorphism of~$\RG$.

\begin{Theorem}[{Rostam\cite[Theorem 4.14, Corollary 4.16]{Rostam:Grpn}}]
  \label{rpnIso}
  Assume that $R = K$ is a field.
  We can choose the isomorphism $f\colon\Hln\bijection\RG$ of
  \autoref{T:BKIsomorphism} so that the following diagram commutes
    \begin{center}
      \begin{tikzpicture}[>=stealth,->,shorten >=2pt,looseness=.5,auto]
        \matrix (M)[matrix of math nodes,row sep=1cm,column sep=16mm]{
          \Hln & \RG \\
          \Hln & \RG \\
         };
         \draw[->](M-1-1)--node[above]{$f$}(M-1-2);
         \draw[->](M-2-1)--node[below]{$f$}(M-2-2);
         \draw[->](M-1-1)--node[left]{$\sigma$}(M-2-1);
         \draw[->](M-1-2)--node[right]{$\sigma$}(M-2-2);
      \end{tikzpicture}
    \end{center}
  Consequently, $f$ induces an isomorphism
  $\displaystyle\Hlpn\cong\bigoplus_{[\alpha]\in\cIsig}\left(\RGa\right)^\sigma$.
\end{Theorem}

\begin{Definition}[{Rostam~\cite{Rostam:Grpn}, \cite[\Sec 1.4]{Rostam:PhD}}]\label{D:Rlpn}
  The \textbf{quiver Hecke algebra of type $G(\ell,p,n)$} of weight $\bLam$ is the $R$-algebra
  \[ \Rlpn = \left(\RG\right)^\sigma = \bigoplus_{[\alpha]\in\cIsig}\left(\RGa\right)^\sigma.\]
  For $[\alpha]\in\cIsig$ let $\Rlpn[\alpha] =\left(\RGa\right)^\sigma$.
\end{Definition}

The algebra $\Rlpn$ inherits a $\Z$-grading from $\RG$ since $\sigma$ is
a homogeneous automorphism of $\RG$.  The aim of this paper is to better
understand the algebra $\Rlpn$. Our main tool is the diagrammatic
Cherednik algebra introduced by
Webster~\cite{Webster:RouquierConjecture} and
Bowman~\cite{Bowman:ManyCellular}.

\subsection{Loadings, multicharges and \texorpdfstring{$\ell$}{l}-partitions}\label{S:configurations}

This section introduces the combinatorics that underpins Webster's
diagrammatic Cherednik algebras.

After \autoref{D:Hlpn} we fixed integers $(d, e, e', p, p', m, n, \rho_1,\dots,\rho_d)$
subject to \autoref{D:dCharge} that determine $\q$ and $\bvQ$.
Using this data we now fix a choice of multicharge
$\charge=(\kappa_1,\dots,\kappa_\ell)$ that we use in
\autoref{S:Webster} to single out a diagrammatic Cherednik algebra that
is particularly well adapted to studying $\Hlpn$.

\begin{Definition} \label{D:multicharge}
  The  \textbf{multicharge} of $\Hlpn$ is the sequence
  $\charge=(\kappa_1,\dots,\kappa_\ell)$ with
  \[
  \kappa_{l} = \rho_{a+1} + \frac{be'}{p'},
  \]
  where $l = ap + b+1$ with $0 \leq a < d$ and $0 \leq b < p$.
\end{Definition}

\begin{Remark}
If $d = 1$, the multicharge $\charge$ is an example of a ``FLOTW charge'' (see~\cite[Example 1.6]{Bowman:ManyCellular}).
\end{Remark}


For the rest of this paper we fix the multicharge $\charge$ of
\autoref{D:multicharge} and we identify~$\eps$ and~$q^{e'/p'}$.
Recalling \autoref{D:dCharge}, \autoref{D:multicharge} implies that
$\bvQ = (\eps^{\kappa_1},\dots, \eps^{\kappa_\ell})$.

Note that $b\tfrac{e'}{p'}<\tfrac{pe'}{p'}=me'=e$, for $0\le b<p$.
Therefore, $0=\kappa_1<\kappa_2<\dots<\kappa_\ell$ since $\rho_{a+1}-\rho_a\ge e$
by \autoref{D:dCharge}.

\begin{Examples}\leavevmode\newline
  \begin{itemize}
    \item\vskip-5mm Consider $G(p,p,n)$ and suppose that $\gcd(e,p)=1$. Then
    $\ell=p=p'$ and $m=d=1$ so that $e'=e$ and
    $\charge=(0,\frac ep,\frac{2e}{p}\dots,\frac{(p-1)e}p)$.
    \item Suppose that $p'=3$ and $m=2$ with $e=2e'<\infty$. Then
    \[
      \charge=\bigl(0,\tfrac{e'}3,\dots,\tfrac{5e'}{3},\rho_2,\dots,\rho_2+\tfrac{5e'}{3},
                  \dots, \rho_d,\dots,\rho_d+\tfrac{5e'}3\bigr).
    \]
  \end{itemize}
\end{Examples}

The set of \textbf{$\ell$-nodes} is the set of all ordered triples
\begin{equation}\label{nodes}
  \Nodes=\set{(r,c,l)\in\N^3|0\le r,c\leq n \text{ and } 1\leq l\leq\ell}.
\end{equation}


A \emph{partition} of $n$ is a sequence $\mu =(\mu_1,\dots,\mu_h)$ of non-negative integers satisfying $\mu_1\geq\dots\geq\mu_h$ and $|\mu|\coloneqq \mu_1+\dots+\mu_h = n$. Let $(0)$ be the empty partition (where $h = 0$) and use exponentiation for repeated parts.
An \emph{$\ell$-partition} of $n$ is an $\ell$-tuple
$\bmu=(\mu^{(1)}|\dots|\mu^{(\ell)})$ of partitions such that
$|\mu^{(1)}|+\dots+|\mu^{(\ell)}|=n$. Let $\Parts$ be the set of
$\ell$-partitions of~$n$.   An $\ell$-partition $\bmu\in\Parts$ is
identified with its \textbf{diagram}, which is the set of nodes
\[
     \bmu=\set{(r,c,l)\in\Nodes| 1\le r\le\mu^{(l)}_c}.
\]
We draw $\ell$-partitions as an array of boxes in plane using Bowman's
variation of the \textit{Russian convention}, as in the following
example.

\begin{Example}\label{Ex:RussianDiagrams}
  Let $\bmu=(3,1|2^2|1^3)$. The diagram of $\bmu$ is:
  \[
    \left(\space
      \RussianTableau{{,,},{\space}}\hspace*{2mm}\middle|\hspace*{2mm}
      \RussianTableau{{,},{,}}\hspace*{2mm}\middle|\hspace*{2mm}
      \RussianTableau{{\space},{\space},{\space}}
    \space\right)
  \]
\end{Example}

Fix an integer $N$ with
\begin{equation}
\label{E:inequality_N}
N>2nep'(\ell+1).
\end{equation}
Using the multicharge
$\charge=(\kappa_1, \dots, \kappa_\ell)$ define a \textbf{loading} function
$\xcoord\map\Nodes\Q$ by
\begin{equation}\label{E:xcoord}
  \xcoord(r, c, l) = c-r + \tfrac1e\bigl(\kappa_l-\tfrac{l-1}{\ell+1}\bigr)
                         - \tfrac{r+c}{N}
\end{equation}

\begin{Remark}\label{R:xcoord}
  Bowman~\cite[\Sec1.3]{Bowman:ManyCellular} defines his
  loading function as
  $\xcoord(r,c,l)=\kappa_l-\tfrac{l}{\ell}+c-r-\NN(r+c)$. The term
  $\frac{l-1}{\ell+1}$ is there to separate the nodes $(r,c,l)$ and
  $(r,c,l')$, when $l\ne l'$. We divide by $\ell+1$, rather than $\ell$
  as Bowman does, precisely because if $m=1=d$ then $p=p'=\ell$ so
  $\kappa_l=\rho_{a+1}+\tfrac{be'}{\ell}$. The $l-1$ in the numerator is a
  convenient renormalisation so that $\xcoord(0,0,1)=\kappa_1=0$.
\end{Remark}

\begin{Lemma}
\label{L:kappal_l-1_increasing}
The function $\set{1, \dots, \ell}\to\R; l \mapsto \kappa_l - \frac{l-1}{\ell+1}$ is strictly increasing. Moreover, if $l = ap  + b + 1$, with $0 \leq a < d$ and $0 \leq b < p$, then
\[
\rho_{a+1} - \frac{ap}{\ell + 1} \leq \kappa_l - \frac{l-1}{\ell+1} < \rho_{a+1} + e - \frac{(a+1)p}{\ell + 1}.
\]
\end{Lemma}

\begin{proof}
Let $l, l' \in \set{1, \dots, \ell}$ and write $l = ap + b+1$ and $l' =
a'p+b'+1$ with $0 \leq a, a' < d$ and $0 \leq b, b' < p$. Without loss
of generality, assume that $l < l'$. If $a < a'$ then by the inequality
of \autoref{D:multicharge}, together with the observation that
$p\frac{e'}{p'} = me' = e$, we obtain
\[
\kappa_{l'} - \kappa_l - \frac{l'-l}{\ell+1}= \rho_{a'+1} - \rho_{a+1} + (b'-b)\frac{e'}{p'} - \frac{l'-l}{\ell+1} \geq (2n+3)e -e -1> 0.
\]
Now if $a = a'$ and $b < b'$ we have
\[
\kappa_{l'} - \kappa_l - \frac{l'-l}{\ell+1} = (b'-b)\left(\frac{e'}{p'} - \frac{1}{\ell + 1} \right) = (b'-b)\left(\frac{e}{p} - \frac{1}{\ell+1}\right) > 0,
\]
since $p \leq \ell$, proving the first claim. If $l = ap + b+1$, we have $ap + 1 \leq l \leq (a+1)p$ and we deduce that
\[
\rho_{a+1} - \frac{ap}{\ell + 1} \leq \kappa_l - \frac{l-1}{\ell + 1} \leq \rho_{a+1} + e - \frac{e}{p}-\frac{(a+1)p-1}{\ell + 1}.
\]
Thus, we deduce the result since $\frac{1}{\ell+1} < \frac{e}{p}$.
\end{proof}

The key properties of the $\xcoord$-coordinate function are given by the
following lemma

\begin{Lemma}\label{L:xcoord}
  Let $\gamma=(r,c,l),\gamma'=(r',c',l')\in\Nodes$. Let $a, a' \in \set{0, \dots, d-1}$ such that $l - ap , l' - a'p \in \set{1, \dots, p}$.
  \begin{enumerate}
    \item If $\gamma\ne\gamma'$ then
    $\xcoord(\gamma)\notin \set{\xcoord(\gamma'), \xcoord(\gamma')\pm 1}$. 
    \item If $a>a'$ then $\xcoord(\gamma)>\xcoord(\gamma') + 1$.
    \item If $a = a'$ and $c - r > c' - r'$ then $\xcoord(\gamma) > \xcoord(\gamma')$.
  \end{enumerate}
\end{Lemma}

\begin{proof}
  First consider part~(a). Write $l=ap+b+1$, $l'=a'p+b'+1$, where  $0\le b,b'<p$.
Suppose that  $\xcoord(\gamma)= \xcoord(\gamma')\pm 1$. Then we have $$
\frac{1}{e}(\rho_{a+1}-\rho_{a'+1})+(c-r-c'+r')+(b-b')\frac{1}{p}-\frac{l-l'}{(\ell+1)e}-\frac{r+c-r'-c'}{N}=\pm 1 .
$$
Applying \ref{D:dCharge}, \ref{nodes} and \ref{E:inequality_N}, we can deduce from the above equality that $a=a'$ and $c-r=c'-r'\pm 1$. Since
$|\frac{r+c-r'-c'}{N}|<\frac{1}{ep'(\ell+1)}$, it follows that $b=b'$ and hence $c+r=c'+r'$, which is impossible because $c-r=c'-r'\pm 1$. This proves that
$\xcoord(\gamma)\neq\xcoord(\gamma')\pm 1$. In the case where $\pm 1$ is replaced by $0$, then a similar argument shows that $a=a'$, $b=b'$, $c-r=c'-r'$ and $c+r=c'+r'$. Thus $c=c'$ and $r=r'$, which contradicts the fact that $\gamma\ne\gamma'$. This proves that $\xcoord(\gamma)\neq\xcoord(\gamma')$.
This completes the proof of~(a).

  For part~(b), as above write $l'=a'p+b'+1$, where
  $0\le b'<p$. Recall that $\rho_{a+1}-\rho_{a'+1}\ge(2n+3)e$
  by \autoref{D:dCharge} since $a>a'$. Therefore,
  \begin{align*}
     \xcoord(\gamma)-\xcoord(\gamma')
       &=\mathrlap{(c-r-c'+r')+\tfrac1e\bigl(\rho_{a+1}-\rho_{a'+1}\bigr)
          +\tfrac{(b-b')e'}{ep'}
          -\tfrac{l-l'}{e(\ell+1)}
          +\tfrac{r'+c'-r-c}{N}}\\
       &\ge-2n+(2n+3)+\tfrac{1-p}e\tfrac{e'}{p'} -\tfrac{\ell-1}{e(\ell+1)}
          -\tfrac{2n}{N},
          &&\text{since $|b-b'|<p$,} \\
       &\ge3+ \tfrac{e'}{ep'} - \tfrac{pe'}{ep'} -\tfrac{\ell}{e(\ell+1)},
          &&\text{since $N\ge 2ne(\ell+1)$},\\
       &> 2-\tfrac1{p}-\tfrac1{e},
       &&\text{since $\tfrac{e}{e'}=m=\tfrac{p}{p'}$ and $\ell=pd$.}
  \end{align*}
  Hence, $\xcoord(\gamma)>\xcoord(\gamma') + 1$ since $p\ge2$ and $e\ge2$.

  Finally, part~(c) is immediate from the definition of $\xcoord$ because
  $0 \leq \kappa_l - \rho_{a+1} < p\frac{e'}{p'} = me' = e$.
\end{proof}

In particular, \autoref{L:xcoord}(a) shows that $\xcoord$ defines a total
order on the set of nodes.



Let $\blam\in\Parts$ be an $\ell$-partition. Abusing notation slightly,
the \textbf{loading} of $\blam$ is the set
\begin{equation}\label{E:lamLoading}
    \xcoord(\blam) = \set{\xcoord(r,c,l)|(r,c,l)\in\blam}.
\end{equation}

\begin{Example}
  Let $\blam=(2,1,1|2^2)$ and suppose that $p=1$ so that
  $\charge=\brho$. The following diagram shows the loadings $\xcoord(\blam)$ for two
  different choices of $\brho$.
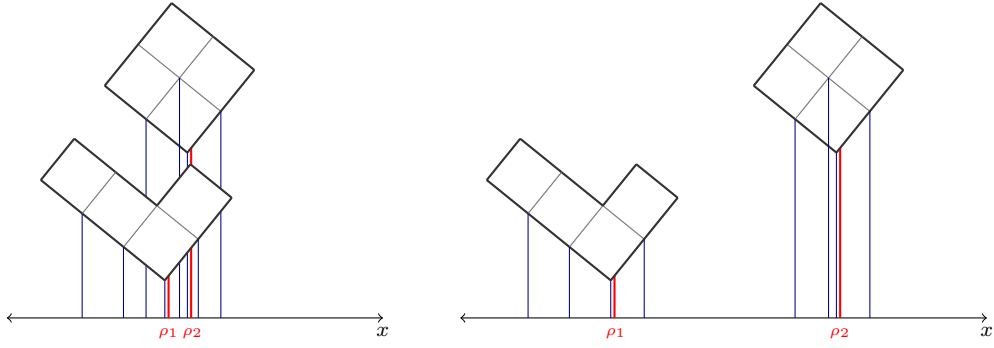
  \begin{center}
    \begin{tikzpicture}[every node/.append style={font=\footnotesize}]
       \coordinate (A) at (-1,-1);
       \coordinate (B) at (4,-1);
       \draw[<->, thin](A)--(B)node[below]{$x$};
       \draw[red, thick](1.15,-1)node[below]{$\scriptstyle\rho_1$} -- ++(0,0.56);
       \draw[red, thick](1.45,-1)node[below]{$\scriptstyle\,\rho_2$} -- ++(0,2.55);
       \pic (one) at (1.4,1.2) {RussianTableau={scale=0.5,rotate=-4}{{,},{,}}};
       \draw[semisolid](one-1-1) -- ($(A)!(one-1-1)!(B)$);
       \draw[semisolid](one-1-2) -- ($(A)!(one-1-2)!(B)$);
       \draw[semisolid](one-2-1) -- ($(A)!(one-2-1)!(B)$);
       \draw[semisolid](one-2-2) -- ($(A)!(one-2-2)!(B)$);
       \pic(two) at (1.1,-0.5) {RussianTableau={scale=0.5,rotate=-4}{{,},{\space},{\space}}};
       \draw[semisolid](two-1-1) -- ($(A)!(two-1-1)!(B)$);
       \draw[semisolid](two-1-2) -- ($(A)!(two-1-2)!(B)$);
       \draw[semisolid](two-2-1) -- ($(A)!(two-2-1)!(B)$);
       \draw[semisolid](two-3-1) -- ($(A)!(two-3-1)!(B)$);
     \end{tikzpicture}
     \qquad
    \begin{tikzpicture}[every node/.append style={font=\footnotesize}]
       \coordinate (A) at (-1,-1);
       \coordinate (B) at (6,-1);
       \draw[<->, thin](A)--(B)node[below]{$x$};
       \draw[red, thick](4.05,-1)node[below]{$\scriptstyle\rho_2$} -- ++(0,2.55);
       \draw[red, thick](1.05,-1)node[below]{$\scriptstyle\,\rho_1$} -- ++(0,0.56);
       \pic (one) at (4,1.2) {RussianTableau={scale=0.5,rotate=-4}{{,},{,}}};
       \draw[semisolid](one-1-1) -- ($(A)!(one-1-1)!(B)$);
       \draw[semisolid](one-1-2) -- ($(A)!(one-1-2)!(B)$);
       \draw[semisolid](one-2-1) -- ($(A)!(one-2-1)!(B)$);
       \draw[semisolid](one-2-2) -- ($(A)!(one-2-2)!(B)$);
       \pic(two) at (1,-0.5) {RussianTableau={scale=0.5,rotate=-4}{{,},{\space},{\space}}};
       \draw[semisolid](two-1-1) -- ($(A)!(two-1-1)!(B)$);
       \draw[semisolid](two-1-2) -- ($(A)!(two-1-2)!(B)$);
       \draw[semisolid](two-2-1) -- ($(A)!(two-2-1)!(B)$);
       \draw[semisolid](two-3-1) -- ($(A)!(two-3-1)!(B)$);
     \end{tikzpicture}
  \end{center}
  In both diagrams, the line from a node $(r,c,l)$ to
  the $x$-axis gives the loading $\xcoord(r,c,l)$.  The different
  components of $\blam$ are drawn with different heights to make it
  easier to distinguished between them.  The next section explains the
  significance of this diagram and the red strings.
\end{Example}

%

Extending this notation slightly, define a \emph{generalised partition} to be a finite subset $\gbconf\subseteq \R\times\R\times\set{1,\dots,\ell}$ such that
$\xcoord(\gbconf)$ has the same cardinality as $\gbconf$ and for any $1\leq l\leq\ell$,
\begin{equation}\label{E:gen-box-config}
  \xcoord(0,0,l) \notin \xcoord(\gbconf), \quad
  \xcoord(0,0,l)-1 \notin\xcoord(\gbconf)\quad\text{and}\quad
  x+1\notin\xcoord(\gbconf) \text{ for all }x\in\xcoord(\gbconf).
\end{equation}
By \autoref{L:xcoord}, if $\blam \in \Parts$ then $\blam$ is a
generalised partition. Conversely, if $\gbconf\notin\Parts$ is a
generalised partition then $\gbconf$ need not satisfy the
conclusions of \autoref{L:xcoord}. When we consider generalised partitions
below we will only be interested in the set~$\xcoord(\gbconf)$.

In order to define a partial order on $\Parts$, define
the \textbf{residue} of $(r,c,l)\in\Nodes$ to be
\begin{equation}\label{E:residues}
  \res(r,c,l)\coloneqq q^{\kappa_l+c-r}\in\cI 
\end{equation}
Recall that $\eps=q^{e'/p'}$, so  $\res(r,c,l)=\eps^bq^{\rho_{a+1}+c-r}$,
where $l=ap+b+1$ for $0\le a<d$ and $0\le b< p$.
%

If $\blam\subseteq\Nodes$ write
$\blam=\set{\gamma^\blam_1,\dots,\gamma^\blam_n}$ so that the nodes
$\gamma^\blam_1,\dots,\gamma^\blam_n$ are sorted by decreasing loading function, that is,
$\xcoord(\gamma^\blam_n)<\dots<\xcoord(\gamma^\blam_1)$.
 The \emph{residue sequence} of~$\blam$ is
\begin{equation}\label{E:ResidueSequence}
   \res(\blam) =\bigl(\res(\gamma^\blam_1), \res(\gamma^\blam_2),\dots,
                    \res(\gamma^\blam_n)\bigr)\in\cI^n.
\end{equation}

If $\alpha\in\Comp$ is an $\cI$-composition of~$n$ set
\[
  \Parts[\alpha]=\set{\blam\in\Parts|\res(\blam)\in\cI^\alpha}.
\]
We have the decomposition $\Parts=\bigsqcup_\alpha\Parts[\alpha]$  (disjoint union).


The following definition plays a key role in this paper. In particular,
it defines the partial order  that appears in our main result, which gives a
skew cellular basis for~$\Hlpn$.

\begin{Definition}[Webster's ordering {\cite[Definition 1.3, Proposition 1.4]{Bowman:ManyCellular}}]
  \label{D:DominanceOrder}
  Let $\blam,\bmu\in\Parts$ be two $\ell$-partitions.
  Then $\blam\tedom\bmu$ if there exists a bijection
  $\theta\map\blam\bmu$ such that
  \[
      \res\left(\theta(\gamma)\right)=\res(\gamma)
      \qquad\text{and}\qquad
      \xcoord\left(\theta(\gamma)\right)\le\xcoord(\gamma),
      \hspace*{20mm}\text{ for all }\gamma\in\blam.
  \]
  If $\blam\tedom\bmu$ and $\blam\ne\bmu$ write $\blam\tdom\bmu$.
\end{Definition}

\begin{Example}\label{Ex:dominance}
  By \autoref{D:dCharge}, $0=\rho_1<\dots< \rho_d$. Therefore,
  $\blam\tedom(1^n|0|\dots|0)$ whenever
  $\res(\blam)=\res(1^n|0|\dots|0)$ and $(0|\dots|0|n)\tedom\bmu$
  whenever $\res(\bmu)=\res(0|\dots|0|n)$.
\end{Example}

\subsection{Diagrammatic Cherednik algebras}\label{S:Webster}

Webster realises the quiver Hecke algebras of type~$G(\ell,1,n)$ as
idempotent subalgebras of his diagrammatic Cherednik
algebras~\cite{Webster:RouquierConjecture}. Following
Bowman~\cite{Bowman:ManyCellular}, we now recall these results, extending
them to the slightly more general quiver~$\Gamma$ as we go.  We start by
defining Webster diagrams.

A \textbf{string} in $\R^2$ is a diffeomorphism of the form
$[0,1]\longrightarrow\R^2; t\mapsto (\str(t),t)$. By definition, a
string is a smooth curve in~$\R^2$ with no loops. We sometimes identify a string with
the corresponding map $t \mapsto \str(t)$.  We regard a string as a
directed path from bottom ($t=0$) to top $(t=1)$.

Every string that we consider will be labelled by a \textbf{residue}
$i\in\cI$. An \textbf{$i$-string} is a string of residue~$i$.

A \textbf{crossing} of two strings is a point where they intersect.  A
\textbf{dot} on a string is a distinguished point in the image of the
string that is not on any crossing or on the start or end points of the string.
We will frequently refer to the following configuration of strings when
they occur in sufficiently small \textbf{local} neighbourhoods of
diagrams: 
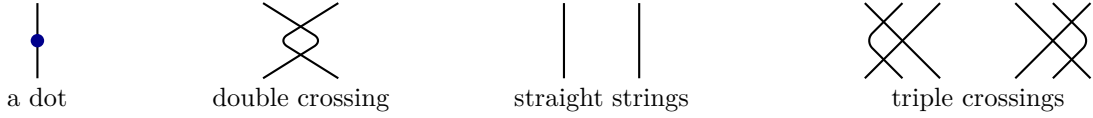
\begin{center} \begin{tikzpicture}[thick]
    \draw[dot](0,1)--++(0,-1)node[below]{a dot};
    \draw[rounded corners](3,0)--++(0.8,0.5)--++(-0.8,0.5);
    \draw[rounded corners](4,0)--++(-0.8,0.5)--++(0.8,0.5); \node at
    (3.5,0)[below]{double crossing};
    \draw[rounded corners](7,0)--++(0,1); \draw[rounded
    corners](8,0)--++(0,1); \node at (7.5,0)[below]{straight strings};
    \draw(11,0)--++(1,1); \draw(12,0)--++(-1,1); \draw[rounded
    corners](11.5,0)--++(-0.5,0.5)--++(0.5,0.5); \draw(13,0)--++(1,1);
    \draw(14,0)--++(-1,1); \draw[rounded
    corners](13.5,0)--++(0.5,0.5)--++(-0.5,0.5); \node at
    (12.5,0)[below]{triple crossings}; \end{tikzpicture} \end{center}
    Pulling \textbf{apart} the strings in a double crossing gives
    straight strings whereas pulling the string \textbf{through} the
    crossing in one of the triple crossings gives the other triple
    crossing. We apply this terminology below to red, solid and ghosts
    strings, which we now define.

Recall that \autoref{D:multicharge} fixes the multicharge $\charge\in \mathbb{Q}^\ell$.

\begin{Definition}[{Webster~\cite[Definition~4.1]{Webster:RouquierConjecture},
                    Bowman~\cite[Definition~4.1]{Bowman:ManyCellular}}]
  \leavevmode\newline
  Let $\blam,\bmu\subseteq\Nodes$. A \textbf{Webster diagram} with multicharge
  $\charge$ and type $(\blam,\bmu)$ and top residue sequence $\bi\in\cI^n$ consists of
  the following:
  \begin{enumerate}
    \item \textbf{Red} strings $\rstr_1,\dots,\rstr_\ell$ such that
    $\rstr_l$ has residue $q^{\kappa_l}$ and
    $\rstr_l(t)=\xcoord(0,0,l)$, for $t\in[0,1]$.
    \item \textbf{Solid} strings $\str_1,\dots,\str_n$, ordered so
    that $\str_1(1)>\str_2(1)>\dots>\str_n(1)$, such that
    \[
        \xcoord(\blam) = \set{\str_k(1)|1\le k\le n}
        \quad\text{and}\quad
        \xcoord(\bmu) = \set{\str_k(0)|1\le k\le n}
    \]
    and $\str_k$ is an $i_k$-string, for $1\le k\le n$.
    \item Each solid $i$-string has a \textbf{ghost} $i$-string that is
    obtained by translating the corresponding solid string one unit
    to the right.
  \end{enumerate}
  The solid strings in a Webster diagram are decorated with finitely
  many \textbf{dots} on the solid strings, with each dot having a
  \textbf{ghost dot} one unit to the right on the corresponding ghost
  string. Exactly two strings in a Webster diagram intersect at each
  crossing and no (red, solid or ghost) string can be tangential to any
  other string.

  Given a Webster diagram $D$, set $\top(D)=(\blam,\bi)$ and $\bot(D)=(\bmu,\bj)$, where
  $\bj=(j_1,\dots,j_n)\in\cI^n$ is the residue sequence of the solid strings
  when read in order from right to left along the bottom of~$D$.
  Then $\bi$ and $\bj$ are the \textbf{top residue sequence} and
  \textbf{bottom residue sequence} of $D$, respectively. If $\bi=\bj$ then
  $\res(D)=\bi=\bj$ is the \textbf{residue sequence} of $D$.
\end{Definition}

To help distinguish between the different types of strings in Webster
diagram we draw red strings as thick red strings and ghost strings as
dashed gray strings.

\begin{Remark}
  Ghost dots do not appear in Webster's paper
  \cite{Webster:RouquierConjecture} but can be found in
  Bowman~\cite[Remark 4.8]{Bowman:ManyCellular}. Including the ghost dots does not change the algebras up to isomorphism and makes the relations easier to write because they are more symmetrical with respect to the dots and ghost dots.
\end{Remark}

Two Webster diagrams of type $(\blam,\bmu)$ are equivalent if they have
the same residues, same number of dots on each string, when ordered from
right to left at the top of the diagram, and they differ by an isotopy,
which is a continuous deformation in which all of the intermediate
diagrams are Webster diagrams. In particular, the red strings are fixed
by isotopy.

Let $\Web(\blam,\bmu)$ be the set of (isotopy classes) of Webster
diagrams of type $(\blam,\bmu)$. If $\alpha$ is an $\cI$-composition
then set
\begin{equation}\label{E:WebsterDiagrams}
  \Web(\alpha) = \bigcup_{\blam,\bmu\in\Parts[\alpha]} \Web(\blam,\bmu).
  \quad\text{and}\quad
  \Web(n)=\bigcup_{\blam,\bmu\in\Parts}\Web(\blam,\bmu).
\end{equation}
Of course, $\Web(n) = \bigcup_{\alpha\in\Comp}\Web(\alpha)
                    = \bigcup_{\blam,\bmu\in\Parts}\Web(\blam,\bmu)$.

Let  $D, E\in\Web(n)$ be Webster diagrams such that $\bot(D)=\top(E)$.
Define $D\circ E$ to be the Webster diagram obtained by identifying the
southern points of~$D$ with the northern points of~$E$ and then
rescaling.

There is a distinguished Webster diagram
$\1_\blam^\bi\in\Web(\blam,\blam)$, for each $\ell$-partition
$\blam\in\Parts$ and each residue sequence $\bi\in\cI^n$, in which all of
the strings are vertical. By definition,
$\1_\blam^\bi\circ\1_\blam^\bi=\1_\blam^\bi$. If a Webster diagram $D$
of type $(\blam,\bmu)$ has top residue sequence $\bi$ and bottom residue
sequence $\bi'$ then $D=\1_\blam^\bi D\1_{\bmu}^{\bi'}$.

\begin{Example}\label{Ex:WebsterDiagram}
  Let $\ell=4$, $e=3$ and $p=2$ so that $e'=3$, $p'=2$ and let $\rho_2=9$
  so that $\charge=(0,1.5,9,10.5)$.
  Let $\blam=(2^2|1|0|2)$ and fix $\bi\in\cI^n$. Then $\1_\blam^\bi$ is
  the Webster diagram
  \[
     {\1}_{\blam}^{\bi}\quad=\quad \WebsterIdempotent3(0,1.5,9,10.5){{2,2},{1},{},{2}}
  \]
  where the solid strings have residues $i_1,\dots,i_7$ when read from
  right to left.  The $\xcoord$--coordinates of the solid strings are
  given by the loadings $\xcoord(\gamma)$, for $\gamma\in\blam$. By
  \autoref{E:xcoord}, the $\xcoord$--coordinate of the $l$th red string is
  $\frac1e\bigl(\kappa_l-\frac{l-1}{\ell+1}\bigr)$, for $1\le l\le\ell=4$.
  In particular, the leftmost red string has $x$-coordinate
  $\xcoord(0,0,1)=\frac{\kappa_1}{3}=0$.
\end{Example}

\begin{Lemma}
\label{L:no_crossing_diagram_is_idempotent}
Let $\gbconf$ be a generalised partition and $\bi \in \cI^n$.
Let $D$ be a Webster diagram of type $(\gbconf,\gbconf)$ that does not contain any crossings.
Then $D$ is isotopic to $\1_{\gbconf}^\bi$, where $\bi$ is the (top) residue sequence of $D$.
\end{Lemma}

\begin{proof}
  In order to construct an isotopy from $D$ to $\1_{\gbconf}^\bi$,
  write $\xcoord(\gbconf)=\set{x_1<\dots<x_n}$ and let
  $\str_1,\dots,\str_n$ be the solid strings in~$D$.
  For $u \in [0,1]$, let $D^{(u)}$ be the Webster
diagram of type $(\gbconf,\gbconf)$ and residue
$\bi$, which has solid strings $\str_1^{(u)},\dots,\str_{n}^{(u)}$ given by
\[
\str^{(u)}_k\map{[0,1]}\R\times[0,1];t\mapsto\bigl((1-u)\str_k(t)+ux_k,t\bigr),\qquad
\text{for }1\le k\le n.
\]
By construction, $D=D^{(0)}$ and $\1_{\gbconf}^\bi=D^{(1)}$, so to
complete the proof it suffices to prove that the strings in~$D^{(u)}$
never intersect, for $u\in[0,1]$. By assumption, the solid strings in
$D$ do not intersect, so $\str_k(t) < \str_l(t)$ for $1\le k<l\le\ell$.
Therefore, if $u\in[0,1]$ then $\str^{(u)}_k(t) < \str^{(u)}_l(t)$ for
$1\le k<l\le\ell$, so the solid strings in $D^{(u)}$ do not intersect.
Essentially the same argument show that there are no intersections
between any of the solid, ghost and red strings in $D^{(u)}$, completing
the proof.
\end{proof}

We can now define Webster's diagrammatic Cherednik algebras.

\begin{Definition}[{Webster~\cite[Definition~4.2]{Webster:RouquierConjecture},
                    Bowman~\cite[Definition~4.5]{Bowman:ManyCellular}}]
\label{D:RationalCherednik}
 The
\textbf{diagrammatic Cherednik algebra} $\WA$ is the $R$-algebra
\[\WA=\bigoplus_{\alpha\in\Comp}\WA[\alpha],\]
where, for each $\cI$-composition $\alpha$, the $R$-algebra $\WA[\alpha]$
is the unital associative algebra generated by the Webster diagrams
in~$\Web(\alpha)$ such that
\[
        DE =   \begin{cases*}
                    D\circ E,& if $\bot(D)=\top(E)$,\\
                    0, & otherwise,
               \end{cases*}
\]
and the following bilocal relations hold:
\begin{enumerate}[label=\Alph*\upshape)]
  \Item(Dots and crossings) Solid and ghost dots can pass through any crossing except:
    \begin{center}
      \Relation{R:DotCrossing}
      \begin{tikzpicture}
          \draw[solid](1,1)--(0,0) node[below]{$i$};
          \draw[solid,dot=0.25](0,1)--(1,0) node[below]{$i$};
          \node at (1.5,0.5){$-$};
          \draw[solid](3,1)--(2,0) node[below]{$i$};
          \draw[solid,dot=0.75](2,1)--(3,0) node[below]{$i$};
          \node at (4.0,0.5){$=$};
          \draw[solid](4.5,1)--(4.5,0) node[below]{$i$};
          \draw[solid](5.5,1)--(5.5,0) node[below]{$i$};
          \node at (6.0,0.5){$=$};
          \draw[solid,dot=0.75](8,1)--(7,0) node[below]{$i$};
          \draw[solid](7,1)--(8,0) node[below]{$i$};
          \node at (8.5,0.5){$-$};
          \draw[solid, dot=0.25](10,1)--(9,0) node[below]{$i$};
          \draw[solid](9,1)--(10,0) node[below]{$i$};
        \end{tikzpicture}
    \end{center}

  \Item(Double crossings)
  \label{I:double_cross} A double crossing between any
  two strings can be pulled apart except in the following cases:
  \begin{center}
    \Relation{R:SolidSolid}
    \begin{tikzpicture}[rounded corners]
      \draw[solid](1.5,1) --++(0.8,-0.5) -- ++(-0.8,-0.5) node[below]{$i$};
      \draw[solid](2.5,1) --++(-0.8,-0.5) --++(0.8,-0.5) node[below]{$i$};
      \node at (3.0,0.5){$=0$};
    \end{tikzpicture},
    \qquad
    \Relation{R:RedSolid}
    \begin{tikzpicture}[rounded corners]
      \draw[redstring](0,1)--(0,0)node[below]{$i$};
      \draw[solid](0.5,1)--(-0.5,0.5)--(0.5,0)node[below]{$i$};
      \node at (1.0,0.5){$=$};
      \draw[redstring](1.5,0)node[below]{$i$}--++(0,1);
      \draw[solid, dot](2.0,0)node[below]{$i$}--++(0,1);
    \end{tikzpicture},
    \qquad
    \Relation{R:SolidRed}
    \begin{tikzpicture}[rounded corners]
      \draw[redstring](0,1)--(0,0)node[below]{$i$};
      \draw[solid](-0.5,1)--(0.5,0.5)--(-0.5,0)node[below]{$i$};
      \node at (1,0.5){$=$};
      \draw[redstring](2.0,0)node[below]{$i$}--++(0,1);
      \draw[solid, dot](1.5,0)node[below]{$i$}--++(0,1);
    \end{tikzpicture}
    \\[2mm]
    \Relation{R:GhostSolid}
    \begin{tikzpicture}[rounded corners]
      \draw[ghost](0,1) --++(0.8,-0.5) --++(-0.8,-0.5) node[below]{$i$};
      \draw[solid](1,1) --++(-0.8,-0.5) --++(0.8,-0.5) node[below]{$j$};
      \node at (1.5,0.5){$=$};
      \draw[ghost](2.0,1)--++(0,-1)node[below]{$i$};
      \draw[solid,dot](3,1)--++(0,-1)node[below]{$j$};
      \node at (3.5,0.5){$-$};
      \draw[ghost,ghostdot](4,1)--++(0,-1)node[below]{$i$};
      \draw[solid](5,1)--++(0,-1)node[below]{$j$};
    \end{tikzpicture}
    \quad and\quad
    \Relation{R:SolidGhost}
    \begin{tikzpicture}[rounded corners]
      \draw[solid](0,1) --++(0.8,-0.5) --++(-0.8,-0.5) node[below]{$j$};
      \draw[ghost](1,1) --++(-0.8,-0.5) --++(0.8,-0.5) node[below]{$i$};
      \node at (1.5,0.5){$=$};
      \draw[solid,dot](2.0,1)--++(0,-1)node[below]{$j$};
      \draw[ghost](3,1)--++(0,-1)node[below]{$i$};
      \node at (3.5,0.5){$-$};
      \draw[solid](4,1)--++(0,-1)node[below]{$j$};
      \draw[ghost,ghostdot](5,1)--++(0,-1)node[below]{$i$};
     \end{tikzpicture}
   \end{center}
   where $j=i\q\in\cI$
  \Item(Triple crossings)
  \label{I:triple_crossings}
  A string can be pulled through
  a  crossing  except in the cases:
  \begin{center}
    \Relation{R:BraidGSG}
     \begin{tikzpicture}
       \draw[ghost](1,1)--++(1,-1)node[below]{$i$};
       \draw[ghost](2,1)--++(-1,-1)node[below]{$i$};
       \draw[solid,rounded corners](1.5,1)--++(-0.5,-0.5)--++(0.5,-0.5)node[below]{$j$};
       \node at (2.5,0.5){$=$};
       \draw[ghost](3,1)--++(1,-1)node[below]{$i$};
       \draw[ghost](4,1)--++(-1,-1)node[below]{$i$};
       \draw[solid,rounded corners](3.5,1)--++(0.5,-0.5)--++(-0.5,-0.5)node[below]{$j$};
       \node at (4.5,0.5){$+$};
       \draw[ghost](5,1)--++(0,-1)node[below]{$i$};
       \draw[ghost](6,1)--++(0,-1)node[below]{$i$};
       \draw[solid](5.5,1)--++(0,-1)node[below]{$j$};
     \end{tikzpicture}
     \qquad
     \Relation{R:BraidSGS}
     \begin{tikzpicture}
       \draw[solid](1,1)--++(1,-1)node[below]{$j$};
       \draw[solid](2,1)--++(-1,-1)node[below]{$j$};
       \draw[ghost,rounded corners](1.5,1)--++(-0.5,-0.5)--++(0.5,-0.5)node[below]{$i$};
       \node at (2.5,0.5){$=$};
       \draw[solid](3,1)--++(1,-1)node[below]{$j$};
       \draw[solid](4,1)--++(-1,-1)node[below]{$j$};
       \draw[ghost,rounded corners](3.5,1)--++(0.5,-0.5)--++(-0.5,-0.5)node[below]{$i$};
       \node at (4.5,0.5){$-$};
       \draw[solid](5,1)--++(0,-1)node[below]{$j$};
       \draw[solid](6,1)--++(0,-1)node[below]{$j$};
       \draw[ghost](5.5,1)--++(0,-1)node[below]{$i$};
     \end{tikzpicture}
     \\[2mm]
     \Relation{R:BraidSRS}
     \begin{tikzpicture}
       \draw[solid,rounded corners](1,1)--(1.8,0.7)--(2,0)node[below]{$i$};
       \draw[solid,rounded corners](2,1)--(1.8,0.3)--(1,0)node[below]{$i$};
       \draw[redstring](1.5,1)--(1.5,0)node[below]{$i$};
       \node at (2.5,0.5){$=$};
       \draw[solid,rounded corners](3,0)node[below]{$i$}--(3.2,0.7)--(4,1);
       \draw[solid,rounded corners](4,0)node[below]{$i$}--(3.2,0.3)--(3,1);
       \draw[redstring](3.5,1)--(3.5,0)node[below]{$i$};
       \node at (4.5,0.5){$-$};
       \draw[solid](5,1)--++(0,-1)node[below]{$i$};
       \draw[solid](6,1)--++(0,-1)node[below]{$i$};
       \draw[redstring](5.5,1)--++(0,-1)node[below]{$i$};
     \end{tikzpicture}
  \end{center}
   where $j=i\q\in\cI$

  \Item(Unsteady diagrams)
    A Webster diagram is \textbf{unsteady} if it contains a solid string
    that at any point is $n$ units or more to the right of the rightmost
    red string. Any unsteady diagram is zero.
\end{enumerate}
\end{Definition}

Solid and ghost strings always occur in pairs, so any solid or ghost
strings that are not drawn in the relations above are still part of the
relations even though they do not appear. All of the relations in
\autoref{D:RationalCherednik} are \textit{bilocal} in the sense that the
relations need to be applied locally in the regions around the solid
strings and their ghost strings. In particular, strings may appear between the
solid strings and their ghosts in the double and triple crossing
relations.

The relations drawn in (A), (B) and (C) of \autoref{D:RationalCherednik},
are the \emph{exceptional relations}. The remaining relations are the
\emph{non-exceptional} relations of $\WA$.  When they are applied, none
of the non-exceptional relations introduce additional diagrams. Explicitly, the
non-exceptional relations in (A) allow a dot to be pulled through a
crossing, those in (B) allow a double crossing to be pulled apart, and
those in~(C) allow a string to be pulled through a triple crossing.

As with $\RG$, the algebra $\WA$ is $\Z$-graded with the grading defined
on the Webster diagrams by summing over the contributions from each dot
and crossing in the diagram according to the following rules:
\begin{align*}
  \deg \tikz{\draw[solid,dot](0,1)--++(0,-1)node[below]{$\phantom{j}$};} &= 2 &
  \deg\tikz{\draw[solid](1,1)--++(1,-1)node[below]{$j$};
            \draw[solid](2,1)--++(-1,-1)node[below]{$i$};} &=-2\delta_{i,j} &
  \deg\tikz{\draw[solid](1,1)--++(1,-1)node[below]{$j$};
            \draw[ghost](2,1)--++(-1,-1)node[below]{$i$};} &=\delta_{j,qi}\\
  \deg\tikz{\draw[ghost](1,1)--++(1,-1)node[below]{$j$};
            \draw[solid](2,1)--++(-1,-1)node[below]{$i$};} &=\delta_{j,q^{-1}i} &
  \deg\tikz{\draw[solid](1,1)--++(1,-1)node[below]{$j$};
            \draw[redstring](1.5,1)--++(0,-1)node[below]{$i$};} &=\delta_{i,j} &
  \deg\tikz{\draw[redstring](1.5,1)--++(0,-1)node[below]{$j$};
            \draw[solid,rounded corners](2,1)--++(-1,-1)node[below]{$i$};} &=\delta_{j,i}.
\end{align*}
All other crossings, and the ghost dots, have degree $0$. The algebra $\WA[\alpha]$ is
$\Z$-graded because all of the relations in
\autoref{D:RationalCherednik} are homogeneous with respect to this
degree function.

\begin{Remark} \label{R:extended_start_end_sets}
To make some proofs easier to read, we sometime require the diagrams to
have their solid strings starting or ending in a set $\xcoord(\gbconf)$,
where $\gbconf$ is a generalised partition (cf.
\autoref{E:gen-box-config}).  In particular, if $D \in \Web(\blam,
\bmu)$ is a Webster diagram with no dots or intersecting strings on the
line $y=H$, where $H \in (0,1)$, then we can factor $D$ as $D = D^+
\circ D^-$ where the diagrams $D^+$ and $D^-$ are the restrictions of
$D$  to $\R\times[H, 1]$ and $\R\times[0, H]$, respectively.
We then have $D^+ \in \Web(\blam,\gbconf)$
and $D^- \in \Web(\gbconf,\bmu)$ where $\gbconf$ is a generalised partition.
\end{Remark}

Following \cite{Bowman:ManyCellular} we now describe a basis of $\WA$.
Recall that we identify an $\ell$-partition $\blam\in\Parts$ with its
diagram. Let $\label{E:omega}\bom=(1^n|0|\dots|0)\in\Parts$; compare with
\autoref{Ex:dominance}. The $\ell$-partition $\bom$ is the unique
$\ell$-partition such that $\xcoord(\gamma)<0$, for all $\gamma\in\bom$.

\begin{Definition}\label{D:standard}
  Let $\blam\in\Parts$. A \textbf{$\blam$-tableau} is a bijection
  $\t\map\blam \set{1, \dots, n}$.  If $\t$ is a $\blam$-tableau then $\t$
  has \textbf{shape}~$\blam$ and we write $\Shape(\t)=\blam$. A tableau
  $\t$ is \textbf{standard} if its entries increase along the rows and
  columns. In other words,
  \begin{enumerate}
  \item If $(r,c,l),(r-1,c,l)\in\blam$ then $\t(r,c,l)\ge\t(r-1,c,l)+1$;
  \item If $(r,c,l),(r,c-1,l)\in\blam$ then $\t(r,c,l)\ge\t(r,c-1,l)+1$.
  \end{enumerate}
  Let $\Std(\blam)$ be the set of standard $\blam$-tableaux.
\end{Definition}

As in \autoref{Ex:RussianDiagrams}, think of standard tableaux as
labelled Russian diagrams.

\begin{Example}\label{Ex:RussianTableau}
  Let $\bmu=(3,1|2^2|1^3)$. Then one tableau in $\Std(\bmu)$ is
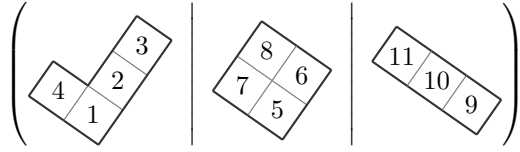
  \[
    \left(\space\RussianTableau{{1,2,3},{4}}\hspace*{2mm}\middle|\hspace*{2mm}
                \RussianTableau{{5,6},{7,8}}\hspace*{2mm}\middle|\hspace*{2mm}
                \RussianTableau{{9},{10},{11}}
    \space\right)
  \]
\end{Example}

Let $\t\in\Std(\blam)$ and $1\le k\le n$. The \textbf{residue} of $k$ in~$\t$ is
$\res_k(\t)=\res(\t^{-1}(k))$ and
\[
    \res(\t) = \bigl(\res_1(\t),\res_2(\t),\dots,\res_n(\t)\bigr)\in\cI^n
\]
is the \textbf{residue sequence} of~$\t$.

Let $\blam\in\Parts$. A node $\alpha\notin\blam$ is an \textbf{addable}
node of $\blam$ if $\blam\cup\set{\alpha}$ is (the diagram of) an
$\ell$-partition. Similarly, $\alpha\in\blam$ is a \textbf{removable}
node of $\blam$ if $\blam\setminus\set{\alpha}$ is an $\ell$-partition.
Let $\Add(\blam)$ and $\Rem(\blam)$ be the sets of addable and removable
nodes of~$\blam$.

Let $\t\in\Std(\blam)$.
If $1\le l\le n$ let $\t_{\downarrow k}$ be the restriction of $\t$ to
$\set{1,\dots,k}$ and let $\blam_k=\Shape(\t_{\downarrow k})$. Since
$\t$ is standard, $\t_{\downarrow k}$ is a standard $\blam_k$-tableau.
Define
\begin{align*}
 \Add_k(\t) &= \set{\gamma\in\Add(\blam_k)|\res(\gamma)=\res_\t(k)\text{ and }
    \xcoord(\t^{-1}(k))>\xcoord(\gamma)},\\
 \Rem_k(\t) &= \set{\gamma\in\Rem(\blam_k)|\res(\gamma)=\res_\t(k)\text{ and }
    \xcoord(\t^{-1}(k))>\xcoord(\gamma)}.
\end{align*}
Following \cite[Definition~1.11]{Bowman:ManyCellular}, and
\cite[(3.5)]{BKW:GradedSpecht}, the \textbf{degree} of $\t$ is the integer
\begin{equation}\label{E:TableauDegree}
    \deg\t = \sum_{k=1}^n \bigl(\#\Add_k(\t)-\#\Rem_k(\t)\bigr).
\end{equation}

In order to attach a Webster diagram to a standard tableau $\s$ let
$\cross(D)$ be the number of crossings in any diagram $D \in \Web(n)$.
The number $\cross(D)$ is preserved by isotopy but when we apply the
relations in $\WA$ diagrams with a different number of crossings can appear.

\begin{Definition}
[{\cite[\Sec4.3]{Webster:RouquierConjecture},
                    \cite[Definition 6.1]{Bowman:ManyCellular}}] \label{D:CTdiag}
  Let $\t\in\Std(\blam)$, for $\blam\in\Parts$. Let
  $C_\t\in\WA$ be any Webster diagram in $\Web(\blam,\bom)$ such that:
  \begin{enumerate}
    \item \label{I:CT_residue}
    For each $\gamma\in\blam$, there is a solid string of residue
    $\res(\gamma)$ from $(\xcoord[\t](\gamma),0)$ to
    $(\xcoord(\gamma),1)$, where
    $\xcoord[\t](\gamma)=\xcoord(\t(\gamma),1,1)$ and
    $(\t(\gamma),1,1)\in\bom$.
    \item \label{I:CT_no_dots}
    The diagram $C_\t$ has no dots on any strings.
    \item \label{I:CT_minimal_number_crossings}
    If $C_\t'$ is another diagram satisfying (a) and then
    $\cross(C_\t)\le\cross(C_\t')$.
  \end{enumerate}
\end{Definition}

A \emph{generalised double crossing} in a Webster diagram $D$ is a pair of strings in $D$ that cross twice. In particular, if a diagram $C_\t$ satisfies part (a) of \autoref{D:CTdiag} and has no generalised double crossing then $C_\t$ satisfies (c).

In general, the diagram $C_\t$ is not uniquely determined by
\autoref{D:CTdiag}.  In \autoref{subsection:particular_class_diagrams}  we give an
explicit construction of such diagrams but, for now, we let $C_\t$ be any
Webster diagram satisfying \autoref{D:CTdiag}. Unless stated otherwise, the
results that follow do not depend on the choice of diagram for $C_\t$.

By construction the bottom residue sequence of $C_\t$ is $\res(\t)$.
Moreover, by \cite[Theorem~7.1]{Bowman:ManyCellular},
$\deg C_\t=\deg \t$.

Let $\ast\map{\Web(\blam,\bmu)}\Web(\bmu,\blam)$ be the map that
reflects a Webster diagram in the line $y=\frac12$. Using
\autoref{D:RationalCherednik} it is easy to see that $\ast$ extends to
an involution $\ast\map{\WA[\alpha]}{\WA[\alpha]}$, for any
$\cI$-composition~$\alpha$. Hence, we can consider $\ast$ as a
homogeneous automorphism of~$\WA$ of order~$2$.

\begin{Definition}\label{D:CST}
  Let $\s, \t \in\Std(\blam)$.
  Define $C_{\s\t} = C_\s^*C_\t\in\WA$.
\end{Definition}

By the remarks above, $C_{\s\t}$ is homogeneous of degree
$\deg\s+\deg\t$.
The next result shows that a certain idempotent truncation of $\WA$, that will turn out to be isomorphic
to~$\RG[\alpha]$, is a graded cellular algebra in the
sense of Graham and Lehrer~\cite{GL,HuMathas:GradedCellular}.
Recall the idempotents $\1_\blam^\bi$ from before
\autoref{Ex:WebsterDiagram}, where $\blam\in\Parts$ and $\bi\in\cI^n$.
Let $\alpha$ be an $\cI$-composition of~$n$. Define the idempotents
\begin{align}\label{Eoma}
   \Eoma &=\sum_{\bi\in\cIa}\1_{\bom}^{\bi},
&
   \Eom &=\sum_{\alpha\in\Comp}\Eoma,
\\
\intertext{and define the algebras}\label{E:WAbom}
     \WA[\alpha](\bom) &= \Eoma\WA[\alpha]\Eoma,
&
     \WA(\bom) &= \Eom\WA\Eom = \bigoplus_{\alpha\in\Comp}\WA[\alpha](\bom).
\end{align}
We can now state one of the main results of
\cite{Webster:RouquierConjecture,Bowman:ManyCellular}.

\begin{Theorem}[{\cite[Theorem~4.11]{Webster:RouquierConjecture},
                 \cite[Theorem~7.1]{Bowman:ManyCellular}}]
  \label{T:CSTBasis}
  Let $\alpha$ be an $\cI$-composition of~$n$. Then the algebra
  $\WA[\alpha](\bom)$ is a graded cellular algebra with graded cellular basis
  \[
    \set[\big]{C_{\s\t}|\s, \t \in\Std(\blam) \text{ for }\blam\in\Parts[\alpha]},
  \]
  with respect to the poset $(\Parts[\alpha], \tedom)$ and homogeneous
  cellular algebra anti-isomorphism~$\ast$.
\end{Theorem}

In fact, Bowman and Webster give a cellular basis for the algebra
$\WA[\alpha]$. We state only this special case of their result because
this is all that we need and because it saves us from having to
introduce additional notation.

\begin{Remark}
  In type~$A$, Bowman~\cite{Bowman:ManyCellular} and
  Webster~\cite{Webster:RouquierConjecture, Webster:WeightedKLR}
  only consider algebras that are attached to the cyclic
  quiver~$\Gamma_e$ whereas we are considering the more general
  quiver~$\Gamma$ from \autoref{D:Gamma}. As the relations in $\WA$ are
  local and depend only on the quiver and the choice of residues, it is easy to
  see that the arguments of these papers apply without change for the
  quiver~$\Gamma$. The Webster diagrams, and hence the algebras $\WA$
  also depend on the choice of loading function. In his papers Webster
  considers arbitrary loadings whereas Bowman fixes a loading that is
  different from ours; see \autoref{R:xcoord}. If $\mathsf{x}$ and
  $\mathsf{y}$ are two loading functions then it is easy to see that
  the corresponding Webster algebras $\mathbb{A}^{\mathbf{x}}_n$ and
  $\mathbb{A}^{\mathbf{y}}_n$ are isomorphic if
  \[
     \mathsf{x}(\gamma)<\mathsf{x}(\gamma')\quad\text{if and only if}\quad
     \mathsf{y}(\gamma)<\mathsf{y}(\gamma')\qquad\text{for all } \gamma,\gamma'\in\Nodes,
  \]
  where an isomorphism is given by conjugating by the ``straight line''
  diagrams that have strings connecting $\mathsf{x}(\gamma)$ to
  $\mathsf{y}(\gamma)$, for all $\gamma\in\Nodes$.  Comparing the
  definition of our loading function $\xcoord$ from~\autoref{E:xcoord}
  with Bowman's (see \autoref{R:xcoord}), it is straightforward to see
  that our loading function is equivalent to one of Bowman's in the
  sense that they give isomorphic Webster algebras.
\end{Remark}

\begin{Theorem}[{\cite{Webster:RouquierConjecture},
                 \cite[Theorem 6.17]{Bowman:ManyCellular}}]
                 \label{ISO1}
  Let $\alpha$ be an $\cI$-composition of~$n$.
  There is a unique isomorphism of $\Z$-graded $R$-algebras
  $g\colon\RG[\alpha]\bijection \WA[\alpha](\bom)$ such that
  \begin{align*}
   e(\bi) & \mapsto
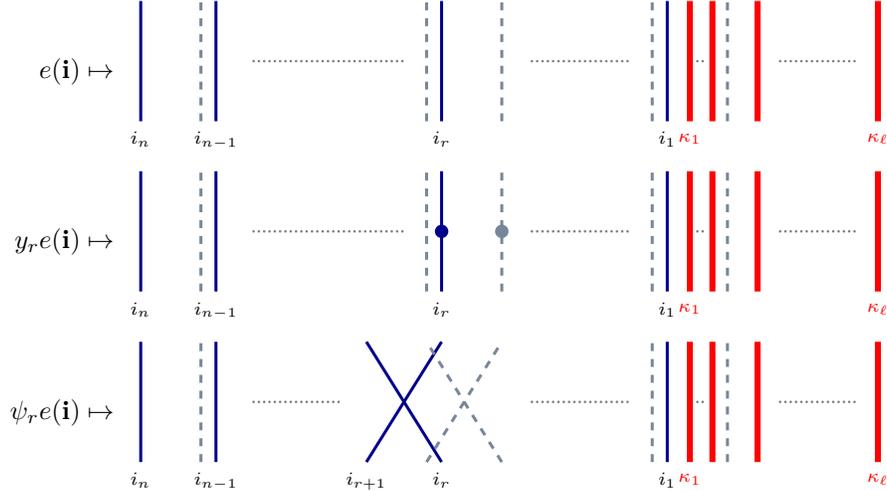
   \begin{tikzpicture}[yscale=0.8, every node/.style={font=\scriptsize}]
      \draw[solid](1,0)node[below]{$i_n$}--++(0,2);
      \draw[ghost](1.8,0)--++(0,2);
      \draw[solid](2,0)node[below]{$i_{n-1}$}--++(0,2);
      \draw[dots](2.5,1)--++(2,0);
      \draw[ghost](4.8,0)--++(0,2);
      \draw[solid](5,0)node[below]{$i_r$}--++(0,2);
      \draw[ghost](5.8,0)--++(0,2);
      \draw[dots](6.2,1)--++(1.3,0);
      \draw[ghost](7.8,0)--++(0,2);
      \draw[solid](8,0)node[below]{$i_1$}--++(0,2);
      \draw[redstring](8.3,0)node[below]{$\kappa_1$}--++(0,2);
      \draw[dots](8.4,1)--++(.1,0);
      \draw[redstring](8.6,0)--++(0,2);
      \draw[ghost](8.8,0)--++(0,2);
      \draw[redstring](9.2,0)--++(0,2);
      \draw[dots](9.5,1)--++(1,0);
      \draw[redstring](10.8,0)node[below]{$\kappa_\ell$}--++(0,2);
    \end{tikzpicture}\\
   y_re(\bi) & \mapsto
   \begin{tikzpicture}[yscale=0.8, every node/.style={font=\scriptsize}]
      \draw[solid](1,0)node[below]{$i_n$}--++(0,2);
      \draw[ghost](1.8,0)--++(0,2);
      \draw[solid](2,0)node[below]{$i_{n-1}$}--++(0,2);
      \draw[dots](2.5,1)--++(2,0);
      \draw[ghost](4.8,0)--++(0,2);
      \draw[solid,dot](5,0)node[below]{$i_r$}--++(0,2);
      \draw[ghost,ghostdot](5.8,0)--++(0,2);
      \draw[dots](6.2,1)--++(1.3,0);
      \draw[ghost](7.8,0)--++(0,2);
      \draw[solid](8,0)node[below]{$i_1$}--++(0,2);
      \draw[redstring](8.3,0)node[below]{$\kappa_1$}--++(0,2);
      \draw[dots](8.4,1)--++(.1,0);
      \draw[redstring](8.6,0)--++(0,2);
      \draw[ghost](8.8,0)--++(0,2);
      \draw[redstring](9.2,0)--++(0,2);
      \draw[dots](9.5,1)--++(1,0);
      \draw[redstring](10.8,0)node[below]{$\kappa_\ell$}--++(0,2);
    \end{tikzpicture}\\
   \psi_re(\bi) & \mapsto
   \begin{tikzpicture}[yscale=0.8, every node/.style={font=\scriptsize}]
      \draw[solid](1,0)node[below]{$i_n$}--++(0,2);
      \draw[ghost](1.8,0)--++(0,2);
      \draw[solid](2,0)node[below]{$i_{n-1}$}--++(0,2);
      \draw[dots](2.5,1)--++(1.2,0);
      \draw[solid](4,0)node[below]{$i_{r+1}$}--++(1,2);
      \draw[solid](5,0)node[below]{$i_{r}$}--++(-1,2);
      \draw[ghost](4.8,0)--++(1,2);
      \draw[ghost](5.8,0)--++(-1,2);
      \draw[dots](6.2,1)--++(1.3,0);
      \draw[ghost](7.8,0)--++(0,2);
      \draw[solid](8,0)node[below]{$i_1$}--++(0,2);
      \draw[redstring](8.3,0)node[below]{$\kappa_1$}--++(0,2);
      \draw[dots](8.4,1)--++(.1,0);
      \draw[redstring](8.6,0)--++(0,2);
      \draw[ghost](8.8,0)--++(0,2);
      \draw[redstring](9.2,0)--++(0,2);
      \draw[dots](9.5,1)--++(1,0);
      \draw[redstring](10.8,0)node[below]{$\kappa_\ell$}--++(0,2);
    \end{tikzpicture}\\
   \end{align*}
   for $\bi\in\cIa$.
\end{Theorem}

The idea behind the proof of \autoref{ISO1} is to use the relations to
pull all of the solid strings in a diagram $D\in\Web(\bom,\bom)$ to the
left of all of the red strings and then check that the relations are
preserved by the map $g\map{\RG[\alpha]}\WA[\alpha](\bom)$ given in the
statement of \autoref{ISO1}.  For us the important point is that instead
of working in $\RG[\alpha]$ we can apply the isomorphism~$g$ and work in
$\WA[\alpha](\bom)$.  We use similar ideas in
\autoref{subsection:main_results} to prove our main results.

\section{Regular diagrams and shifted tableaux combinatorics}
\label{S:ShiftedRegularity}

  This chapter is the technical heart of this paper. It prepares all the tools we will need in the next section to prove that the graded Hecke algebras of type $G(\ell,1,n)$
  have a shift automorphism, which will imply that the Hecke algebras of type
  $G(\ell,p,n)$ are skew cellular by
  \autoref{proposition:sigma_cellular_implies_extended}. All of the
  calculations take place inside the diagrammatic
  Cherednik algebra ~$\WA[\alpha]$. The key point is that the special choice of loading
  made in \autoref{D:multicharge} ensures that the dominance order $\tedom$ for the
  cellular basis of $\RG[\alpha]\cong\Eoma\WA[\alpha]\Eoma$ in
  \autoref{ISO1} is compatible with the shifted tableaux combinatorics
  that we introduce later in this chapter.

\subsection{Regular diagrams}
\label{subsection:regular_diagrams}

This section defines a class of diagrams that are easy to work with and
which play a key role in the arguments of \autoref{subsection:particular_class_diagrams}.

\begin{Definition} Let $D$ be a Webster diagram.
A \emph{singular} crossing in $D$ is a crossing between a solid
$i$-string and either:
\begin{itemize}
\item another solid $i$-string
\item a red $i$-string, or
\item a ghost $iq^{-1}$-string.
\end{itemize}
A crossing is \emph{regular} if it is not singular. A diagram $D$ is a
\emph{regular diagram} if $D$ has no dots and all crossings in $D$ are
regular. A \emph{singular} diagram is any diagram that is not regular.
\end{Definition}

In particular, any crossing that does not involve a solid string is
regular. Note that regular crossings are preserved by the relations in
\autoref{D:RationalCherednik} and by isotopy. An element of $\WA$ is
\emph{regular} if it is the image of a regular Webster diagram. By
assumption, regular diagrams do not contain any of the exceptional
crossings in \autoref{D:RationalCherednik}, so the span of the regular
diagrams in  $\WA$ is a subalgebra of $\WA$.

The next result is the analogue for regular diagrams of the algorithm for reducing words in the symmetric group.

\begin{Proposition}
\label{P:resolve_bigon}
Let $\gbconf,\gbconf'$ be two generalised partitions and let $C \in
\Web(\gbconf,\gbconf')$ be a regular diagram. Then~$C$ is equal to a
regular diagram with no generalised double crossings.
\end{Proposition}

\begin{proof}  Number the strings in $C$ from left to right along the top of the diagram as
  $s_1, s_2, \dots, s_{2n+\ell}$.

   Each string cuts the diagram into two
  pieces, say $\Left(s_k)$  and $\Right(s_k)$. In $\WA$, we claim that the diagram $C$ is equal to a diagram that does not have any generalised double crossings in $\Left(s_k)$, for $1\leq k\leq 2n+\ell$, and if $s_k$ is a not a ghost
  string then the only crossings of non-ghost strings in $\Left(s_k)$ are between
  non-ghost strings $s_a$ and $s_b$ with $\min\set{a,b}<k$. Of course, a
  non-ghost string is either a red string or a solid string.

  We prove the claim by arguing by induction on $k=1,2,\dots,2n+\ell$.
  If $k=1$ then we can use the non-exceptional triple
  crossing relation~(C) from \autoref{D:RationalCherednik} to pull the  string $s_1$ to the left through
  any crossings in $\Left(s_1)$. By induction we assume that the claim is true for the strings
  $s_1,\dots,s_{k-1}$. If $s_k$ is a ghost string there is
  nothing to prove so we may assume that $s_k$ is a non-ghost string. To
  show that the claim holds for $s_k$, use the non-exceptional triple
  crossing relation~(C) from \autoref{D:RationalCherednik} to pull the  string $s_k$ to the left through
  any crossing in $\Left(s_k)$ that involve two larger strings. Pulling
  $s_k$ through a crossing does not destroy any generalised double
  crossings in the diagram. Moreover, for any generalised double crossing $D$ in $\Left(s_k)$ that involve $s_k$, any string which goes through the region
  surrounded by $D$ can be moved away from the region by using the non-exceptional triple
  crossing relation~(C) from \autoref{D:RationalCherednik} as $C$ is regular. As a result, we can apply the non-exceptional relation (B) from \autoref{D:RationalCherednik} to
  pull apart any generalised double crossings in $\Left(s_k)$ that involve $s_k$.
  Note that the two strings do not have the same residue since all
  crossings are regular. Observe that the crossings involving string
  $s_j$ for $j<k$ are unchanged in this process. After a finite
  number of steps we will show that all of the crossings
  between non-ghost crossings in $\Left(s_k)$ will involve a string
  $s_j$, with $j<k$, and $s_k$ does not meet any other string
  twice. This completes the proof of the inductive step and hence proves
  the lemma.
\end{proof}

%
%
%
%
%
%
%

\begin{Corollary}
\label{C:begin_end_points}
Let $\gbconf,\gbconf'$ be two generalised partitions and let $C, D \in
\Web(\gbconf,\gbconf')$ be two regular diagrams such that the string
starting from $\xcoord(\gamma)$ in $C$ and $D$ has the same residue and
the same end points, for all $\gamma \in \gbconf$. Then $C = D$
in~$\WA(\gbconf,\gbconf')$.
\end{Corollary}

\begin{proof}
By~\autoref{P:resolve_bigon},  in $\WA$ the diagram $C^*D$ is equal to a
diagram $E$ that does not contain any generalised double crossings. By
assumption, each string in $E$ starts and ends at the same point, so $E$
does not contain any crossings. Hence, $E$ is an idempotent diagram by
\autoref{L:no_crossing_diagram_is_idempotent}. By the same argument,
$CC^*$ is also equal to an idempotent diagram. Thus, in $\WA$,
\[
D = (CC^*)D = C (C^* D) = C,
\]
which completes the proof.
\end{proof}


Hence, once we fix the start and end positions of the $n$ solid strings, together with their residues, then the set of regular diagrams can be identified with a subgroup of~$\Sym_n$.

%

  \subsection{Shifted tableaux combinatorics}\label{SS:ShiftedComb}
  Recall from \autoref{D:Rlpn} that $\Rlpn$ is defined as an algebra
  of $\sigma$-fixed points:
  \[\Rlpn = \left(\RG\right)^\sigma.\]
  Motivated by \autoref{ISO1}, we want to consider the $\sigma$-fixed
  point subalgebra of the diagrammatic Cherednik algebra but it is not
  clear how to extend $\sigma$ to an automorphism of~$\WA(\bom)$.  This
  section introduces a combinatorial \textit{shift operator} on the set
  of nodes that will allow us to extend $\sigma$ to an automorphism of
  $\WA(\bom)$.

  Define a shift operation $\sigma_\Nscr$ on the set of nodes $\Nodes$ from
  \autoref{nodes} by
  \begin{equation}\label{shiftaction1}
    \sigma_\Nscr(r, c, l) \coloneqq\begin{cases*}
        (r, c, l+1), &if $l\not\equiv 0\pmod{p}$,\\
        (r, c, l+1-p), &if $l\equiv 0\pmod{p}$.
    \end{cases*}
  \end{equation}
  We usually abuse notation and write $\sigma=\sigma_\Nscr$, as the
  meaning will be clear from context.
  Equivalently, if we write $l=ap+b+1$, where $0\le a<d$ and $0\le b<p$,
  then
  \[
    \sigma(r, c, ap+b+1) =\begin{cases*}
        (r, c, ap+b+2), &if $b\ne p-1$,\\
        (r, c, ap+1), &if $b=p-1$.
    \end{cases*}
  \]

%

  Let $\Std(\Parts)=\bigcup_{\blam\in\Parts} \Std(\blam)$.

  \begin{Definition}\label{D:sigmaP}
    Define a map $\sigma_{\mathscr{P}}\map{\Parts}{\Parts}$ by
    $\sigma_{\mathscr{P}}(\blam)=\set{\sigma_\Nscr(\gamma)|\gamma\in\blam}$, for
    $\blam\in\Parts$. Similarly, let
    $\sigma_{\Std}\map{\Std(\Parts)}{\Std(\Parts)}$ be given by
    $\sigma_{\Std}(\t)=\t\circ\sigma_{\mathscr{P}}$, for $\t\in\Std(\Parts)$.
  \end{Definition}

  The definitions readily imply that $\sigma_{\mathscr{P}}(\blam)\in\Parts$
  if $\blam\in\Parts$, so $\sigma_{\mathscr{P}}$ is well-defined. Namely, the partition $\sigma_{\mathscr{P}}(\blam)$ is obtained from $\blam$ by a certain permutation of its  components.
  Similarly, if $\t\in\Std(\blam)$ then
  $\sigma_{\Std}(\t)\in\Std(\sigma_{\mathscr{P}}(\blam))$. As with the
  automorphism $\sigma=\sigma^\bLam_n$ of~$\RG[\alpha]$, we usually omit
  the subscript and write $\sigma=\sigma_{\mathscr{P}}$ and
  $\sigma=\sigma_{\Std}$. In this way, we think of~$\sigma$ as:
  \begin{itemize}
    \item the automorphism $\sigma^\bLam_n$ of $\RG[\alpha]$
    \item the map $\sigma_{\mathscr{P}}$ of order~$p$ on the set~$\Parts$ of
          $\ell$-partitions
    \item the map $\sigma_{\Nscr}$ on the set $\Nodes$ of nodes
    \item the map $\sigma_{\Std}$ on the set $\Std(\Parts)$ of standard tableaux.
  \end{itemize}
  This should not cause any ambiguity because the meaning will always be
  clear from the context. It is not yet clear how the map
  $\sigma^\bLam_n$ is related to the other three combinatorially
  defined maps but we will ultimately see that the triple of maps
  $(\sigma^\bLam_n,\sigma_{\mathscr{P}}, \sigma_{\Std})$ is a
  shift-automorphism in the sense of \autoref{D:sigma_cellular} (see~\autoref{T:sigmaCst}).

%

  \begin{Example}
   Suppose that $p=p'=2$ and $n=4$. Two standard tableaux
   $\t\in\Std(2|1^2)$ and $\s\in\Std(1^2|2)$ are:
    \[
       \t=\Biggl(\space\RussianTableau{{1,3}}\hspace*{2mm}\Bigg|\hspace*{2mm}
                \RussianTableau{{2},{4}}\space\Biggr)
       \qquad\text{and}\qquad
       \s=\Biggl(\space\RussianTableau{{2},{4}}\hspace*{2mm}\Bigg|\hspace*{2mm}
                \RussianTableau{{1,3}}\space\Biggr)
    \]
  Then $\s=\sig\t$ and $\t=\sig\s$.
  \end{Example}

\begin{Lemma}\label{L:SigmaResidue}
Suppose that $\gamma\in\Nodes$. Then
   $\res(\sigma(\gamma))=\eps\res(\gamma)$.

\end{Lemma}

\begin{proof}
  Let $\gamma=(r,c,l)$ and write $l=ap+b+1$, where $0\le a<d$ and
  $0\le b<p$. By the remarks after \autoref{E:residues},
  $\res(\gamma)=\eps^bq^{\rho_{a+1}+c-r}$, which implies the result.
\end{proof}

\begin{Lemma} \label{P:SigmaNondecreasing}
  Suppose that $\xcoord(\gamma)>\xcoord(\gamma') + k$, where
  $\gamma,\gamma'\in\blam$ with $\blam\in\Parts$ and $k \in \set{-1,0,1}$.
  Write $\gamma=(r,c,l)$ and $\gamma'=(r',c',l')$, where $l=ap+b+1$ and
  $l'=a'p+b'+1$, with $0\le a,a'<d$ and $0\le b, b'<p$.
  Then
  $\xcoord\left(\sigma(\gamma)\right)>\xcoord\left(\sigma(\gamma')\right) + k$
  unless the following four conditions hold:
  \[
     a=a',
       \quad 0\le b'<b=p-1,
       \quad c-r=c'-r' + k
       \quad \text{and}
       \quad q^k \res(\gamma')=\eps^{b'+1}\res(\gamma).
  \]
  In particular, if $\res(\gamma)=q^k\res(\gamma')$ then
  $\xcoord\left(\sigma(\gamma)\right)>\xcoord\left(\sigma(\gamma')\right) + k$.
\end{Lemma}

\begin{proof}
  By \autoref{L:xcoord}(b),
  $\xcoord(\sigma(\gamma))>\xcoord(\sigma(\gamma')) + 1$ if $a>a'$, so
  we may assume that $a=a'$.  Recalling the definition of the loading
  $\xcoord$ from~\autoref{E:xcoord},
  if $b<p-1$ or $b = b'$ then $\xcoord(\sigma(\gamma))-\xcoord(\gamma)
     \ge\xcoord(\sigma(\gamma'))-\xcoord(\gamma')$ so we can, and do,
  assume that $b' < b=p-1$.

  Notice that
  $\tfrac{e'}{ep'}=\tfrac1{mp'}=\tfrac1p$ and that
  $\tfrac{2n}N\le\tfrac1{e(\ell+1)}$ since $N\ge 2ne(\ell+1)$ by~\autoref{E:inequality_N}. Using
  these facts for the third and fourth inequalities,
  \begin{align*}
     k & < \xcoord(r,c,l)-\xcoord(r',c',l')
         = c-r-c'+r'+\tfrac{p-1-b'}{e}\bigl(\tfrac{e'}{p'}-\tfrac1{\ell+1}\bigr)
            + \tfrac{r'+c'-c-r}N\\
       & \le c-r-c'+r'+(p-1)\bigl(\tfrac1p-\frac1{e(\ell+1)}\bigr)
            + \tfrac{2n}N\\
       & \le c-r-c'+r'+1-\tfrac1p-\tfrac{p}{e(\ell+1)}
            +\tfrac1{e(\ell+1)}+\tfrac1{e(\ell+1))}\\
       & \le c-r-c'+r'+1-\tfrac1p.
  \end{align*}
  Hence, $c-r\ge c'-r' + k$. 

  A similar calculation, replacing $b=p-1$ with $0$ and $b'$ with
  $b'+1$, shows that
  \[
     \xcoord\left(\sigma(\gamma)\right)-\xcoord\left(\sigma(\gamma')\right)
     =\tfrac{b'+1}{e}\bigl(\tfrac{1}{\ell+1}-\tfrac{e'}{p'}\bigr)
       +c-r-c'+r' - \tfrac{r+c-c'-r'}N
       > c - r - c' + r'-\tfrac{b'+1}{p}.
  \]
  Consequently, if $c-r>c'-r' + k$ then
  $\xcoord\left(\sigma(\gamma)\right)>\xcoord\left(\sigma(\gamma')\right) + k$.

Therefore, $\xcoord\left(\sigma(\gamma)\right)\leq
\xcoord\left(\sigma(\gamma')\right) + k$ only if $c-r=c'-r' + k$. In
this case, by the remarks following \autoref{E:residues},
  \[
  \res(\gamma')=\eps^{b'}q^{v_{a+1}+c'-r'}=\eps^{b'}q^{v_{a+1}+c-r - k}
                =q^{-k}\eps^{b'+1}\res(\gamma),
  \]
  completing the proof of the first part of the lemma.

  Finally, if $\res(\gamma) = q^k \res(\gamma')$ then for any $0 \leq b' < p-1$ we have $q^k \res(\gamma') \neq \eps^{b'+1}\res(\gamma)$ since $\eps$ has order $p$, thus $\xcoord(\sigma(\gamma)) > \xcoord(\sigma(\gamma'))$.
\end{proof}

\begin{Remark}
  It follows from the proof of \autoref{P:SigmaNondecreasing} that if the four
  conditions of the proposition are satisfied then
  $\xcoord\left(\sigma(\gamma)\right)-\xcoord\left(\sigma(\gamma')\right)
      <\tfrac{p}{e(\ell+1)}$.
  Hence, it can still happen that
  $\xcoord\left(\sigma(\gamma)\right)>\xcoord\left(\sigma(\gamma')\right)$.
\end{Remark}

\autoref{P:SigmaNondecreasing} implies that $\sigma$
respects the $\tedom$ partial order. More precisely, we have:

\begin{Corollary} \label{nondescrease1}
  Let $\blam, \bmu\in\Parts$ and suppose that $\blam\tedom\bmu$.  Then
  $\sig\blam \tedom \sig\bmu$.
\end{Corollary}

\begin{proof}
  By \autoref{D:DominanceOrder}, $\blam\tedom\bmu$ if and only if there
  exists a bijection $\theta\map\blam\bmu$ such that
  \[ \res(\theta(\gamma))=\res(\gamma)\qquad\text{and}\qquad
     \xcoord(\theta(\gamma))\le\xcoord(\gamma),
     \qquad\text{for all }\gamma\in\blam.
  \]
  Let $\theta'=\sigma\circ\theta\circ\sigma^{-1}$. Then $\theta'$ is a bijection from
  $\sig\blam$ to $\sig\bmu$ and if $\gamma\in\sig\blam$ then
  $\res(\theta'(\gamma))=\res(\gamma)$, by \autoref{L:SigmaResidue}, and
  $\xcoord(\theta'(\gamma))\le\xcoord(\gamma)$, by
  \autoref{P:SigmaNondecreasing}. Hence, $\sig\blam\tedom\sig\bmu$.
\end{proof}


Abusing notation we extend $\sigma$ to a map on the set of standard
tableaux.

\begin{Definition}\label{D:sigT}
  Let $\t\in\Std(\blam)$, for $\blam\in\Parts$. Let
  $\sig\t\map{\sig\blam}\set{1, \dots, n}$ be the standard $\sig\blam$-tableau
  given by $\sig\t=\t\circ\sigma^{-1}$.
\end{Definition}

Note that the $\sig\blam$-tableau $\sig\t$ is a standard $\sig\blam$-tableau since $\sigma$ just
permutes the components of~$\t$. We now extend the dominance order $\tedom$ on $\Parts$ to the set of standard
tableaux by defining
\[
  \s\tedom\t\qquad\text{if}\qquad \Shape(\s_{\downarrow k}) \tedom \Shape(\t_{\downarrow k})
      \quad\text{ for all}\, 1\le k\le n.
\]
As before, write $\s\tdom\t$ if $\s\tedom\t$ and $\s\ne\t$.

\begin{Lemma}\label{L:SigmaOnTableaux}
  Let $\s\in\Std(\bmu)$ and $\t\in\Std(\blam)$, for $\blam,\bmu\in\Parts$. Then
  $\deg\sig\t=\deg\t$ and if $\s\tdom\t$ then $\sig\s\tdom\sig\t$.
\end{Lemma}

\begin{proof}
  It is immediate from \autoref{P:SigmaNondecreasing} that
  $\sig\s\tedom\sig\t$ if $\s\tedom\t$. Recall the definition of
  $\deg\t$ from \autoref{E:TableauDegree} and
  observe that $\sigma$ induces bijections
  $\Add_k(\t)\bijection\Add_k(\sig\t)$ and $\Rem_k(\t)\bijection\Rem_k(\sig\t)$
  by \autoref{L:SigmaResidue} and \autoref{P:SigmaNondecreasing}, for
  $1\le k\le n$. Hence, $\deg\sig\t=\deg\t$.
\end{proof}

%

The results in this section suggest that it is not unreasonable to
expect that $\sigma(C_{\s\t})=C_{\sig\s,\sig\t}$. Before we can prove
this we first need to extend $\sigma$ to an automorphism of $\WAA(\bom)$
and then construct explicit diagrams $C_\t$ for which this is true. We
do this in the next sections.

\subsection{Applying \texorpdfstring{$\sigma^\bLam_n$}{sigma} in the diagrammatic Cherednik algebra}
\label{SS:Embedding}
This section identifies the image of the automorphism $\sigma^\bLam_n$ of $\Rlpn$ of~\autoref{sigmaIso} under the
isomorphism of \autoref{ISO1}. A key ingredient is the diagram automorphism $D \mapsto D^\cyc$ introduced in~\cite{Bowman:ManyCellular}, and its generalisations below, that play a crucial role in the proof of our main results in~\autoref{S:SkewCellularity}.

Let $\alpha\in\Comp$ be an $\cI$-composition of~$n$. By \autoref{ISO1}, there is a $\Z$-graded
$R$-algebra isomorphism $\RG[\alpha]\cong \Eoma \WA[\alpha]\Eoma$, so we
now identify these two algebras.
Let $[\alpha]$ be the orbit of $\alpha$ under the action of the finite group $\<\sigma\>$, as described in~\autoref{E:action_sigma_compositions}. Define
\begin{equation}\label{E:WAA}
  \RGa=\bigoplus_{\beta\in[\alpha]}\RG[\beta],\quad
  \EomA=\bigoplus_{\beta\in[\alpha]}\Eoma[\beta]\quad\text{and}\quad
  \WAA=\bigoplus_{\beta\in[\alpha]}\WA[\beta].
\end{equation}

The isomorphism $\sigma$ in \autoref{sigmaIso} induces an
automorphism of $\RGa$. Hence, we can regard~$\sigma$ as a homogeneous
$R$-algebra automorphism of $\WAA(\bom)=\EomA\WAA\EomA$. The next result
gives a more precise description of $\sigma$ considered as an
automorphism of $\WAA(\bom)$, or equivalently, of $\WA(\bom)$.

Let $\tsig$ be  the automorphism of the set of Webster diagrams that
sends a diagram~$D$ to the diagram that has the same strings but where
the residues of the solid and ghost strings are multiplied by~$\eps$.
Note that $\tsig$ is not well-defined as an automorphism of $\WA(\bom)$,
or even as a map $\WA(\bom)\to\WA(\bom)$, since it is not compatible
with the exceptional defining relations (for $\tsig$ does not change the residues of the red strings).  The next lemma is immediate
from the definitions.


\begin{Lemma}
\label{L:tsig_star_commute}
If $D$ is any Webster diagram then $\left(\tsig(D)\right)^* = \tsig\left(D^*\right)$.
\end{Lemma}

Let $\gbconf$ be a generalised partition and $a \in \{0,\dots,d-1\}$.
Then $\gbconf$ is \emph{$a$-bounded} if $\xcoord(\gamma) <
\xcoord(0,0,ap+1)$ for all $\gamma \in \gbconf$. That is, all the nodes
in $\gbconf$ are to the left of the $(ap+1)$th red string
$\rstr_{ap+1}$. In other words, the nodes are to the right of at most
$ap$ red strings. For example, the $\ell$-partition $\bom$ is
$0$-bounded.


Let $a \in \{0,\dots,d-1\}$ and let $\gbconf,\gbconf'$ be two $a$-bounded generalised partitions. Let $D \in \Web(\gbconf,\gbconf')$ be a Webster diagram. Then $D$ is an \emph{$a$-bounded diagram} if all of its solid string are to the left of $\rstr_{ap+1}$. Further, $D$ is an \emph{$\bom$-diagram} if $\gbconf=\gbconf'=\bom$ and $D$  is $0$-bounded. In particular,  not all
Webster diagrams of type $(\bom,\bom)$ are $\bom$-diagrams.  By
\autoref{ISO1}, $\WA(\bom)$ is spanned by $\bom$-diagrams.


Now let $a \in \{0,\dots,d-1\}$ and let $\gbconf,\gbconf'$ be two
$a$-bounded generalised partitions. Let $D \in \Web(\gbconf,\gbconf')$.
Let $\str$ be a solid $i$-string in $D$ and let $l \in
\{ap+1,\dots,\ell\}$. Suppose that $q^{\kappa_{l}}=i$ and that $\str$
crosses the $l$th red string $\rstr_{l}$. We pull all solid (and ghost) strings
that are to the right of $\rstr_{ap+1}$ to the left of $\rstr_{ap+1}$
while still staying to the right of $\rstr_{ap}$ if $a > 0$, and for each crossing involving a solid $i$-string (from southwest to northeast)
and a red string (from southeast to northwest) place a dot at the position of the
crossing in $D$. The argument used in the proof of \cite[Proposition~6.19, Figure 18]{Bowman:ManyCellular} now shows that by iterating this process we get
a single diagram $D^{a,\cyc}$,  which is denoted by $\overline{D}$ in \cite[Proposition~6.19]{Bowman:ManyCellular}. When $a = 0$ we simply write $D^\cyc \coloneqq D^{0,\cyc}$. Moreover, by construction,
\begin{equation}\label{E:composition_cyc}
  (D^{{a+1},\cyc})^{a,\cyc} = D^{a,\cyc}, \qquad\text{ whenever }0\le a<d-2.
\end{equation}


The following result is implicit in the proof
of\cite[Proposition~6.19]{Bowman:ManyCellular}.

\begin{Lemma}
  \label{L:acyc}
  Let $a \in \{0,\dots,d-1\}$, let $\gbconf,\gbconf'$ be two $a$-bounded generalised partitions and let $D \in \Web(\gbconf,\gbconf')$. In $\WA(\gbconf,\gbconf')$ we have $D = D^{a,\cyc}$. In particular, if $D$ is of type $(\bom,\bom)$ then $D = D^\cyc$ in $\WA(\bom)$.
\end{Lemma}

\begin{proof} This can be proved using the argument of \cite[Proposition~6.19, Figure 18]{Bowman:ManyCellular}.
\end{proof}

\begin{Definition}
Let $a \in \{0,\dots,d-1\}$.
\begin{itemize}
\item The \emph{$a$-comb} in a diagram is the set of red strings
$\rstr_{ap+1},\dots,\rstr_{(a+1)p}$.
\item A diagram $D$ is \emph{$a$-greedy} if each solid string $\str$
either does not cross $\rstr_{ap+1}$, or $\str$ crosses all of the
strings $\rstr_{ap+1},\dots,\rstr_{(a+1)p}$ in the $a$-comb but does
not cross $\rstr_{(a+1)p+1}$.
\end{itemize}
\end{Definition}
Note that an $a$-bounded diagram is $a$-greedy. In \ref{Ex:regions},
$\set{\rstr_1,\rstr_2}$ is the $1$-comb and $\set{\rstr_3,\rstr_4}$ is the $2$-comb.


\begin{Lemma}
\label{L:tsig_acyc}
 Let $a \in \{0,\dots,d-1\}$, let $\gbconf,\gbconf'$ be two $a$-bounded generalised partitions and let $D \in \Web(\gbconf,\gbconf')$. If $D$ is $a$-greedy  then in $\WA(\gbconf,\gbconf')$ we have
  $\tsig\bigl(D^{a,\cyc}\bigr) = \tsig(D)^{a,\cyc}$.
\end{Lemma}

\begin{proof}
Recall that the red strings $\rstr_{ap+1},\dots,\rstr_{(a+1)p}$ of the $a$-comb have
residues $q^{\kappa_{ap+1}},\eps q^{\kappa_{ap+1}},\dots,\eps^{p-1}
q^{\kappa_{ap+1}}$. Let $\str$ be a solid $i$-string in $D$. If $\str$ does not cross $\rstr_{ap+1}$,
then $D=D^{a,\cyc}$ and $\tsig\bigl(D^{a,\cyc}\bigr) = \tsig(D)^{a,\cyc}$ clearly holds.
Suppose that $\str$ crosses $\rstr_{ap+1}$. By assumption, the diagram $D$ is $a$-greedy, so
$\str$ crosses $\rstr_{ap+1},\dots,\rstr_{(a+1)p}$ and it does not cross
$\rstr_{(a+1)p+1}$. If $i \notin \eps^{\Z} q^{\kappa_{ap+1}}$ then in
$D^{a,\cyc}$ the string $\str$ of $D$ does not gain a dot, so neither
will its image in $\tsig(D)$. On the other hand, if $i = \varepsilon^b
q^{\kappa_{ap+1}}$ with $0\le b<p$, if the string $\str$ gains  $N \in \N$
dots when it is pulled past $\rstr_{ap+b+1}$ during of the operation
$D\mapsto D^{a,\cyc}$, then the image of $\str$ in $\tsig(D)$ also
gains $N$ nodes at the same positions, up to isotopy, during the operation $\tsig(D)\mapsto
\tsig(D)^{a,\cyc}$.
\end{proof}

\autoref{L:tsig_acyc} can be thought as a commutation rule for the operations $D \mapsto D^{a,\cyc}$ and
$D \mapsto \tsig(D)$.

\begin{Proposition}
\label{P:diagram_crosses_group_p_strings}
Let $D \in \Web(\bom,\bom)$ be any diagram such that a solid string crosses a red string
in the $a$-comb only if it cross every red string in the $a$-comb, for  $0\le a<d$.
Then
\(
\tsig(D^\cyc) = \tsig(D)^\cyc
\)
in $\WA(\bom)$.
\end{Proposition}

\begin{proof}
The result is clear if there are no solid-red crossings in $D$ since in
this case $D = D^\cyc$ as Webster diagrams.  Let now $\overline a \in
\{0,\dots,d-1\}$ be maximal such that a solid string crosses a red
string of $\region {\overline{a}} 0$.  We will prove by reverse
induction on $a \in \{0,\dots,\overline{a}\}$ that
\[
\tsig(D^{a,\cyc}) =
\tsig(D)^{a,\cyc}.
\] The initialisation $a = \overline a$ follows
directly from \autoref{L:tsig_acyc} since the assumption of $D$ implies
that $D$ is $\overline a$-greedy. Now let $a \in
\{0,\dots,\overline{a}-1\}$ and assume that $\tsig(D^{a+1,\cyc}) =
\tsig(D)^{a+1,\cyc}$. The diagram $D^{a+1,\cyc}$ is $(a+1)$-bounded by
construction. By assumption of $D$, the diagram $D^{a+1,\cyc}$ is
$a$-greedy thus by \autoref{L:tsig_acyc} we have
\begin{equation}
\label{E:tsig_acyc_induction}
\tsig\bigl((D^{a+1,\cyc})^{a,\cyc}\bigr) = \tsig(D^{a+1,\cyc})^{a,\cyc}.
\end{equation}
By~\eqref{E:composition_cyc} and the induction hypothesis we deduce that
\[
\tsig(D^{a+1,\cyc})^{a,\cyc} = (\tsig(D)^{a+1,\cyc}))^{a,\cyc} = \tsig(D)^{a,\cyc},
\]
thus the hereditary property follows from~\autoref{E:tsig_acyc_induction}. The result now follows from the case $a = 0$.
\end{proof}

\begin{Proposition}\label{P:SigmaDiagram}
Let $D\in\Web(\bom,\bom)$. Then $\sigma(D)=\tsig(D^\cyc)$ in $\WA(\bom)$.
\end{Proposition}

\begin{proof}
By \autoref{sigmaIso} and \autoref{ISO1}, the maps $\tilde{\sigma}$
  and $\sigma$ coincide on the images of the KLR generators in
  $\Web(\bom,\bom)$. Consequently, if  $D$ and $E$ are $\bom$-diagrams
  then
  \[ \tilde{\sigma}(DE)=\sigma(DE)=\sigma(D)\sigma(E)=\tilde{\sigma}(D)\tilde{\sigma}(E),
  \]
  in $\WA$.
  Therefore, in view of \autoref{L:acyc}, writing $D^\cyc \in \WA$ as a product of $\bom$-diagrams, the result follows because
  $\sigma$ is an algebra automorphism of $\WA(\bom)$.
\end{proof}


\subsection{Regions}

Recall that the definition of the set of  nodes $\Nodes$ implies that if
$(r,c,l) \in \Nodes$ then $r - c \in \set{-n+1,\dots,n-1}$. The next
definition allows us to write the set of nodes, and their
$\xcoord$--coordinates, as a disjoint set of \emph{regions}. These regions are the
key to defining an explicit diagrammatic basis of~$\WA(\bom)$ that is
compatible with $\sigma$, which we will use to prove that~$\Rlpn$ is a
graded skew cellular algebra.

\begin{Definition}
\label{D:regions}
For integers $0\le a<d$ and $-n\le\dia\le n$ define
\begin{align*}
\regionnodes &\coloneqq \set{\gamma \in \Nodes | \gamma = (r,c,l), c - r = \dia, l = ap + b + 1 \text{ with } 0 \leq b < p},
\\
\regionblam &\coloneqq \regionnodes \cap \blam, \quad \text{for } \blam \in \Parts,
\\
\region{a}{\dia} &\coloneqq \left[\dia + \frac{1}{e}\left(\rho_{a+1} - \frac{ap}{\ell + 1}\right) - \frac{2n}{N},  \dia + 1 + \frac{1}{e}\left(\rho_{a+1} - \frac{(a+1)p}{\ell + 1}\right)\right] \subseteq \mathbb{R}.
\end{align*}
\end{Definition}

\begin{Example}\label{Ex:regions}
  Let $\ell=4,e=3,p=2,d=2,n=1,N=61,\rho=(0,15)$. An illustration of
  regions in a diagram is given by:
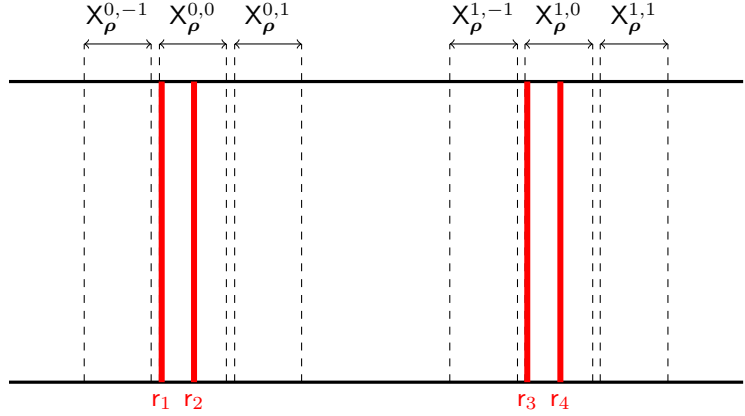
\begin{center}
\begin{tikzpicture}
\foreach \i in {-2,2}
	\draw[very thick] (-2.03,\i) -- (7.73,\i);
\foreach \d in {-1,0,1}
	{
	\draw[dashed] (-.03+\d,-2) --++ (0,4.5);
	\draw[dashed] (.86+\d,-2) --++ (0,4.5);
	\draw[<->] (-.03+\d,2.5) -- (.86+\d,2.5) node[midway, above]{$\region{0}{\d}$};
	}
\foreach \d in {-1,0,1}
	{
	\draw[dashed] (4.83+\d,-2) --++ (0,4.5);
	\draw[dashed] (5.73+\d,-2) --++ (0,4.5);
	\draw[<->] (4.83+\d,2.5) -- (5.73+\d,2.5) node[midway, above]{$\region{1}{\d}$};
	}
\draw[redstring] (0,-2) node[below]{$\rstr_1$} --++ (0,4);
\draw[redstring] (.43,-2) node[below]{$\rstr_2$} --++ (0,4);
\draw[redstring] (4.86,-2) node[below]{$\rstr_3$} --++ (0,4);
\draw[redstring] (5.3,-2) node[below]{$\rstr_4$} --++ (0,4);
\end{tikzpicture}
\end{center}
In this diagram the $0$-comb is $\set{\rstr_1,\rstr_2}$ and the $1$-comb is
$\set{\rstr_3,\rstr_4}$. Hence, the
$a$-comb is contained in the region~$\region{a}{0}$, for
$a\in\set{0,1}$.
\end{Example}

The node $(0, 0, ap + b + 1)$, which corresponds to a red
string, belongs to $\regionnodes[0]$. In particular, note that the red strings in the $a$-comb are precisely the red strings contained in the region $\region{a}{0}$.  For the following lemma, for any
$(a,\dia),(a',\dia') \in \mathbb{Z}^2$, write $(a,\dia) \lex (a',\dia')$
if $(a, \dia)$ is smaller than $(a',\dia')$ in the lexicographic order.
That is, $(a,\dia) \lex (a',\dia')$ if $a < a'$, or $a = a'$ and $\dia<\dia'$.

\begin{Lemma}
\label{lemma:interval}
Suppose that $0\le a<d$, $-n\le\dia\le n$ and
$\gamma = (r,c,l) \in \regionnodes$.
\begin{enumerate}
\item Fix $0\le a'<d$ and $-n\le\dia'\le n$ such that $(a, \dia) \lex (a',\dia')$. Then $\region{a}{\dia} \cap \region{a'}{\dia'} = \emptyset$. More precisely,
if $x \in \region{a}{\dia}$ and $x' \in \region{a'}{\dia'}$ then $x < x'$.
\item We have $\gamma \in \regionnodes$ if and only if $\xcoord(\gamma) \in \region{a}{\dia}$.
\item If $x \in \region{a}{\dia}$ and $\dia < n$ then $x+1 \in \region{a}{\dia+1}$.
\item If  $\xcoord(\gamma) \in \region{a}{\dia}$ then $\xcoord(\sigma(\gamma)) \in \region{a}{\dia}$. That is,  each region is stable under $\sigma$.
\item If $\gamma' = (r',c',l') \in \Nodes$, $\gamma,\gamma' \in \regionnodes$ and $l \neq l'$ then $\res(\gamma) \neq \res(\gamma')$.
\item If $\gamma' = (r',c',l') \in \Nodes$, $\gamma \in \regionnodes$ and $\gamma'\in\regionnodes[\dia-1]$ and $l \neq l'$ then $\res(\gamma) \neq q\res(\gamma')$.
\end{enumerate}
\end{Lemma}

\begin{proof}
Let us prove (a). If $\dia' = \dia + 1$ then $\min \region{a}{\dia+1} -
\max\region{a}{\dia} = -\frac{2n}{N}+\frac{p}{e(\ell + 1)} > 0$, since $N
> 2ne(\ell+1)$ by~\autoref{E:inequality_N}. Hence, the result follows
when $a=a'$. To prove the result when $a < a'$, it suffices to consider
the case when $a' = a+1$ and $\dia = n = -\dia'$. We have
\begin{align*}
\min \region{a+1}{-n} - \max\region{a}{n}
&=
-n + \frac{1}{e}\left(\rho_{a+2}-\frac{(a+1)p}{\ell+1}\right) - \frac{2n}{N}
- \left[n+1 + \frac{1}{e}\left(\rho_{a+1}-\frac{(a+1)p}{\ell+1}\right)\right]
\\
&=
-2n-1 + \frac{1}{e}(\rho_{a+2} - \rho_{a+1}) - \frac{2n}{N}.
\end{align*}
Hence,
\[
\min \region{a+1}{-n} - \max\region{a}{n} > 0
\iff
\rho_{a+2} - \rho_{a+1} > \left(2n+1+\frac{2n}{N}\right)e.
\]
The latter inequality holds since $\rho_{a+2} - \rho_{a+1} \geq (2n+3)e$, by \autoref{D:dCharge},
and $\frac{2n}{N} < \frac{1}{ep'(\ell+1)} < 1$
by~\autoref{E:inequality_N}. This completes the proof of~(a)

We now prove (b).  If $\gamma = (r,c,l) \in \regionnodes$ then $\xcoord(\gamma) = \dia + \frac{1}{e}\left(\kappa_l - \frac{l-1}{\ell+1}\right) - \frac{r+c}{N}$, thus \autoref{L:kappal_l-1_increasing} we obtain
\[
\dia + \frac{1}{e}\left(\rho_{a+1} - \frac{ap}{\ell + 1}\right) - \frac{2n}{N} \leq \xcoord(\gamma) \leq \dia + 1 + \frac{1}{e}\left(\rho_{a+1} - \frac{(a+1)p}{\ell + 1}\right).
\]
That is, $\xcoord(\gamma) \in \region{a}{\dia}$ so (b) now follows
in view of (a).

Parts (c) and (d) are clear from the definitions. For part~(e), writing
$l = ap + b+1$, with $0 \leq b < p$, we have $\res(\gamma)=
q^{\rho_{a+1}+\dia} \eps^b$ by the remark after \autoref{E:residues}.
Thus, writing $l' = ap + b'+1$, with $0 \leq b' < p$, we have $b \neq
b'$ by assumption and so $\res(\gamma') = q^{\rho_{a+1}+\dia} \eps^{b'}
\neq \res(\gamma)$ since $\eps$ has order $p$. The proof of~(f) is
similar.
\end{proof}

%
%

The next lemma will allow us to consider generalised partitions constructed from $\min \region a \dia$.

\begin{Lemma}
\label{L:min_interval_not_red}
Let $0\le a<d$ and $-n\le\dia\le n$.
\begin{enumerate}
\item Let $\gamma = (r,c,ap+1) \in \R \times \R \times \set{1,\dots,\ell}$, where
  $c-r = \dia$ and $c+r = 2n$. Then $\xcoord(\gamma) = \min \region a \dia$.
\item\label{I:min_region_not_red}
If  $l = ap + b + 1$ with  $0\le b<p$, then
  $\min \region{a}{0} < \xcoord(0,0,l)$ and $\min \region{a}{-1} + 1 < \xcoord(0,0,l)$.
  \item We always have
      $\min \region{a}{\dia} \notin \set{\xcoord(0,0,l),\xcoord(0,0,l)-1}$.
\end{enumerate}
\end{Lemma}

\begin{proof}
Part~(a) is clear.  Part~(b) follows from~\autoref{L:kappal_l-1_increasing} and, using~\autoref{lemma:interval}, we deduce part~(c).
\end{proof}

The final result in this section proves a lemma that describes the
$\xcoord$-coordinates of nodes in~$\regionblam[0]$ in terms of the
$\xcoord$--coordinates of the adjacent red strings.

\begin{Lemma}
\label{L:node_region0_between_reds}
Let $0\le a<d$ and $\gamma = (r,c,l) \in \regionblam[0]$.
Then
\[
\xcoord(0,0,l-1) < \xcoord(\gamma) < \xcoord(0,0,l).
\]
\end{Lemma}

\begin{proof} The second inequality follows from \ref{L:min_interval_not_red}.
If $l=ap+1$, then the first equality is clear by \ref{D:dCharge}.
If $l\not\equiv1\pmod p$ then, since $N>2nep'(\ell+1)$ by~\eqref{E:inequality_N},
\begin{align*}
\xcoord(\gamma) - \xcoord(0,0,l-1)
&=
\frac{1}{e}\left(\kappa_l - \frac{l-1}{\ell+1}\right) - \frac{r+c}{N} - \frac{1}{e}\left(\kappa_{l-1} - \frac{l-2}{\ell+1}\right)
\\
&=
\frac{1}{p} - \frac{1}{e(\ell+1)} - \frac{r+c}{N}
\\
&> \frac{1}{p} - \frac{1}{e(\ell+1)} - \frac{1}{ep'(\ell+1)}.
\end{align*}
The last quantity is non-negative since $p \leq \ell$, proving the
result.
\end{proof}

\subsection{Regular and singular diagrams for tableaux}
\label{subsection:particular_class_diagrams}

We are now ready to start constructing the diagrammatic basis of
$\WA(\bom)$ that we will use to prove our main results. In this section we
take the first step of fixing a choice of diagrams $\set{B_\t}$, for
$\t$ a standard tableau, in accordance with \autoref{D:CTdiag}. We construct
the basis elements $B_\t$ by gluing smaller diagrams together.  The next
lemma studies the residues and possible positions of strings at the top
of a $B_\t$ diagram.

\begin{Lemma}
\label{L:positions}
Let $\t \in \Std(\blam)$ with $\blam \in \Parts$ and let $\gamma \in \regionblam$, for $0\le a<d$ and $-n\le\dia\le n$.
\begin{description}[style=multiline,leftmargin=14mm]
\item[(SS)] If $\gamma' \in \regionblam$ and $\res(\gamma) = \res(\gamma')$ and $\t(\gamma) > \t(\gamma')$ then $\xcoord(\gamma) < \xcoord(\gamma')$.
\item[(SG)] If $\gamma' \in \blam^{a, \dia-1}$ and $\res(\gamma) = q\res(\gamma')$ and $\t(\gamma) > \t(\gamma')$ then $\xcoord(\gamma) < \xcoord(\gamma') + 1$.
\item[(GS)] If $\gamma' \in \blam^{a, \dia+1}$ and $q\res(\gamma) = \res(\gamma')$ and $\t(\gamma) > \t(\gamma')$ then $\xcoord(\gamma) + 1 < \xcoord(\gamma')$.
\item[(SR)] If $\dia=0$ and $\gamma' \in \regionnodes[0]$ corresponds to a red string with $\res(\gamma) = \res(\gamma')$ then $\xcoord(\gamma) < \xcoord(\gamma')$.
\item[(GR)] If $\dia=-1$ and $\gamma' \in \regionnodes[0]$ corresponds to a red string with $q\res(\gamma) = \res(\gamma')$ then $\xcoord(\gamma) + 1 < \xcoord(\gamma')$.
\end{description}
\end{Lemma}

\begin{proof}
Write $\gamma = (r,c,l)$ and $\gamma'= (r',c',l')$.
For (SS), which stands for ``solid-solid'', since $\gamma,\gamma' \in \regionblam$ we have $c-r = c'-r'=\dia$, so $q^{\kappa_l} = q^{\kappa_{l'}}$ since $\res(\gamma) = \res(\gamma')$. If $b, b' \in \set{0, \dots, p-1}$ are such that $l - b - 1 = l' - b' - 1 = ap$ we deduce that $\eps^b = \eps^{b'}$ and thus $b = b'$ since $\eps$ has order $p$ and thus $l = l'$. Thus, since $\t$ is standard, the box corresponding to $\gamma$ is higher than the box corresponding to $\gamma'$. Therefore, $r + c > r' + c'$, so $\xcoord(\gamma) - \xcoord(\gamma')= \frac{r'+c'-r-c}{N} < 0$ as desired.

For  (SG), which corresponds to solid-ghost positions, $\gamma, \gamma'$ are necessarily in the same component of $\blam$. Hence, since $\t(\gamma) > \t(\gamma')$ we have $r + c > r' + c'$ and thus $\xcoord(\gamma') + 1 - \xcoord(\gamma) = \frac{r+c-r'-c'}{N} > 0$. The (GS) case, or ghost-solid case, is similar.

Finally, the solid-red (SR) case is deduced from the solid-solid (SS),
and (GR) is deduced from (GS), using the convention that $\t(\gamma')
\coloneqq 0$.
\end{proof}

Recall the definition of a generalised partition,
from~\autoref{E:gen-box-config} and the definition of a regular
diagram from~\autoref{subsection:regular_diagrams}.
In order to make the inductive arguments below clearer, we write
$\bom_n = (1^n|0|\dots|0)$, instead of $\bom$, and let $\region{a}{\dia}(n)$ be
the interval $\region{a}{\dia}$ defined in \autoref{D:regions}, for any
$0\le a<d$ and $-n\le\dia\le n$.

\begin{Definition}
\label{D:gbconf_Ct}
Let $\blam \in \Parts$ and fix a standard tableau $\t \in \Std(\blam)$.
Set $\gamma_k \coloneqq \t^{-1}(k)$, for $1\le k\le n$, and let
$a_k,\dia_k$ be the integers such that $\gamma_k \in
\region{a_k}{\dia_k}(k)$. Define $\gbconf_{\t,k} \coloneqq \min
\region{a_k}{\dia_k}(k)$ and let $\gbconf_\t$ be a generalised partition
such that
\[
  \xcoord(\gbconf_\t) = \set{ \gbconf_{\t,k} | 1\le k \le n}.
\]
\end{Definition}

The inequalities $\min \region{a}{\dia}(n) < \min
\region{a}{\dia}(n-1)$, for any $a$ and $\dia$, together with
\autoref{lemma:interval} and \autoref{L:min_interval_not_red}, ensure
that $\gbconf_\t$ satisfies the requirements of~\eqref{E:gen-box-config}
and hence is a generalised partition.

For the rest of this section we fix a partition $\blam \in \Parts$ and a standard tableau~$\t \in \Std(\blam)$ together
with the associated notation from \autoref{D:gbconf_Ct}.


%
%
%

\begin{Proposition}\label{L:BTreg}
There exists a regular
diagram $\Btreg \in \Web(\blam,\gbconf_\t)$ such that
 $\Btreg$ contains a solid string $\str_k$ of
residue~$\res_\t(k)$ that starts at~$\str_k(0)=\gbconf_{\t,k}$ and ends
at $\str_k(1) = \xcoord(\gamma_k)$, for $1\le k\le n$.
\end{Proposition}

\begin{proof}
  We argue by induction on $n=|\blam|$ to show that such a diagram
  $\Btreg$ exists for any $\blam \in \Parts$ and $\t \in \Std(\blam)$.
  The base case $n = 0$ is immediate because in this case $\t$ is the
  empty tableau and we can define $\Btreg$ to be the diagram with no
  solid (and ghost) strings.

  Now suppose that $n>0$ and let $\s$ be the (standard) tableau obtained
  from $\t$ by removing the (removable) box $\gamma_n=\t^{-1}(n)$. Let
  $\bmu=\Shape(\s)\in\Parts[n-1]$. By induction, there exists a regular
 diagram $\Btreg[\s] \in \Web(\bmu,\gbconf_\s)$ that
  satisfies the conditions of the proposition.

  Define $\Btreg$ to be the diagram obtained by adding the solid string
  $\str_n$ of residue $\res_\t(n)$ to the diagram $\Btreg[\s]$, together
  with its ghost, where the string $\str_n$ is given by
  \[
  \str_n(t) = \begin{cases*}
            \gbconf_{\t,n}, & if $0\le t<1-\epsilon$,\\
            \frac1{\epsilon}\bigl(\gbconf_{\t,n}(1-t) + \xcoord(\gamma_n)(t-1+\epsilon)
                \bigr),
               & if $1-\epsilon\le t\le 1$,
  \end{cases*}
  \]
  where $\epsilon$ is sufficiently small. That is, the solid string
  $\str_n$ is vertical for $0\le t\le 1-\epsilon$ after which it is
  an almost horizontal line connecting the points
  $(\gbconf_{\t,n},1-\epsilon)$ and $(\xcoord(\gamma_n),1)$.

  By induction, $\Btreg$ has the required
  endpoints, so it remains to show that $\Btreg$ is regular. By
  construction, the string $\str_n$ does not cross any other strings
  when $0\le t<1-\epsilon$ because
  $\min \region{a}{\dia}(n) < \min\region{a}{\dia}(n-1)$.
  By definition, $\gamma_n=(r,c,l)\in\region{a}{\dia}$, where for convenience
  we write $a=a_n$ and $\dia=\dia_n$.
  By \autoref{lemma:interval}, when $1-\epsilon\le t\le 1$ the solid string (resp., the ghost string corresponding to) $\str_n$ crosses:
  \begin{itemize}
    \item any solid string $\str_{\gamma'}\in\Btreg[\s]$ for $\gamma' \in
    \bmu^{a,\dia}$ with $\xcoord(\gamma') < \xcoord(\gamma_n)$;
    \item any ghost (resp. solid) string corresponding to
    $\str_{\gamma'}\in\Btreg[\s]$   for $\gamma' \in \bmu^{a, \dia-1}$ (resp.
    $\gamma' \in \bmu^{a,\dia+1}$) with $\xcoord(\gamma')+1 <
    \xcoord(\gamma_n)$ (resp. $\xcoord(\gamma') < \xcoord(\gamma) + 1$);
    \item any red string starting at
    $(\xcoord(0,0,ap+b'+1),1)$ with $0 \leq b' < b$ if $\dia = 0$ (resp.
    $\dia = -1$), where  $b \in \set{0, \dots, p-1}$ is such that $l =
    ap + b + 1$.
  \end{itemize}
  Since $\gamma=\gamma_n$ is the box with the highest label in $\t$, if
  $\gamma' \in \bmu$ then have $\t(\gamma_n) = n > \t(\gamma')$.
  Therefore, by \autoref{L:positions}, all of the crossings above are
  regular. Hence, $\Btreg$ is regular and the proof is complete.
\end{proof}

\begin{figure}[htbp]
\begin{center}
\begin{tikzpicture}
\foreach \i in {-2,0,2}
	\draw[very thick] (-4.02,\i) -- (9.733,\i);

\foreach \d in {-1,0,1}
	{
	\draw[dashed] (-.0331+\d,-2) --++ (0,4.5);
	\draw[dashed] (.8666+\d,-2) --++ (0,4.5);
	\draw[<->] (-.0331+\d,2.5) -- (.8666+\d,2.5) node[midway, above]{$\region{0}{\d}$};
	}

\foreach \d in {-1,0,1}
	{
	\draw[dashed] (6.8334+\d,-2) --++ (0,4.5);
	\draw[dashed] (7.7333+\d,-2) --++ (0,4.5);
	\draw[<->] (6.8334+\d,2.5) -- (7.7333+\d,2.5) node[midway, above]{$\region{1}{\d}$};
	}

\draw[redstring] (0,-2) node[below]{$\rstr_1$} --++ (0,4);
\draw[redstring] (.4333,-2) node[below]{$\rstr_2$} --++ (0,4);
\draw[redstring] (6.8666,-2) node[below]{$\rstr_3$} --++ (0,4);
\draw[redstring] (7.3,-2) node[below]{$\rstr_4$} --++ (0,4);

\usetikzlibrary{decorations.pathreplacing,calligraphy}
\coordinate (nu1) at (-.0082,0);
\coordinate (nu2) at (-.0165,0);
\coordinate (nu4) at (6.8417,0);
\coordinate (nu3) at (.9668,0);
\draw [thick,decorate, decoration = {calligraphic brace, raise=-1.5ex,amplitude=1.5ex}] (-4.5,0) --++ (0,2) node[left,midway]{$\Btreg$};
\draw[solid] (nu1) -- (-.0082,2);
\draw[solid] (nu2) --++ (0,1) -- (.4250,2);
\draw[solid] (nu4) -- (7.2917,2);
\draw[solid] (nu3) -- (1.4208,2);

\draw [thick,decorate, decoration = {calligraphic brace, raise=-1.5ex,amplitude=1.5ex}] (-4.5,-2) --++ (0,2) node[left,midway]{$\Btsing$};
\draw[solid] (-.0082,-2) -- (nu1);
\draw[solid] (-1.0124,-2) -- (nu2);
\draw[solid] (-2.0165,-2) -- (nu3);
\draw[solid] (-3.0207,-2) --++ (0,1.75) -- (nu4);

\draw[densely dotted,orange,very thick] (-.2,-3) rectangle (.5,3.2);
\draw (.5,-3) node[above right,orange]{see \ref*{F:Btzoom}};
\end{tikzpicture}
\end{center}
\caption{An example of a diagram $B_\t = \Btreg\Btsing$ (the distance between $\rstr_2$ and $\rstr_3$ was shorten for expository means)}
\label{F:Bt}
\end{figure}
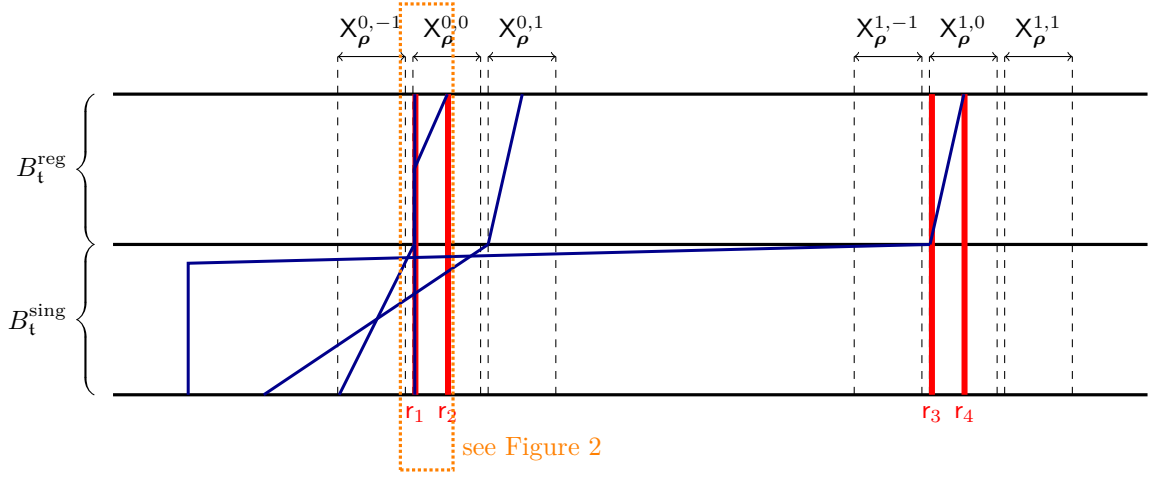

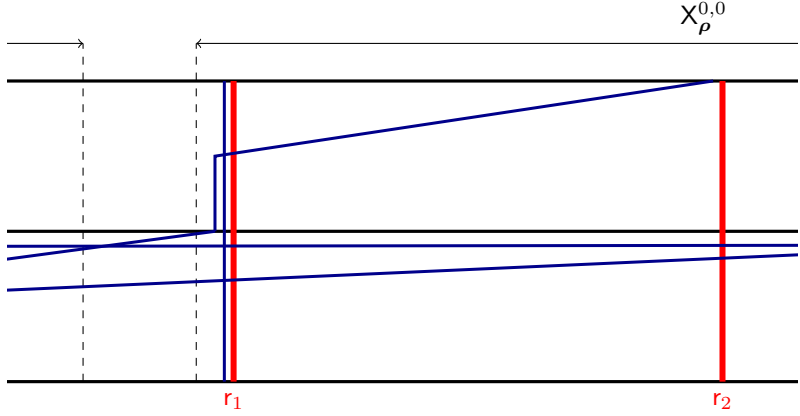
\begin{figure}[htbp]
\begin{center}
\begin{tikzpicture}[xscale=15]
\clip (-.2,-3) rectangle (.5,3.2);
\foreach \i in {-2,0,2}
	\draw[very thick] (-4.02,\i) -- (13.733,\i);

\foreach \d in {-1,0,1}
	{
	\draw[dashed] (-.0331+\d,-2) --++ (0,4.5);
	\draw[dashed] (.8666+\d,-2) --++ (0,4.5);
	\draw[<->] (-.0331+\d,2.5) -- (.8666+\d,2.5) node[midway, above]{$\region{0}{\d}$};
	}

\foreach \d in {-1,0,1}
	{
	\draw[dashed] (10.8334+\d,-2) --++ (0,4.5);
	\draw[dashed] (11.7333+\d,-2) --++ (0,4.5);
	\draw[<->] (10.8334+\d,2.5) -- (11.7333+\d,2.5) node[midway, above]{$\region{1}{\d}$};
	}

\draw[redstring] (0,-2) node[below]{$\rstr_1$} --++ (0,4);
\draw[redstring] (.4333,-2) node[below]{$\rstr_2$} --++ (0,4);
\draw[redstring] (10.8666,-2) node[below]{$\rstr_3$} --++ (0,4);
\draw[redstring] (11.3,-2) node[below]{$\rstr_4$} --++ (0,4);

\coordinate (nu1) at (-.0082,0);
\coordinate (nu2) at (-.0165,0);
\coordinate (nu4) at (10.8417,0);
\coordinate (nu3) at (.9668,0);
\draw [thick,decorate, decoration = {calligraphic brace, raise=-1.5ex,amplitude=2ex}] (-4,0) --++ (0,2) node[left,midway]{$\Btreg$};
\draw[solid] (nu1) -- (-.0082,2);
\draw[solid] (nu2) --++ (0,1) -- (.4250,2);
\draw[solid] (nu4) -- (11.2917,2);
\draw[solid] (nu3) -- (1.4208,2);

\draw [thick,decorate, decoration = {calligraphic brace, raise=-1.5ex,amplitude=2ex}] (-4,-2) --++ (0,2) node[left,midway]{$\Btsing$};
\draw[solid] (-.0082,-2) -- (nu1);
\draw[solid] (-1.0124,-2) -- (nu2);
\draw[solid] (-2.0165,-2) -- (nu3);
\draw[solid] (-3.0207,-2) --++ (0,1.75) -- (nu4);
\end{tikzpicture}
\end{center}
\caption{Focusing on the orange rectangle of \ref{F:Bt}}
\label{F:Btzoom}
\end{figure}

An example of a diagram $\Btreg$ can be found in \ref{F:Bt} and \ref{F:Btzoom}.
Now \autoref{D:gbconf_Ct} and~\autoref{L:BTreg} immediately imply the following result.

\begin{Lemma}
\label{L:string_Btreg_stays_region}
Each string of $\Btreg$ is contained in a single region.
\end{Lemma}


We now define the second part of the diagram $B_\t$.

\begin{Proposition}
\label{L:BTsing}
There exists a  diagram  $\Btsing \in \Web(\gbconf_\t,\bom)$ with solid strings $\str_1,\dots\str_n$, and corresponding ghost strings $\gstr_1,\dots,\gstr_n$, such that:
\begin{enumerate}
  \item\label{I:CTsing_res}
  The string  $\str_k$ has residue $\res_\t(k)$, starts at
  $\xcoord(k,1,1)$  and ends at $\gbconf_{\t,k}$, for $1\le k\le n$.
  \item \label{I:CTsing_a_ignoring}
  In particular, the string $\str_k$ crosses a red string in $\region{a}{0}$ for $a \in \{0,\dots,d-1\}$ if and only if it crosses all the red strings in~$\region{a}{0}$.
  \item \label{I:CTsing_no_crossing_same_region}
  If $\str_k$ crosses $\str_j$ then $(a_k,\dia_k)\ne(a_j,\dia_j)$.
  If $\str_k$ crosses $\gstr_j$ then $(a_k,\dia_k)\ne(a_j,\dia_j+1)$.
  If $\gstr_k$ crosses $\str_j$ then $(a_k,\dia_k+1)\ne(a_j,\dia_j)$.
  \item \label{I:CTsing_red_not_from_same_region}
  If $\str_k$ crosses the red string in $\region{a}{0}$ then
  $\gamma_k\notin\regionblam[0]$. Similarly, if $\gstr_k$ crosses the
  red string in $\region{a}{0}$ then $\gamma_k\notin\regionblam[-1]$.
\end{enumerate}
\end{Proposition}


\begin{proof} Again we argue by induction on the number of solid strings
  $n \geq 0$. When $n=0$ the result is vacuously true so suppose $n>0$.  Let $\s$ be the (standard) $\bmu$-tableau obtained from $\t$ by removing the (removable) box $\gamma_n=\t^{-1}(n)$, where $\bmu \in \Parts[n-1]$ is the $\ell$-partition obtained from $\blam \in \Parts$ by removing~$\gamma_n \in \regionblam[\dia_n][a_n]$. By induction, there exists a diagram $\Btsing[\s] \in \Web(\gbconf_{\s},\bom_{n-1})$ that satisfies the conditions of \autoref{L:BTsing}. As in the proof of \autoref{L:BTreg}, define $\Btsing\in \Web(\gbconf_\t,\bom_{n})$ to be the diagram obtained from $\Btsing[\s]$ by adding the solid string $\str_n$ of residue $\res_\t(n)$, and its ghost, going vertically from $(\xcoord(n,1,1),0)$ until almost the top of the diagram after which $\str_n$ goes almost horizontally to the point $(\gbconf_{\t,n}, 1)$. More explicitly,
\[
\str_n(t) = \begin{cases*}
          \xcoord(n,1,1), & if $0\le t<1-\epsilon'$,\\
          \frac1{\epsilon'}\bigl(\xcoord(n,1,1)(1-t) + \gbconf_{\t,n}(t-1+\epsilon')
              \bigr),& if $1-\epsilon'\le t\le 1$,
\end{cases*}
\]
for some sufficiently small $\epsilon'$.

Note that $\xcoord(n,1,1) \in \region{0}{1-n}(n)$ and $\xcoord(n,1,1)<\xcoord(c,1,1)$ for any $1\leq c<n$.  By definition,
$\dia_n \geq 1-n$. Hence, by \autoref{lemma:interval}(a),
$\xcoord(n,1,1) \leq \gbconf_{\t,n}= \min \region{a_n}{\dia_n}(n)$.  In particular, by induction, the diagram
$\Btsing$
satisfies~\autoref{I:CTsing_res}. Recalling that $\xcoord(n,1,1) < \xcoord(0,0,l)$ for any $l \in \{1,\dots,\ell\}$, we deduce that if $\delta_n \neq 0$ then $\str_n$ satisfies~\autoref{I:CTsing_a_ignoring}, and if $\delta_n = 0$ then $\str_n$ also satisfies~\autoref{I:CTsing_a_ignoring} by~\autoref{L:min_interval_not_red}\autoref{I:min_region_not_red}.  By the same argument, $\str_n$ crosses
no solid or ghost strings ending in $\region{a_n}{\dia_n}$ and
similarly for the ghost string $\gstr_n$ ending in $\region{a_n}{\dia_n + 1}$,
so~\autoref{I:CTsing_no_crossing_same_region} holds. Finally,
condition~\autoref{I:CTsing_red_not_from_same_region} follows from~\autoref{L:min_interval_not_red} and the explicit construction of $\str_n$.
\end{proof}

An example of a diagram $\Btreg$ can be found in \ref{F:Bt}.
Composing the diagrams $\Btreg$ and $\Btsing$ from \autoref{L:BTreg} and
\autoref{L:BTsing} we can now define a diagram $B_\t$ satisfying the
conditions of \autoref{D:CTdiag}.

\begin{Definition}\label{D:RegSing}
Let $B_\t \coloneqq
\Btreg \Btsing \in \Web(\blam,\bom)$.
\end{Definition}

\begin{Proposition}
\label{P:CT}
The diagram $B_\t \in
\Web(\blam,\bom)$
satisfies the assumptions of~\autoref{D:CTdiag}. 
\end{Proposition}

\begin{proof}
By construction, for $1\le k\le n$ the diagram $B_\t$
  has a solid string $\str_k$ of residue $\res_\t(k)$ from
  $(\xcoord(\gamma),1)$ to $(\gbconf_{\t,\t(\gamma)},\frac12)$ and from
  $(\gbconf_{\t,\t(\gamma)},\frac12)$ to
  $\bigl(\xcoord(\t(\gamma),1,1),0\bigr) = (\xcoord[\t](\gamma),0)$,
  so~$B_\t$ satisfies~\autoref{D:CTdiag}\autoref{I:CT_residue}. The
  diagram $B_\t$ has no dots so it
  satisfies~\autoref{D:CTdiag}\autoref{I:CT_no_dots}. To prove that
  $B_\t$
  satisfies~\autoref{D:CTdiag}\autoref{I:CT_minimal_number_crossings},
  it suffices to show that $B_\t$ does not have a generalised double
  crossing.  By \ref{L:string_Btreg_stays_region}, the crossings in $\Btreg$ are
between strings that belong to the same region, while by \autoref{L:BTsing},  the crossings in $\Btsing$ are between strings that end in different regions. Hence, $B_\t$ does not contain any generalised double crossings.
\end{proof}

An illustration of \ref{D:RegSing} is given in~\ref{F:Bt} and~\ref{F:Btzoom} for $\ell=4,e=3,p=2,d=2,N=241,\rho=(0,33),\t=\Bigl(\RussianTableauscale{{1}} \Big| \RussianTableauscale{{2,3}} \Big| \emptyset \Big| \RussianTableauscale{{4}}\Bigr)$.


\subsection{Orbit diagrams}\label{S:OrbitDiagrams}
Recall from \autoref{SS:ShiftedComb} that $\sigma$ is an automorphism of
order~$p$ that acts on $\Parts$ and $\Std(\blam)$, for $\blam\in\Parts$.
This section uses the diagrams from \autoref{D:RegSing} to construct
diagrams indexed by orbits of tableaux under $\sigma$. Ultimately, this will make easier to compute the image under $\sigma$ of the diagrams of the cellular basis.

\begin{Lemma}
\label{P:map_fI}
Let $I \subseteq \R$ be an interval that is neither empty nor a
singleton. We can find a map $f_I : \Nodes \to I$ such that if $\gamma
\in \regionnodes$ and $\gamma' \in \regionnodes[\dia']$, for $0\le a<d$
and $-n\le \dia,\dia'\le n$, then
\[
f_I(\gamma) < f_I(\gamma')  \iff  \xcoord(\gamma) - \dia+\dia' < \xcoord(\gamma').
\]
\end{Lemma}

\begin{proof} This is obvious because we can define
  $f'_I(\gamma)\coloneqq\xcoord(\gamma)-\dia$, where $\gamma = (r,c,l)$ with
  $c-r = \dia$ and then suitably renormalise $f'_I$ so that it has the
  required properties.
\end{proof}

\begin{Corollary}
\label{C:maps_fI_gI}
Let $I \subseteq [ 0, 1]$ be an interval that is neither empty nor a singleton. We can find two maps $f_I,g_I : \Nodes \to I$ so that for any $(a, \dia)$ and any $\gamma, \gamma' \in \Nodes$ with $\gamma' \in \regionnodes$ we have
\begin{align*}
f_I(\gamma) < f_I(\gamma') &\iff g_I(\gamma) > g_I(\gamma') \iff  \xcoord(\gamma) < \xcoord(\gamma'),
&
&\text{if } \gamma \in \regionnodes,
\\
f_I(\gamma) < f_I(\gamma') &\iff g_I(\gamma) > g_I(\gamma') \iff \xcoord(\gamma) + 1 < \xcoord(\gamma'),
&
&\text{if } \gamma \in \regionnodes[\dia-1],
\end{align*}
\end{Corollary}

\begin{proof}
The existence of $f_I$ follows directly from \autoref{P:map_fI}. The existence of $g_I$ follows from the existence of the map $1 - f_{1 - I}$.
\end{proof}

For each $\blam\in\Parts$ let $\so_\blam \in \set{1,\dots,p}$ be the
order of $\blam$ under the action of $\langle \sigma\rangle$ on $\Parts$
and let $\Partss$ be a fixed set of representatives in~$\Parts$ under
the action of~$\sigma$.  The integer $\so_\blam$ divides~$p$, so
$\sp_\blam \coloneqq \frac{p}{\so_\blam} \in \set{1, \dots, p}$ is an
integer.  The cyclic group $\Z/p_\blam\Z\cong\<\sigma^{\so_{\blam}}\>$
generated by $\sigma^{\so_\blam}$ acts on the set $\Std(\blam)$ of
standard $\blam$-tableaux. Let $\Stds$ be a fixed set of
$\Z/\sp_\blam\Z$-orbit representatives with respect to this action.

Let $\blam \in \Parts$ and fix $k \in \set{1, \dots, p-1}$.
Define a decomposition of $\blam = \blaml \sqcup \blamr$ by
\begin{align*}
\blaml &\coloneqq \set{\gamma \in \blam | \xcoord(\sigma^k\gamma) < \xcoord(\gamma)}
\\
&= \set{\gamma = (r,c,l) \in \blam | l = ap + b + 1, 0 \leq a < p \text{ and } p-k \leq b < p},
\\
\blamr &\coloneqq \set{\gamma \in \blam | \xcoord(\sigma^k\gamma) > \xcoord(\gamma)}
\\
&= \set{\gamma = (r,c,l) \in \blam | l = ap + b + 1, 0 \leq a < p \text{ and } 0 \leq b < p - k}.
\end{align*}
In other words, the nodes in $\blaml$ and $\blamr$ move to the left and right, respectively, after applying $\sigma^k$. In the Example of \ref{F:Clamk}, we have $\gamma_1,\gamma_3 \in \blaml[1]$ and $\gamma_2 \in \blamr[1]$.

Now fix a real number $H \in (0, 1)$.  Several applications of
\autoref{C:maps_fI_gI} show that there exist families of real numbers
$(x^{[k]}_\gamma)_{\gamma \in \blam}, (y^{[k]}_\gamma)_{\gamma \in \blamr}$ and
$(z^{[k]}_\gamma)_{\gamma \in \blamr}$ such that:
\begin{itemize}
\item For any $\gamma, \gamma' \in \blaml$ with $\gamma' \in \regionnodes$ we have
\begin{gather*}
\frac{H}{2} < x^{[k]}_{\gamma} < H,
\\
\text{if } \gamma \in \regionnodes \text{ then } \xcoord(\gamma) < \xcoord(\gamma') \iff x^{[k]}_\gamma < x^{[k]}_{\gamma'},
\\
\text{if } \gamma \in \regionnodes[\delta-1] \text{ then } \xcoord(\gamma) + 1< \xcoord(\gamma') \iff x^{[k]}_\gamma  < x^{[k]}_{\gamma'}.
\end{gather*}

\item For any $\gamma, \gamma' \in \blamr$ with $\gamma' \in \regionnodes$ we have
\begin{gather*}
0 < x^{[k]}_\gamma < \frac{H}{2},
\\
\xcoord(0,0,l+k-1) - \frac{1}{N} \leq y^{[k]}_\gamma < \xcoord(0,0,l+k-1), \text{ where } \gamma = (r,c,l),
\\
H < z^{[k]}_\gamma < 1,
\\
\text{if } \gamma \in \regionnodes \text{ then } \xcoord(\gamma) < \xcoord(\gamma') \implies x^{[k]}_\gamma > x^{[k]}_{\gamma'} \text{ and } y^{[k]}_{\gamma} < y^{[k]}_{\gamma'} \text{ and } z^{[k]}_\gamma > z^{[k]}_{\gamma'},
\\
\text{if } \gamma \in \regionnodes[\delta-1] \text{ then } \xcoord(\gamma) + 1 < \xcoord(\gamma') \implies x^{[k]}_\gamma  > x^{[k]}_{\gamma'}  \text{ and } y^{[k]}_{\gamma} < y^{[k]}_{\gamma'} \text{ and } z^{[k]}_\gamma > z^{[k]}_{\gamma'}.
\end{gather*}
\end{itemize}

Given two points $(x_1,y_1)$ and $(x_1,y_2)$, with $x_2\ge x_1$ and
$y_1\le y_2$, let
$(x_1,y_1)\rightsquigarrow(x_2,y_2)$ be the straight line string that:
\[
  \begin{cases*}
    \text{goes from $(x_1,y_1)$ to $(x_2,y_2)$},          & if $y_1\ne y_2$,\\
    \text{goes from $(x_1,y_1)$ to $(x_2,y_2+\epsilon)$}, & if $y_1=y_2$,\\
  \end{cases*}
\]
where $\epsilon$ is sufficiently small.

\begin{Definition}
\label{D:Cslam}
Let $\blam \in \Partss$ and $1\le k<p$. Let $\Cslam \in
\Web(\sigma^k\blam,\blam)$ be the Webster diagram with the following
solid strings:
\begin{itemize}
\item For any $\gamma \in \blaml$, there is a solid string of residue $\res(\gamma)$ given by
\[ (\xcoord(\gamma), 0)\rightsquigarrow(\xcoord(\gamma), x^{[k]}_\gamma)\rightsquigarrow(\xcoord(\sigma^k\gamma), x^{[k]}_\gamma)\rightsquigarrow(\xcoord(\sigma^k\gamma), 1).
\]
\item For any $\gamma \in \blamr$, there is a solid string of residue $\res(\gamma)$ given
by
\[(\xcoord(\gamma), 0)\rightsquigarrow(\xcoord(\gamma), x^{[k]}_\gamma)\rightsquigarrow(y^{[k]}_\gamma, x^{[k]}_\gamma)\rightsquigarrow(y^{[k]}_\gamma, z^{[k]}_\gamma)\rightsquigarrow(\xcoord(\sigma^k\gamma), z^{[k]}_\gamma)\rightsquigarrow(\xcoord(\sigma^k\gamma), 1).
\]
\end{itemize}
\end{Definition}

An illustration of a diagram $\Cslam$ is given  in \ref{F:Clamk}. Note the following result, which follows from~\autoref{lemma:interval}.

\begin{figure}[htbp]
\begin{center}
\begin{tikzpicture}[rounded corners]
\draw[very thick] (-8,2) -- (6,2);
\draw[very thick] (-8,-2) -- (6,-2);
\draw[thick] (-8,0.8) --++ (14,0) node[right]{$H$};
\draw[thick] (-8,-0.6) --++ (14,0) node[right]{$\frac{H}{2}$};
\draw[redstring] (-2,2) --++ (0,-4);
\draw[redstring] (0,2) --++ (0,-4);
\draw[redstring] (2,2) --++ (0,-4);
\coordinate (g1) at (1.7,0);
\draw[solid] (-2.3,2) -- (-2.3,0)  -- (1.7,0)  -- (1.7,-2);

\node (fin1) at (1.7,-2){};
\draw[->] (1.7,-2.5) node[below]{$\xcoord(\gamma_1)$} -- (fin1);
\node (deb1) at (-2.3,2){};
\draw[-<] (-2,2.25) node[right]{$\xcoord(\sigma\gamma_1)$} -| (deb1);

\draw[thin, <-,gray]([xshift=0.1]g1)--++(1.1,0) node[right]{$(\xcoord(\gamma_1), x^{[1]}_{\gamma_1})$};
\node (g3) at (-2.4,-0.2) {};
\draw[ghost] (-2.5,2) -- (-2.5,-0.2) -- (1.5,-0.2) -- (1.5,-2);
\draw (-7.5,1.95) node[above right]{$\xcoord(\sigma\gamma_3)$};
\draw[solid] (-7.5,2) -- (-7.5,-0.2) -- (-3.5,-0.2) --
 (-3.5,-2) node[below]{$\xcoord(\gamma_3)$};
\draw[thin, <-,gray] (g3) -- (-3.5,0.4) node[left] {$(\xcoord(\sigma\gamma_3) + 1, x^{[1]}_{\gamma_3})$};
\draw[solid] (-2.6,2) -- (-2.6,-0.4) -- (1.4,-0.4) -- (1.4,-2);
\draw[solid] (-0.3,2) -- (-0.3,1.5) -- (-2.1,1.5) -- (-2.1,-1.5)-- (-2.4,-1.5) -- (-2.4,-2);

\node (fin2) at (-2.4,-2){};
\draw[->] (-2.4,-2.5) node[below]{$\xcoord(\gamma_2)$} -- (fin2);
\node (deb2) at (-.3,2){};
\draw[-<] (0,2.25) node[right]{$\xcoord(\sigma\gamma_2)$} -| (deb2);

\draw[thin,
<-,gray](-2.1,1.5)--++(-3.0,0)node[left]{$(y^{[2]}_{\gamma_2},z^{[2]}_{\gamma_2})$};
\draw[thin, <-,gray](-0.2,1.5)--++(3,0)node[right]{$(\xcoord(\sigma\gamma_2),z^{[1]}_{\gamma_2})$};
\draw[thin, <-,gray](-2.1,-1.5)--++(4.9,0)node[right]{$(y^{[1]}_{\gamma_2},x^{[1]}_{\gamma_2})$};
\draw[solid] (-0.5,2) -- (-0.5,1.7) -- (-2.2,1.7) -- (-2.2,-1.4) -- (-2.6,-1.4) -- (-2.6,-2);;
\draw[<->] (-7.8,2.5) -- (-3,2.5) node[midway, above]{$\region{a}{-1}$};
\draw[dotted] (-7.8,2.5) -- (-7.8,-2);
\draw[dotted] (-3,2.5) -- (-3,-2);
\draw[<->] (-2.7,2.5) -- (2.1,2.5) node[midway, above]{$\region{a}{0}$};
\draw[dotted] (-2.7,2.5) -- (-2.7,-2);
\draw[dotted] (2.1,2.5) -- (2.1,-2);
\end{tikzpicture}
\end{center}
\caption{Illustration of a diagram $\Cslam[1]$}
\label{F:Clamk}
\end{figure}

\begin{Lemma}
\label{L:string_Cslam_stays_region}
Each string of $\Cslam$ is contained in a single region.
\end{Lemma}

The next result describes all of the crossings the diagrams $\Cslam$,
for $\blam\in\Parts$. Examples of crossings are depicted in \ref{F:crossings}.

\begin{Proposition}
\label{P:crossings_Cslam}
Let $\blam \in \Partss$ and $1\le k<p$ and fix $0\le a<p$ and
$-n\le\dia\le n$. If~$\gamma\in\blam$ let  $\str_\gamma$ be the solid
string in $\Cslam$ that starts at $\xcoord(\gamma)$ and let
$\gstr_\gamma$ be its ghost. Let $\rstr_1,\dots,\rstr_\ell$ be the red
strings in $\Cslam$.
\begin{enumerate}
  \item The diagram $\Cslam \in \Web(\sigma^k\blam,\blam)$ does not
  contain any generalised double crossings.
  \item If $\gamma\in\blamr\cap\regionblam$ then $\str_\gamma$ crosses
  $\str_{\gamma'}$ and $\gstr_{\gamma''}$, for all
  $\gamma'\in\blaml\cap\regionblam$ and $\gamma''\in \blaml\cap\regionblam[\dia-1]$,
  and~$\gstr_\gamma$ crosses $\str_{\gamma'''}$, for all
  $\gamma'''\in\blaml\cap\regionblam[\dia+1]$. Moreover, all of these
  crossings are regular and every crossing between the solid and ghost
  strings in $\Cslam$ is of this form.
  \item  If $\gamma=(r,c,l) \in \blamr \cap \blam^{a, 0}$ then the solid
strings $\str_\gamma$  crosses the red
  string $\rstr_j$ if and only if $l\le j<l+k$. Moreover, the crossing
  of $\str_\gamma$ with $\rstr_l$ is singular and all of the other
  solid-red crossings involving $\str_\gamma$ are regular. The remaining
  solid strings from $\blamr$ do not cross any red strings.
   \item  If $\gamma=(r,c,l) \in \blamr \cap \blam^{a, -1}$ then  the ghost string $\gstr_\gamma$ crosses the red
  string $\rstr_j$ if and only if $l\le j<l+k$.
  \item If $\gamma=(r,c,l) \in \blaml\cap \regionblam[0]$ then
  $\str_\gamma$ crosses the red string $\rstr_j$ if and only if
  $l-p+k\le j<l$. Moreover, all these crossings are regular. The
  remaining solid strings from $\blaml$ do not cross any red strings.
  \item If $\gamma=(r,c,l) \in \blaml\cap\regionblam[-1]$ then the ghost
  string $\gstr_\gamma$ crosses the red string $\rstr_j$ if and only if
  $l-p+k\le j< l$.  The remaining ghost strings from $\blaml$ do not
  cross any red strings.
\end{enumerate}
\end{Proposition}

\begin{proof}
All the results follow directly from~\autoref{D:Cslam} and~\autoref{lemma:interval} because:
\begin{itemize}
\item if $\gamma \in \blamr$ and $\gamma' \in \blaml$ then $\gamma$ and $\gamma'$ belong to different components of $\blam$, so if $\gamma$ and $\gamma'$ belong to the same $\regionblam$ then $\res(\gamma) \neq \res(\gamma')$ by~\autoref{lemma:interval};
\item if $\gamma = (r,c,l) \in \blaml \cap \regionblam[0]$ then $l = ap + b + 1$ with $p - k \leq b < p$ and $\res(\gamma) = q^{\kappa_l} \neq q^{\kappa_{l-i}} = \res(\rstr_{l-i})$ for all $i \in \set{1,\dots,p-k}$ since $ap + 1 \leq l-i = ap + (b-i) + 1 \leq ap + (b-1) + 1$.
\end{itemize}
\end{proof}

\begin{figure}[htbp]
\begin{center}
\begin{tikzpicture}[rounded corners,xscale=2.5]
\draw[very thick] (-4,2) -- (2.2,2);
\draw[very thick] (-4,-2) -- (2.2,-2);
%
%
\draw[redstring] (-2,2) --++ (0,-4);
\draw[redstring] (0,2) --++ (0,-4);
\draw[redstring] (2,2) --++ (0,-4);
\coordinate (g1) at (1.7,0);
\draw[solid] (-2.3,2) node[above,shift={(0,-.05)}]{$\xcoord(\sigma\gamma_1)$} -- (-2.3,0)  -- (1.7,0)  -- (1.7,-2)node[below]{$\xcoord(\gamma_1)$};

\node (g3) at (-2.4,-0.2) {};
\draw[ghost] (-2.5,2) -- (-2.5,-0.2) -- (1.5,-0.2) -- (1.5,-2);
\draw[solid] 
(-4,-.2) -- (-3.5,-0.2) --
 (-3.5,-2) node[below]{$\xcoord(\gamma_3)$};

\draw[solid] (-2.6,2) -- (-2.6,-0.4) -- (1.4,-0.4) -- (1.4,-2);
\draw[solid] (-0.3,2) node[above,shift={(0,-.05)}]{$\xcoord(\sigma\gamma_2)$} -- (-0.3,1.5) -- (-2.1,1.5) -- (-2.1,-1.5)-- (-2.4,-1.5) -- (-2.4,-2) node[below]{$\xcoord(\gamma_2)$};
\draw[solid] (-0.5,2) -- (-0.5,1.7) -- (-2.2,1.7) -- (-2.2,-1.4) -- (-2.6,-1.4) -- (-2.6,-2);
\draw[->] (-4,2.5) -- (-3,2.5) node[midway, above]{$\region{a}{-1}$};
\draw[dotted] (-3,2.5) -- (-3,-2);
\draw[<->] (-2.7,2.5) -- (2.1,2.5) node[midway, above]{$\region{a}{0}$};
\draw[dotted] (-2.7,2.5) -- (-2.7,-2);
\draw[dotted] (2.1,2.5) -- (2.1,-2);

\draw[olive,rounded corners=0pt] (-2.25,-.7) node[below left,shift={(.1,.2)}]{$b)$} rectangle (-2.05,.3);

\draw[olive,rounded corners=0pt] (-2.05,1.4) rectangle (-1.95,1.8);
\draw[olive] (-1.95,1.4) node[below right,shift={(-.05,.05)}]{$c)$};

\node(l1) at (-2,0){};
\node(l2) at (-2,-.4){};
\node(l1p) at (0,0){};
\node(l2p) at (0,-.4){};
{\color{olive}
\node (cross_e) at (-1,.4){$e)$};

\draw[->] (cross_e) -- (l1);
\draw[->] (cross_e) -- (l1p);
\draw[->] (cross_e) -- (l2);
\draw[->] (cross_e) -- (l2p);
}

\node(l3) at (-2,-.2){};
\node(l3p) at (0,-.2){};
{\color{purple}
\node (cross_f) at (-1,-.6){$f)$};
\draw[->] (cross_f) -- (l3);
\draw[->] (cross_f) -- (l3p);
}
\end{tikzpicture}
\end{center}
\caption{Illustration of \ref{P:crossings_Cslam} for a diagram $\Cslam[1]$}
\label{F:crossings}
\end{figure}

\begin{Corollary}
\label{P:Cslamreg_Cslamsing}
There exists a generalised partition $\gbconf_\blam^k$ and diagrams
$\Cslamsing \in \Web(\sigma^k\blam,\gbconf^k_\blam)$ and
$\Cslamreg\in\Web(\gbconf^k_\blam,\blam)$ such that $\Cslam = \Cslamsing \Cslamreg$, and any crossing in $\Cslam$ between a solid or ghost string from $\blamr$ and a red string is contained in $\Cslamsing$.
In
particular, the diagram $\Cslamreg$ is regular.
\end{Corollary}

\begin{proof} Following \autoref{D:Cslam}, define the diagram $\Cslamreg$ to be the subdiagram of $\Cslam$ that is
  below the line $y=H$ and, similarly, define $\Cslamsing$ to be the subdiagram that is above this line. By \autoref{P:Cslamreg_Cslamsing}, these two diagrams satisfy the requirements of the corollary.
\end{proof}

%

Let $\blam \in \Partss$, let $\t \in \Stds$ and fix $k \in \set{1, \dots, p-1}$. Recall that we constructed a regular diagram $\Btreg$ in \autoref{L:BTreg}. The diagram $\Cslamreg \Btreg$ is non-zero but, in general, it can contain generalised double crossings. However, $\Cslamreg$ and $\Btreg$ are both regular diagrams so, by~\autoref{P:resolve_bigon}, there is a regular diagram $\Ckt \in \Web(\gbconf_\t^k,\gbconf_\t)$ that does not have any generalised double crossings such that
\begin{equation}
\label{E:def_Csigmakblam_sigmakt}
  \Ckt = \Cslamreg \Btreg \quad\text{  in } \WA(\gbconf_\t^k,\gbconf_\t).
\end{equation}

\begin{Lemma}
\label{L:red_crossings_Csigmaktsigmakblam}
  Let $\t\in\Stds$, where $\blam\in\Partss$, and  let $1\le k<p$.
Suppose that a solid $i$-string in $\Ckt$ crosses a red $j$-string. Then
$j=\eps^c i$, where $0<c<k$. In particular, there are no solid-red
crossings in $\Ckt$ when $k=1$.
\end{Lemma}

\begin{proof}
By \autoref{P:crossings_Cslam}(d), the only solid-red  crossings in $\Cslamreg$ are between a
solid string $\str_\gamma$, for $\gamma = (r,c,l) \in \blaml \cap \regionblam[0]$, and the
red strings $\rstr_{l+k-p}, \dots, \rstr_{l-1}$, which all belong to
$\region{a}{0}$. By \autoref{L:BTreg} and~\autoref{L:node_region0_between_reds}, in $\Btreg$
the solid string $\str_\gamma$ crosses only the red strings $\rstr_{ap+1}, \dots,
\rstr_{l-1}$. Since $\gamma \in \blaml$ we have $l+k-p
\geq ap+1$. Therefore, since $\Ckt$ has no
generalised double crossings, the only crossings that remain are for the
red strings $\rstr_{ap+1}, \dots, \rstr_{l+k-p-1}$. If $k = 1$ this set
is always empty  since $l-p \leq ap$, thus we now assume $k \geq 2$. We
also assume $l+k-p > ap+1$, since the latter set of red strings is empty
if $l+k-p = ap+1$.  We now write $l = ap + b + 1$ with $0 \leq b < p$.
Since $l + k - p > ap+1$  we have $b+k-p >  0$. Set $\alpha=q^{\kappa_{ap+1}}$,
then the red strings have residues
$\alpha, \eps \alpha, \dots, \eps^{b+k-p- 1} \alpha$. As
$i=\res(\gamma)=\eps^b \alpha$ these
residues can be written as $\eps^{-b} i,\eps^{1-b}i, \dots, \eps^{k-p-1} i$, or as
$\eps^{p-b}i,\eps^{p-b+1},\dots,\eps^{k-1}i$ since  $\eps$ has
order~$p$. As $0\le b<p$ the result follows.
\end{proof}

%
%

\begin{Remark}
The statement of~\autoref{L:red_crossings_Csigmaktsigmakblam} also holds for ghost-red crossings.
\end{Remark}

\begin{Lemma}
\label{L:Cslamsing_Csigmak_no_bigon}
 The diagram $\Cslamsing \Ckt$ has no generalised double crossings.
\end{Lemma}


\begin{proof}
  By construction, neither of the diagrams $\Cslamsing$ and $\Ckt$ has a
  generalised double crossing, so it suffices to prove that any
  crossing that appears in $\Cslamsing$ does not appear in~$\Ckt$. By
  construction, the Webster diagram $\Cslam = \Cslamsing \Cslamreg$ has
  no generalised double crossings, thus it suffices to prove that any
  crossing in $\Cslamsing$ does not appear in $\Btreg$, since
  $\Ckt = \Cslamreg \Btreg$. As $\Btreg$ is
  regular, we only need to consider the regular crossings in $\Cslamsing$.
  By~\autoref{P:crossings_Cslam}(c), the only regular crossings in
  $\Cslamsing$ are crossings of the solid strings $\str_\gamma$, for
  $\gamma=(r,c,l) \in \blamr \cap \regionblam[0]$, with the
  red strings $\rstr_{l+1},\dots,\rstr_{l+k-1}$. On the other hand,
  by~\autoref{L:BTreg} and~\autoref{L:node_region0_between_reds} the
  corresponding solid strings in $\Btreg$ only cross the red strings
  $\rstr_{l'}$ for $l' < l$. This completes the proof.
\end{proof}

\section{A skew cellular basis for
\texorpdfstring{$\Rlpn$}{R(l,p,n)}}\label{S:SkewCellularity}

We are finally ready to define the basis of $\WA(\bom)$ that we
need to prove our main results. Recall from \autoref{SS:Embedding} that
$\tsig$ is the automorphism on the set of Webster diagrams that
multiplies the residues of all solid and ghost strings by~$\eps$ and
recall the definitions of the diagrams $B_\t$, $\Btreg$ and $\Btsing$ from
\autoref{subsection:particular_class_diagrams}. For the readers' convenience, we summarise
the relationships between the different diagrams we defined in the last
chapter:
\begin{align*}
\Btreg&\in\Web(\blam,\gbconf_\t) &\text{(\autoref{L:BTreg})}
\\
\Btsing&\in \Web(\gbconf_\t,\bom) &\text{(\autoref{L:BTsing})}
\\
B_\t  &= \Btreg \Btsing \in \Web(\blam,\bom) &\text{(\autoref{D:RegSing})}
\\
\Cslam &= \Cslamsing \Cslamreg \in \Web(\sigma^k\blam,\blam)& \text{(\autoref{D:Cslam}  and  \autoref{P:Cslamreg_Cslamsing})}
\\
 \Ckt &= \Cslamreg \Btreg \in \WA(\gbconf_\t^k,\gbconf_\t) &\eqref{E:def_Csigmakblam_sigmakt}
\end{align*}

\subsection{A particular cellular basis of \texorpdfstring{$\RG$}{R(Lambda,n)}.}

To show that $\Rlpn$ is a skew cellular algebra we first show that~$\RG$
has a shift-automorphism. To do this we first use the diagrams defined
in \autoref{S:ShiftedRegularity} to define a particular basis of
$\set{C_{\s\t}}$ of $\RG\cong\WA(\bom)$ that is compatible with
\autoref{D:CTdiag}  and then show that
$\sigma(C_{\s\t})=C_{\sigma(\s)\sigma(\t)}$, for all pairs $(\s,\t)$ of
standard tableaux of the same shape.

\begin{Definition}
\label{D:Csigmat}
Suppose that $\blam \in \Partss$, $\t \in \Stds$ and $0\le k<p$. Set
\begin{align*}
  \Ctreg[\sigma^k\t] &\coloneqq \begin{cases*}
      \tsig^k\bigl(\Cslamsing \Ckt\bigr) \in \Web(\sigma^k\blam,\gbconf_\t),
         & if  $k\ne0$,\\
      \Btreg,& if $k=0$,
  \end{cases*}
\\
\intertext{so that $\Ctreg[\sigma^k\t]\in \Web(\sigma^k\blam,\gbconf_\t)$. Define}
\label{equation:Csigmakt}
C_{\sigma^k\t} &\coloneqq \Ctreg[\sigma^k\t] \circ \tsig^{k}\bigl(\Btsing\bigr) \in \Web(\sigma^{k}\blam,\bom).
\end{align*}
\end{Definition}

Notice that the diagrams $\set{C_{\sigma^k\t}|0\le k<p}$ are in a single
$\Z/p\Z$-orbit under $\tsig$. By definition, $C_\t=C_{\sigma^k\t}{}|_{k=0}=B_\t$.

%

\begin{Lemma}
\label{L:Csigmakt_reg}
Let $\blam \in \Partss$, $\t \in \Stds$ and $0\le k<p$. Then the
diagram $\Ctreg[\sigma^k\t]$ is regular.
\end{Lemma}

\begin{proof}
  If $k=0$ then $\Ctreg=\Btreg$ and the result follows from \autoref{P:CT},
  so we can assume that $k>0$. By definition, since $\tsig$ is compatible with the concatenation of Webster diagrams we have
  $\Ctreg[\sigma^k\t]=\tsig^k\bigl(\Cslamsing\bigr)\circ \tsig^k\bigl(\Ckt\bigr)$. It suffices to prove that both $\tsig^k\bigl(\Cslamsing\bigr)$ and  $\tsig^k\bigl(\Ckt\bigr)$ are regular.
  By~\autoref{P:crossings_Cslam}(c), the only crossings in $\Cslamsing$ are
  solid-red (and ghost-red) where the solid string $\str_\gamma$, for
  $\gamma\in\blamr$, has residue $i$ and the red string has residue
  $\eps^c i$ for $0 \leq c< k$.
%
  So, the only crossings in
  $\tsig^k\left(\Cslamsing\right)$ are solid-red (and ghost-red) where
  the solid string has residue $\eps^k i$ and the red string has
  residue $\eps^c i$. In particular,
  $\tsig^k\bigl(\Cslamsing\bigl)$ is a regular diagram since $0 \leq c < k < p$.

  By definition,
  the diagram $\Ckt$ is regular, so to prove that
  $\tsig^k\bigl(\Ckt\bigr)$ is regular it suffices to consider the
  solid-red crossings in $\Ckt$. By
  \autoref{L:red_crossings_Csigmaktsigmakblam}, if a solid
  $\eps^ki$-string crosses a red $j$-string in $\tsig^k\bigl(\Ckt\bigr)$
  then $j=\eps^{c'}i$, where $0<c'<k$. All of these crossing are
  regular, so the lemma is proved.
\end{proof}

%
%

\begin{Proposition}
\label{P:Csigmakt_satisfies_assumptions}
Let $\blam \in \Partss$, $\t \in \Stds$ and $0\le k<p$. Then the
diagram $C_{\sigma^k\t}$ satisfies the assumptions of \autoref{D:CTdiag}.
\end{Proposition}


\begin{proof}
  When $k=0$ the result is just \autoref{P:CT}, so we can assume that
  $0<k<p$. By construction the diagram $C_{\sigma^k\t}$ has no dots and
  a solid string of residue $\eps^k \res(\gamma)=\res(\sigma^k\gamma)$
  from  $(\xcoord[\t](\gamma), 0)$ to $(\xcoord(\sigma^k\gamma), 1)$, for
  all $\gamma\in\blam$. Hence, $C_{\sigma^k\t}$ satisfies parts (a) and
  (b) of \autoref{D:CTdiag} and it remains to verify~(c). That is, we
  need to show that $C_{\sigma^k\t}$ does not contain a generalised
  double crossing. By~\autoref{L:Cslamsing_Csigmak_no_bigon} and
  \autoref{L:BTsing}, respectively, neither of the diagrams
  $\Ctreg[\sigma^k\t]$ and $\tsig^k\bigl(\Btsing\bigr)$ contains a
  generalised double crossing. Therefore, it suffices to prove that the
  diagrams $\Ctreg[\sigma^k\t]$ and $\tsig^k\bigl(\Btsing\bigr)$ do not
  have any crossings in common.

Recall that $\Ctreg[\sigma^k\t] = \tsig^k\bigl(\Cslamsing \Ckt\bigr)$ in
$\Web(\sigma^k\blam,\gbconf_\t)$, where
$\Ckt=\Cslamreg \Btreg \in \WA(\gbconf_\t^k,\gbconf_\t)$. By~\autoref{L:string_Btreg_stays_region}
and~\autoref{P:crossings_Cslam}, any crossing in $\Ctreg[\sigma^k\t]$ is
between two strings (solid, ghost or red) that begin in a same region
$\region{a}{\dia}$. 
On the contrary, by~\autoref{L:BTsing} any  two strings
 in $\Btsing$ that end in a same region $\region{a}{\dia}$ do not cross, hence do not cross in
$\tsig^k\bigl(\Btsing\bigr)$ either. This completes the proof.
\end{proof}

Now we can apply \ref{P:Csigmakt_satisfies_assumptions} to \ref{D:CTdiag}, \ref{D:CST} and \ref{T:CSTBasis}
to get a particular $\Z$-graded cellular basis $\{C_{\s\t}\}$ of $\RG$. This particular $\Z$-graded cellular basis will play a key role in the proof of our main result in the next section.


\subsection{Main results: proof of skew cellularity}
\label{subsection:main_results}

  We are now ready to prove our main theorem from the introduction,
  which says that $\Rlpn$ is a graded skew cellular algebra. As a
  consequence, we deduce that the (ungraded) Iwahori-Hecke algebras of
  type~$D$ are cellular algebras, under weaker assumption than in the literature.


We can now state the main technical result of this paper, which implies all our main results. Recall from~\autoref{D:sigma_cellular} that a shift automorphism of a graded cellular algebra~$A$ is a triple of maps $\sigma=(\sigma_A,\sigma_\P,\sigma_{\Std})$, that satisfies three
requirements, the most important of which is that $\sigma_A(c_{\s\t})=c_{\sigma_{\Std}(\s)\sigma_{\Std}(\t)}$,
for all $\s,\t\in T(\lambda)$ and $\lambda\in\P$. Recall that
definitions of the map $\sigma^\bLam_n$ from
\autoref{sigmaIso} and the maps
$\sigma_{\mathscr{P}}$ and $\sigma_{\Std}$ from \autoref{D:sigmaP}. The
graded cellular structure that we consider on $\RG$ is the one described
in~\autoref{S:Webster}, with the graded cellular basis obtained
from~\autoref{P:Csigmakt_satisfies_assumptions}.


\begin{Theorem}\label{T:sigmaCst}
  The triple of maps
  $\sigma=(\sigma^\bLam_n,\sigma_{\mathscr{P}}, \sigma_{\Std})$ is
  a shift automorphism of~$\RG$.
\end{Theorem}

\begin{proof}
First, note that by~\autoref{nondescrease1} we know that $\sigma_{\mathscr{P}}$ is a poset automorphism of $(\Parts, \tedom)$.   By \autoref{D:sigma_cellular}, we need to show that:
\begin{enumerate}
  \item If $\s\in \Std(\blam)$ then $\sigma_{\Std}(\s)\in \Std(\sigma_{\mathscr{P}}(\blam))$
    and $\deg(\sigma_{\Std}(\s))=\deg(\s)$.
  \item If $\s,\t\in \Std(\blam)$ then $\sigma_{\RG}(c_{\s\t})=c_{\sigma_{\Std}(s)\sigma_{\Std}(\t)}$.
  \item If $\s,\t\in \Std(\blam)$, for $\blam\in\Parts$, then
  $\sigma_{\Std}^k(\t)=\t$ if and only if $\sigma_{\Std}^k(\s)=\s$, for $k\in\Z$.
\end{enumerate}
The first requirement in part~(a) is immediate from \autoref{D:sigmaP}
and the second requirement, that the fact that $\sigma_{\Std}$ is homogeneous
follows from~\autoref{L:SigmaOnTableaux}. Part~(c) follows because all
tableaux have order~$p$ under the action of $\sigma_{\Std}$.
It remains to check part~(b). That is, we need to show that
\[
\sigma\left(C_{\sigma^k\s,\sigma^l\t}\right) = C_{\sigma^{k+1}\s,\sigma^{l+1}\t},
\qquad\text{for any  $\blam \in \Partss$, $\s,\t \in \Stds$ and $k,l\in\Z/p\Z$,}
\]
Equivalently, we need to show that
$
\sigma\left(C_{\sigma^k\s}^*C_{\sigma^l\t}\right) =
C_{\sigma^{k+1}\s}^*\:C_{\sigma^{l+1}\t}.
$
By~\autoref{D:Csigmat} and \autoref{L:tsig_star_commute},
\begin{align*}
C_{\sigma^k\s}^* &= \tsig^k\left(\Btsing[\s]\right)^* \left(\Ctreg[\sigma^k\s]\right)^*,
&
C_{\sigma^l\t} &=\Ctreg[\sigma^l\t]\: \tsig^l\left(\Btsing\right),
\\
C_{\sigma^{k+1}\s}^* &= \tsig^{k+1}\left(\Btsing[\s]\right)^* \left(\Ctreg[\sigma^{k+1}\s]\right)^*,
&
C_{\sigma^{l+1}\t} &=\Ctreg[\sigma^{l+1}\t]\: \tsig^{l+1}\left(\Btsing\right).
\end{align*}
By~\autoref{L:string_Btreg_stays_region} and~\autoref{L:string_Cslam_stays_region}, any string in ${(\Ctreg[\s])}^* \Ctreg \in
\Web(\gbconf_\s,\gbconf_\t)$ stays inside a single region. Moreover, by \autoref{L:min_interval_not_red}, any solid or ghost string  of ${(\Ctreg[\s])}^* \Ctreg \in
\Web(\gbconf_\s,\gbconf_\t)$ in the region $\region{a}{0}$  starts (and ends)  at the left of all the red strings in $\region{a}{0}$, for $0 \le a < d$.
 By \autoref{D:Csigmat} and~\eqref{E:def_Csigmakblam_sigmakt}, this is also true in the two
diagrams $\left(\Ctreg[\sigma^k \s]\right)^* \Ctreg[\sigma^l\t]$ and
$\left(\Ctreg[\sigma^{k+1} \s]\right)^* \Ctreg[\sigma^{l+1}\t]$. As both
of these diagrams are regular by~\autoref{L:Csigmakt_reg}, we can apply
\autoref{P:resolve_bigon} to find two regular Webster diagrams
$\Dkl$ and $\Dkl[+1]$ that do not contain any generalised double
crossings such that in $\WA(\gbconf_{\sigma^k\s},\gbconf_{\sigma^l\t})$
\begin{align*}
\Dkl =\left(\Ctreg[\sigma^k \s]\right)^* \Ctreg[\sigma^l\t]
\qquad\text{and}\qquad
\Dkl[+1]=\left(\Ctreg[\sigma^{k+1} \s]\right)^* \Ctreg[\sigma^{l+1}\t]
\end{align*}
and all of the solid and ghost strings in any region $\region{a}{0}$ are
to the left of all of the red strings in~$\region{a}{0}$.
 Since $\Dkl$ is a regular diagram that has no crossings involving red strings, the diagram $\tsig(\Dkl)$ is also regular. Hence, the diagrams $\tsig(\Dkl)$ and $\Dkl[+1]$ satisfy the assumptions of \autoref{C:begin_end_points}, and so
$\tsig(\Dkl) = \Dkl[+1]$ in~$\WA(\gbconf_{\sigma^k\s},\gbconf_{\sigma^l\t})$.
We have proved so far that
\begin{align*}
C_{\sigma^k\s,\sigma^l\t} &= \tsig^k\left(\Btsing[\s]\right)^* \Dkl \:\tsig^l\left(\Btsing\right),
\\
C_{\sigma^{k+1}\s,\sigma^{l+1}\t} &= \tsig^{k+1}\left(\Btsing[\s]\right)^* \tsig\left(\Dkl\right) \tsig^{l+1}\left(\Btsing\right).
\end{align*}
Consequently,
\begin{equation}
\label{E:almost_final}
\tsig\left(C_{\sigma^k\s,\sigma^k\t}\right) = C_{\sigma^{k+1}\s,\sigma^{l+1}\t}.
\end{equation}
Recall that in $\Dkl$, all of the solid and ghost strings that start or, equivalently, end in~$\region{a}{0}$ are to the left of all of the red strings in $\region{a}{0}$. Similarly, since $\left(\Ctsing[\s]\right)^* \Ctsing \in \Web(\bom,\bom)$, all the solid strings begin and end to the left of all red strings. Moreover by~\autoref{L:BTsing}\autoref{I:CTsing_a_ignoring}, in this diagram any solid string that crosses a red string in~$\region a 0$ must cross all of the red strings in $\region a 0$, for  $0\le a<d$.  Thus, we can apply \autoref{P:diagram_crosses_group_p_strings}, which gives
\[
\tsig\left(C_{\sigma^k\s,\sigma^l\t}^\cyc\right) = \tsig\left(C_{\sigma^k \s,\sigma^l \t}\right)^\cyc.
\]
Combining \eqref{E:almost_final} with \autoref{L:acyc} and~\autoref{P:SigmaDiagram}, we find that
\[
\sigma\left(C_{\sigma^k\s,\sigma^k\t}\right) = C_{\sigma^{k+1}\s,\sigma^{l+1}\t}^{\cyc} = C_{\sigma^{k+1}\s,\sigma^{l+1}\t},
\]
as desired.
\end{proof}


\begin{Remark}
A very particular case of~\autoref{T:sigmaCst} can be found in~\cite[Example 7.5]{Bowman:ManyCellular}. Namely, for $\ell = e = 2$ and the $2$-charge obtained from $(0,1)$,  the basis elements are of the form
\begin{align*}
c_{\mathfrak{s},\mathfrak{s}} &= e(\bi) & c_{\sigma\mathfrak s, \sigma\mathfrak s}&=e(-\bi)
\\
\cline{1-4}
c_{\mathfrak{t},\mathfrak t} &= y_2 e(\bi)
&
c_{\sigma\mathfrak t, \sigma\mathfrak t} &= y_2 e(-\bi)
\\
c_{\sigma\mathfrak t , \mathfrak t} &= \psi_1 e(\bi)
&
c_{\mathfrak t ,\sigma\mathfrak{t}} &= \psi_1 e(-\bi)
\\
\cline{1-4}
c_{\mathfrak u, \mathfrak u} &= y_2^2 e(\bi)
&
c_{\sigma\mathfrak u,\sigma\mathfrak u} &= y_2^2 e(-\bi)
\end{align*}
where $\mathfrak{s}$ (resp. $\mathfrak{u}$) is the only element in $\Std(1^2|0)$ (resp. $\Std(2|0)$), $\bi \in \{(1,-1),(-1,1)\}$ and the tableau $\mathfrak{t}$ is a particular element of $\Std(1|1)$. Recall that $\sigma$ fixes $y_1,y_2,\psi_1$ and that $\sigma\bigl(e(\bj)\bigr) = e(-\bj)$.
\end{Remark}


We can now prove our \hyperlink{MainTheorem}{Main Theorem} from the introduction.

\begin{Corollary}
\label{C:quiver_hecke_graded-skew-cellular}
Assume that $R$ contains a primitive $p$-th root of unity and $p\cdot 1_R$ is invertible in $R$. The algebra  $\Rlpn$ is a graded skew cellular algebra. In particular, if $2 \cdot 1_R \in R^\times$ then $\RR^{\bLam}_{2d,2,n}$ is graded cellular.
\end{Corollary}

\begin{proof} 
  This is an immediate consequence of~\autoref{T:sigmaCst},~\autoref{C:skew_cellular_cellular} and~\autoref{proposition:sigma_cellular_implies_extended} applied for the graded cellular algebra $\RG$. Note that $\sigma^\bLam_n$ has order $p$ indeed.
\end{proof}

Combining the last result with \autoref{P:unitriangular} yields:

\begin{Corollary}
  The graded decomposition matrix of $\Rlpn$ is unitriangular.
\end{Corollary}


We now explicitly describe the graded skew cellular datum of $\Rlpn$,
following the construction given in
\autoref{proposition:sigma_cellular_implies_extended}.
Recall that $\Partss$ is any set of representatives for the action of $\langle\sigma_{\mathscr{P}}\rangle$ on $\Parts$. We endow $\Partss$ with the partial order $\tedom$ so that $a \tdom b$ for $a,b \in \Partss$ if and only if there exists $\blam \in a$ and $\bmu \in b$ such that $\blam \tdom \bmu$.
Define
\begin{equation}\label{E:Partsp}
\Partsp \coloneqq \Partss \times \mathbb{Z}/\so_\blam \mathbb{Z},
\end{equation}
where $\so_\blam$ is the size of the orbit of $\blam$ under the action of $\langle\sigma_{\mathscr{P}}\rangle$. The order $\tedom$ extends to an order of $\Partsp$ with the rule $(\blam,k) \tdom (\bmu,l) \iff \blam \tdom \bmu$.
Let $\iota$ be the involution on $\Partsp$ given by $(\blam,k) \mapsto (\blam,-k)$. For $(\blam,l)\in\Partsp$ let
$\Std_\sigma(\blam,k) \coloneqq \Std_\sigma(\blam)$.
For $(\blam,k)\in\Partsp$ and
$\s,\t\in\Std_\sigma(\blam,k)$, define
\[
D^{(k)}_{\s\t} \coloneqq
\sum_{j = 0}^{\sp_\blam - 1} \eps^{kj\so_\blam} \sigmavg (C_{\s,\sigma^{j\so_\blam} \t}),\]
where
\[
    \sigmavg  \coloneqq  \sum_{k = 0}^{p-1} \sigma^k\quad
    \text{and}\quad\sp_\blam=\frac{p}{\so_\blam}.
\]
In particular, note that $\olambdatab[\blam] = \sp_\blam$, since $\sp = \sp_\P = p$, so $\frac{p}{\olambdatab[\blam]} = \frac{p}{\sp_\blam} = \so_\blam$ and $\epsl[\blam] = \eps^{\so_\blam}$. Finally, define the map
\[
D\map{ \Std_\sigma(\blam,k)\times \Std_\sigma(\blam,k)}{\Rlpn};(\s,\t)\mapsto D_{s\t}^{(k)},
\]
and let $\iota_{\blam}^{(k)}$ be the identity map on $\Std_\sigma(\blam,k)$. Then
\begin{equation}
\label{E:cellular_datum_grpn}
\bigl(\Partsp, \iota, \Std_\sigma, D, \deg\bigr),
\end{equation}
is a graded skew cell datum for   $\Rlpn$. In particular, $\RG$ has cell
modules $C_{\blam}^{(k)}$, and simple modules the non-zero quotients
$D_\blam^{(k)}=C_\blam^{(k)}/\rad C_{\blam}^{(k)}$, for
$(\blam,k)\in\Partsp$.

Note that the graded skew cell datum for $\Rlpn$ restricts to give a graded
skew cell data for the blocks $\Rlpn[\alpha]$ and that the graded
decomposition matrix of $\Rlpn[\alpha]$ is unitriangular, for
$\alpha\in\cIsig$.

Recall from \autoref{D:ICompositions} that
  $\cI=\set{\eps^j \q^i|0\leq j<p \text{ and } 0\leq i<e}$. Define
  \[
    z \coloneqq \sum_{\bi \in \cI^n} i_1^{-1} e(\bi) \in \RG.
  \]

\begin{Lemma}
\label{L:existence_z}
The element $z\in\RG $ is homogeneous of degree $0$, invertible and $\sigma( z) = \eps  z$.
\end{Lemma}

\begin{proof}
First, recall from~\autoref{SS:CycHecke} that $e$ is finite thus $\cI$ is finite and $ z$ is well defined.
We have $\deg  z = 0$ since $\deg e(\bi) = 0$ for all $\bi \in \cI^n$. Using the fact that $\set{e(\bi)  : \bi \in \cI^n}$ is a complete set of orthogonal idempotents we find that $ z$ is invertible in $\RG$ with inverse $\sum_{(i_1,\dots,i_n) \in \cI^n} i_1 e(i_1,\dots,i_n)$. Finally, we have
\begin{align*}
\sigma( z)
&=
\sum_{(i_1,\dots,i_n) \in \cI^n} i_1^{-1} e(\eps i_1,\dots,\eps i_n)
\\
&=
\sum_{(i_1,\dots,i_n) \in \cI^n} \bigr(\eps^{-1}i_1\bigr)^{-1} e(i_1,\dots, i_n)
\\
&=
\eps  z,
\end{align*}
which completes the proof.
\end{proof}

\begin{Corollary}\label{C:diagonalisable}
  The automorphism $\sigma$ is $\varepsilon$-splittable in the sense of
  \autoref{E:z_strong}.
\end{Corollary}

Hence, all the results in~\autoref{subsection:Clifford} now apply to $\Rlpn$. Note that multiplication by $ z$ induces an analogue of \autoref{L:existence_z} for the blocks $\RG[\alpha]$ of~$\RG$, so the results in \autoref{subsection:Clifford} also apply to the blocks~$\Rlpn[\alpha]$ of~$\Rlpn$, for $[\alpha]\in\cIsig$ .



By~\autoref{rpnIso},   all the results concerning the skew cellularity of $\Rlpn$, over a ring~$R$ containing $\eps$, can be deduced from $\Hlpn$ over the field $K$, which contains $\eps$ by assumption. For example, the following holds:

\begin{Corollary}
\label{C:hecke_graded-skew-cellular}
The algebra  $\Hlpn$ is a graded skew cellular algebra. Moreover, if $p = 2$ then $\HH_{2,n}(q,\bvQ)$ is a graded cellular algebra.
\end{Corollary}


Geck~\cite{Geck:cellular} has proved that the Iwahori--Hecke algebra of a finite Coxeter group is always an (ungraded) cellular algebra. By~\autoref{C:hecke_graded-skew-cellular}, the Iwahori--Hecke algebras of the Coxeter groups of types $A_{n-1}$, $B_n=C_n$, $D_n$ and $I_2(n)$ are graded skew cellular algebras. (These Coxeter groups are the complex reflection groups of types $G(1,1,n), G(2,1,n), G(2,2,n)$ and $G(n,n,2)$, respectively.)

The Iwahori-Hecke algebras of types $A_{n-1}$ and $B_n$ are graded
cellular algebras by \cite{HuMathas:GradedCellular}. For the algebras of
types $D_n$ and $I_2(n)$, Geck's proof of the
cellularity of the Iwahori-Hecke algebras $\mathscr{H}_\q(D_n)$ assumes that
$\q^{1/2}\in K$. The following corollary generalises Geck's result to
the graded case and removes the assumption that $\q^{1/2}\in K$.


\begin{Corollary} Suppose $\cha K\neq 2$. The Iwahori-Hecke algebra $\mathscr{H}_\q(D_n)$ of
  type $D_n$ is a $\Z$-graded cellular algebra.
\end{Corollary}

\begin{proof} The Coxeter group of type $D_n$ is the complex reflection group of type $G(2,2,n)$ in the Shephard--Todd classification. The result thus directly follows from~\autoref{C:hecke_graded-skew-cellular} for $d =1$.
\end{proof}


It is an interesting open question whether the algebras $\Rlpn$ are
graded cellular algebras when $p>2$ (and $n > 2$ if $p = \ell$). 

\subsection{Adjustment matrices}
\label{subsection:adjustment}

As a final application we use the bases in this paper to describe
``adjustment matrices'' for the Hecke algebras of type $G(\ell,p,n)$,
which relate decomposition matrices in different characteristics. For
the Iwahori-Hecke algebras of finite Coxeter groups Geck and Rouquier
used Lusztig's asymptotic Hecke algebra to show that adjustment
matrices exist whenever the Iwahori-Hecke algebra is defined over a field of
``good characteristic'', which depends on the root system of the
underlying Coxeter group; see \cite[Table 1.4 and
Theorem~3.6.3]{GeckJacon:Book}.  When $p=1$ we recover the results of
\cite[\Sec5.6]{BK:GradedDecomp} but, even in the ungraded case, these
results appear to be new when $p>2$.


Throughout this section fix $\bLam$ as in~\autoref{subsection:QHA}.  If $R$ is an integral
domain let $\Rlpn(R)$ be the quiver Hecke algebra over $R$ with
weight~$\bLam$. Similarly, if $M$ is an $\Rlpn(R)$-module then we
write $M=M^R$ to emphasize that $M$ is an $R$-module. By
\autoref{C:quiver_hecke_graded-skew-cellular}, if $R$ contains a
primitive $p$th root of unity and $p \in R^\times$ then $\Rlpn(R)$ is a
skew cellular algebra with cell modules $C_{\blam}^{(k),R}$, for
$(\blam,k)\in\Partsp$. By \autoref{T:SkewCellularSimples}, if $R$ is a
field, the graded simple $\Rlpn(R)$-modules are  shifts of the non-zero
quotients $D_{\blam}^{(k),R}=C_{\blam}^{(k),R}/\rad C_{\blam}^{(k),R}$, for
$(\blam,k)\in\Partsp$.


Fix a field $F$ of characteristic $c>0$ that contains a primitive $p$th
root of unity $\eps_F$. In particular, this implies that $c$ does not
divide $p$. The subfield $F_{p,c} = \Z/c\Z(\eps_F)$ of $F$ is a splitting field for the $p$-th cyclotomic polynomial $\Phi_p$ over $\Z/c\Z$.   By~\cite[Theorem~2.47]{LidlNiederreiter}, the field $F_{p,c}$ is the finite field with $c^r$ elements, where $r>0$
is minimal such that $c^r\equiv 1\pmod p$.

Let $\pi : \Z[\eps] \to F_{p,c}$ be the unique ring homomorphism determined by $\pi(\eps) = \eps_F$. Then $\mathfrak{p} = \ker \pi$ is a prime ideal of $\Z[\eps]$ and the localisation $\Ocal=\Z[\eps]_{\mathfrak{p}}$  of $\Z[\eps]$
at the prime ideal~$\mathfrak{p}$ is a discrete valuation ring with maximal ideal  $\mathfrak{p}\Ocal$. The next result follows from elementary properties of localisation.

\begin{Lemma}
The residue field $\Ocal / \mathfrak{p}\Ocal$ of $\Ocal$ is isomorphic to $F_{p,c}$.
\end{Lemma}

Consequently, if~$M^\Ocal$ is an~$\Rlpn(\Ocal)$-module then, by base
change, $M^{F_{p,c}}=F_{p,c}\otimes_\Ocal M^\Ocal$ is
an~$\Rlpn(F_{p,c})$-module and~$M^F=F\otimes_{F_{p,c}}M^{F_{p,c}}\cong
F\otimes_\Ocal M^\Ocal$ is an~$\Rlpn(F)$-module.

\begin{Remark}
  By \autoref{T:SkewCellularSimples} every field is a splitting field
  for the algebra $\Rlpn$, so there would be no loss of generality in assuming
  that $F=F_{p,c}$.
\end{Remark}

Let $\Qcal = \mathrm{Frac}(\Ocal) = \Q(\eps)$ be the field of fractions of $\Ocal$.
By \autoref{D:CellModule}, $C_{\blam}^{(k),\Qcal}\cong\Qcal\otimes_\Ocal
C_{\blam}^{(k),\Ocal}$ and
$C_{\blam}^{(k),F}\cong F\otimes_\Ocal C_{\blam}^{(k),\Ocal}$, for
$(\blam,k)\in\Partsp$. Now recall from \autoref{E:CellRadical} that
\[
  \rad C_{\blam}^{(k),\Ocal} =
      \set[\big]{y\in C_{\blam}^{(k),\Ocal} |\phi_\blam^{(k)}(x, y)=0
           \text{ for all }x \in C_{\blam}^{(k),\Ocal}}
\]
and that
$D_{\blam}^{(k),\Ocal} = C_{\blam}^{(k),\Ocal}/\rad C_{\blam}^{(k),\Ocal}$ by \autoref{D:CellularSimples}.

%

\begin{Lemma}\label{L:Efree}
  Let $(\blam,k)\in\Partsp$. Then $D_{\blam}^{(k),\Ocal}$ is an
  $\Ocal$-free graded $\Rlpn(\Ocal)$-module. Moreover,
  $D_{\blam}^{(k),\Qcal}\cong \Qcal\otimes_\Ocal D_{\blam}^{(k),\Ocal}$ as
  graded $\Rlpn(\Qcal)$-modules.
\end{Lemma}

\begin{proof}
  Let $\set{\Dkt|\t\in\Std_\sigma(\blam,k)}$ be the standard basis of
  the cell module $C_{\blam}^{(k),\Ocal}$.
  Since $\rad C_{\blam}^{(k),\Ocal}$ is a pure $\Ocal$-submodule of $C_{\blam}^{(k),\Ocal}$,
    $D_{\blam}^{(k),\Ocal}$ is free as an $\Ocal$-module.
  The final claim that
  $D_{\blam}^{(k),\Qcal}\cong \Ocal\otimes_\Ocal D_{\blam}^{(k),\Ocal}$
  is immediate because
  $\rad C_{\blam}^{(k),\Qcal}\cong\Qcal\otimes_\Ocal \rad C_{\blam}^{(k),\Ocal}$
  since~$\Qcal$ is a field of fractions of~$\Ocal$.
\end{proof}

Let $\Partspz=\set{(\blam,k)\in\Partsp|D_{\blam}^{(k),\Qcal}\ne0}$. For
convenience, set
\[
E_{\blam}^{(k),F}=F\otimes_\Ocal D^{(k),\Ocal}_{\blam},
  \qquad\text{ for }(\blam,k)\in\Partsp.
\]
In general, if $(\blam,k)\in\Partsp$ then the $\Rlpn(F)$-modules
$D_{\blam}^{(k),F}$ and~ $E_{\blam}^{(k),F}$, and
$\rad C_{\blam}^{(k),F}$ and~$F\otimes_\Ocal\rad C_{\blam}^{(k),\Ocal}$,
are \textit{not} isomorphic. In particular, the $\Rlpn(F)$-module
$E_{\blam}^{(k),F}$ is not necessarily irreducible.

Suppose that $R$ is a field and let $M^R$ be an $\Rlpn(R)$-module.
Recall from \autoref{E:DecompNum} that if~$(\bmu,l)\in\Partspz$ and
then the graded decomposition multiplicity
of $D_{\bmu}^{(l),R}$ in $M^R$ is
\[
  [M^R: D_{\bmu}^{(l),R}]_t = \sum_{s\in\Z}\, [M^R: D_{\bmu}^{(l),R}\<s\>]\,t^s
    \quad\in\N[t,t^{-1}].
\]
For $(\blam,k)\in\Partsp$ and $(\bmu,l), (\bnu,m)\in\Partspz$
define graded decomposition numbers
\begin{align*}
  d^\Qcal_{(\blam,k)(\bmu,l)}(t)
      & =[C_{\blam}^{(k),\Qcal}: D_{\bmu}^{(l),\Qcal}]_t,\\
  d^F_{(\blam,k)(\bmu,l)}(t)
      & =[C_{\blam}^{(k),F}: D_{\bmu}^{(l),F}]_t,\\
  \alpha^{F}_{(\bmu,l)(\bnu,m)}(t)
      & =[E_{\bmu}^{(l),F}: D_{\bnu}^{(m),F}]_t.
\end{align*}
Let
$D_{\Qcal}(t)=\bigl(d^\Qcal_{(\blam,k)(\bmu,l)}(t)\bigr)$,
$D_{F}(t)=\bigl(d^F_{(\blam,k)(\bmu,l)}(t)\bigr)$ and
$A_{F}(t)=\bigl(\alpha^F_{(\bmu,l)(\bnu,m)}(t)\bigr)$
be the corresponding matrices, where
$(\blam,k)\in\Partsp$ and $(\bmu,l), (\bnu,m)\in\Partspz$.
Then $D_\Qcal(t)$ and $D_F(t)$ are the graded decomposition matrices of
$\Rlpn(\Qcal)$ and $\Rlpn(F)$, respectively. The matrix $A_F(t)$ is the
\emph{graded adjustment matrix} for these two algebras.

Let $\Rep\Rlpn(\Qcal)$ and $\Rep\Rlpn(F)$ be the Grothendieck groups of
finitely generated graded $\Rlpn(\Qcal)$- and $\Rlpn(F)$-modules,
respectively. Both Grothendieck groups are free $\Z[t,t^{-1}]$-modules
where~$t$ acts by grading shift. If~$M$ is a module for one of these
algebras let $[M]$ be the image of~$M$ in the corresponding Grothendieck
group.

\begin{Lemma}\label{L:Lattice}
  There is a unique abelian group homomorphism
   $d^\Qcal_F\map{\Rep\Rlpn(\Qcal)}{\Rep\Rlpn(F)}$ such that
   $d^\Qcal_F[C_{\blam}^{(k),\Qcal}]=[C_{\blam}^{(k),F}]$ and
   $d^\Qcal_F[D_{\blam}^{(k),\Qcal}]=[E_{\blam}^{(k),F}]$,
   for $(\blam,k)\in\Partsp$.
\end{Lemma}

\begin{proof}
  Since $\Ocal$ is a discrete valuation ring, for any
  $\Rlpn(\Qcal)$-module $M^\Qcal$ there exists an $\Rlpn(\Ocal)$-module
  $M^\Ocal$, a full $\Ocal$-lattice, such that
  $M^\Qcal\cong\Qcal\otimes_\Ocal M^\Ocal$. Define $M^F=F\otimes_\Ocal
  M^\Ocal$. The choice of $\Ocal$-lattice is not unique but~$[M^F]$ is
  independent of the choice of $\Ocal$-lattice; compare
  with \cite[Proposition~16.16]{CurtisReiner:VolI}. Hence, we define
  $d^\Qcal_F[M^\Qcal]=[M^F]$. By \autoref{D:CellModule}, if
  $(\bnu,m)\in\Partsp$ then
  $d^\Qcal_F[C_{\blam}^{(k),\Qcal}]=[C_{\blam}^{(k),F}]$. Moreover, we have
  $d^\Qcal_F[D_{\blam}^{(k),\Qcal}]=[E_{\blam}^{(k),F}]$ because
  $E_{\blam}^{(k),\Qcal}\cong D_{\blam}^{(k),\Qcal}$ by \autoref{L:Efree}.
  Finally, this establishes the uniqueness of $d^\Qcal_F$ since
  $\set{D_{\blam}^{(k),\Qcal}|(\blam,k)\in\Partspz}$ is a basis of
  $\Rep\Rlpn(\Qcal)$.
\end{proof}

Let $\overline{\phantom{X}}\map{\Z[t,t^{-1}]}\Z[t,t^{-1}]$ be the unique $\Z$-linear
map such that $\overline{t^k}=t^{-k}$, for $k\in\Z$. Observe that $\dim_t M^*=\overline{\dim_t M}$ if $M$ is a graded module. Hence, we can extend
$\overline{\phantom{X}}$ to a map of the Grothendieck $\Rep\Rlpn(\Qcal)$
and $\Rep\Rlpn(F)$ by setting $\overline{[M]} =[M^*]$.


\begin{Proposition}\label{P:adjustment}
  Let $(\blam,k)\in\Partsp$ and $(\bmu,l)\in\Partspz$. Then
  \begin{enumerate}
    \item $\alpha^F_{(\bmu,l)(\bmu,l)}(t)=1$,
    \item $\alpha^F_{(\blam,k)(\bmu,l)}(t)\ne0$ only if $(\blam,k)\tedom(\bmu,l)$,
    \item $\alpha^F_{(\blam,-k)(\bmu,-l)}(t)=\overline{\alpha^F_{(\blam,k)(\bmu,l)}(t)}$.
  \end{enumerate}
  Moreover, $D_{\bmu}^{(l),F}\ne0$ if and only if $(\bmu,l)\in\Partspz$.
\end{Proposition}

\begin{proof}
  By construction, $F\otimes_\Ocal\rad C_{\blam}^{(k),\Ocal}$ is an
  $\Rlpn(F)$-submodule of $\rad C_{\blam}^{(k),F}$, so parts~(a) and~(b)
  are immediate from \autoref{P:unitriangular}. For part~(c), in
  $\Rep\Rlpn(F)$ we have
  \begin{align*}
        [E_{\blam}^{(k),F}] &= \sum_{(\bmu,l)\in\Partspz}
           \alpha^F_{(\blam,k)(\bmu,l)}(t) [D_{\bmu}^{(l),F}].
    \intertext{Taking duals and using \autoref{P:DualSimples} this becomes}
        \overline{[E_{\blam}^{(k),F}]} &= \sum_{(\bmu,l)\in\Partspz}
        \overline{\alpha^F_{(\blam,k)(\bmu,l)}(t)}[D_{\bmu}^{(-l),F}].
  \end{align*}
  On the other hand,
  $\overline{[E_{\blam}^{(k),F}]}
  =\overline{d^\Qcal_F[D_{\blam}^{(k),\Qcal}]}
     =d^\Qcal_F\overline{[D_{\blam}^{(k),\Qcal}]}
     =d^\Qcal_F[D_{\blam}^{(-k),\Qcal}]
     =[E_{\blam}^{(-k),F}]$,
  where the second last equality comes from the same argument used in the proof of \cite[16.16]{CurtisReiner:VolI} and \autoref{P:DualSimples}, and
  the last equality from \autoref{L:Lattice}. This proves~(c).


  Finally, we prove that $D_{\blam}^{(k),F}\ne0$ if and only if
  $(\blam,k)\in\Partspz$. First, observe that if~$D_{\blam}^{(k),F}\ne0$
  then $\rad C_{\blam}^{(k),F}\ne C_{\blam}^{(k),F}$, implying that
  $\rad C_{\blam}^{(k),\Qcal}\ne C_{\blam}^{(k),\Qcal}$. Hence,
  $D_{\blam}^{(k),\Qcal}\ne0$ and $(\blam,k)\in\Partspz$. Conversely,
  if $(\blam,k)\in\Partspz$ then $D_{\blam}^{(k),\Qcal}\ne0$ so
  that $E_{\blam}^{(k),F}\ne0$. Hence, $E_{\blam}^{(k),F}\ne0$ by part~(a).
\end{proof}

We can now prove the main result of this section.

\begin{Theorem}\label{T:Adjustment}
  Suppose that $F$ is a field of characteristic $c$ that $F$
  contains a primitive $p$th root of unity $\eps_F$.
  Let $(\blam,k)\in\Parts$ and $(\bmu,l)\in\Partspz$.  Then
  \[
        d^F_{(\blam,k)(\bmu,l)}(t) =
             \sum_{(\bnu,m)\in\Partspz}
                 d^\Qcal_{(\blam,k)(\bnu,m)}(t)
                 \alpha^F_{(\bnu,m)(\bmu,l)}(t).
  \]
  That is, $D_{F}(t) = D_{\Qcal}(t)A_{F}(t)$.
\end{Theorem}

\begin{proof}
  Using \autoref{L:Lattice}, and \autoref{P:unitriangular} twice,
  \begin{align*}
    \sum_{(\bmu,l)\in\Partspz} d^F_{(\blam,k)(\bmu,l)}(t)[D_{\bmu}^{(l),F}]
       &= [C_{\blam}^{(k),F}] = d^\Qcal_F[C_{\blam}^{(k),\Qcal}]\\
       & =d^\Qcal_F\Bigl(\sum_{(\bnu,m)\in\Partspz}
            d^\Qcal_{(\blam,k)(\bnu,m)}(t) [E_{\bnu}^{(m),\Qcal}]\Bigr)\\
       & =\sum_{(\bnu,m)\in\Partspz}
            d^\Qcal_{(\blam,k)(\bnu,m)}(t) [E_{\bnu}^{(m),F}]\\
       & =\sum_{(\bnu,m)\in\Partspz} d^\Qcal_{(\blam,k)(\bnu,m)}(t)
       \sum_{(\bmu,l)\in\Partspz}\alpha^F_{(\bnu,m)(\bmu,l)}(t)
            [D_{\bmu}^{(l),F}]\\
       & =\sum_{(\bmu,l)\in\Partspz} \Bigl(\sum_{(\bnu,m)\in\Partspz}
          d^\Qcal_{(\blam,k)(\bnu,m)}(t)\alpha^F_{(\bnu,m)(\bmu,l)}(t)
            \Bigr)[D_{\bmu}^{(l),F}].
  \end{align*}
  Since $\set{[D_{\bmu}^{(l),F}]|(\bmu,l)\in\Partspz}$ is a basis of
  $\Rep\Rlpn(F)$, the result follows by comparing the coefficient of
  $[D_{\bmu}^{(l),F}]$ on both sides of this equation.
\end{proof}


Let $\Fcal$ be a field of characteristic zero that contains a
primitive $p$th root of unity $\eps$ and let $\Hlpn<\Fcal>$ be a Hecke
algebra of type $G(\ell,p,n)$ over $\Fcal$.  Similarly, let $\Hlpn<F>$
be the corresponding Hecke algebra over the field~$F$, which contains a
primitive $p$th root of unity $\eps_F$.  By \autoref{rpnIso},
$\Hlpn<F>\cong\Rlpn(F)$, so $D_F=D_F(1)$ is the decomposition matrix of
$\Hlpn<F>$, where $D_F(1)$ is the graded decomposition matrix of
$\Rlpn(F)$ evaluated at $t=1$. By \autoref{T:SkewCellularSimples}(b),
every field is a splitting field for $\Rlpn(\Qcal)$,
and $\Rlpn(\Fcal)\cong\Fcal\otimes_{\Qcal}\Rlpn(\Qcal)$,
so $D_\Qcal(t)$ is the graded decomposition matrix of $\Rlpn(\Fcal)$. On the other hand,
$\Hlpn<\Fcal>\cong\Rlpn(\Fcal)$ by \autoref{rpnIso}, so
$D_\Fcal=D_\Qcal(1)$ is the decomposition matrix of $\Hlpn<\Fcal>$.
Finally, set $A_F=A_F(1)$. Hence, by \autoref{T:Adjustment} we obtain
the following.

\begin{Corollary}
  Suppose that $\Fcal$ is a field of characteristic~$0$ containing a
  primitive $pth$ root of unity~$\eps$ and that $F$ is a field of
  characteristic $c>0$ that contains a primitive $p$th root of unity
  $\eps_F$.  Then the decomposition matrix $D_F$
  of $\Hlpn<F>$ factorises as $D_F=D_\Fcal A_F$.
\end{Corollary}

\begin{Remark}
  This section does not really use the machinery of skew cellular
  algebras. Rather the results in this section follow from the fact that
  the skew cellular basis of
  \autoref{C:quiver_hecke_graded-skew-cellular} is defined over the ring
  $\Z[\eps]$, which makes it easy to apply standard modular reduction
  arguments. The existence of adjustment matrices usually requires a
  delicate choice of modular system. The beauty of using KLR algebras is
  that we can work over $\Z[\eps]$, which makes this result almost
  trivial.
\end{Remark}

\subsection{Graded simple modules}\label{SS:Simples}
Let $F$ be a field.
The algebra $\Rlpn(F)$ is a skew cellular algebra by
\autoref{C:quiver_hecke_graded-skew-cellular}, so
\(\set{D_{\blam}^{(k),F}\<s\>|(\blam,k)\in\Partspz\text{ and }s\in\Z}\)
is a complete set of pairwise non-isomorphic graded simple modules by
\autoref{T:SkewCellularSimples} (and \autoref{P:adjustment}). The aim of
this section is to explicitly describe the set~$\Partspz$.

By definition, $\Partspz =\set{(\blam,k)\in\Partsp|D_{\blam}^{(k),F}\ne0}$,
where $\Partsp=\Partss \times \mathbb{Z}/\so_\blam \Z$ and~$\Partss$ is
a fixed set of representatives in~$\Parts$ under the action of~$\sigma_\mathscr{P}$.

For $\blam\in\Parts$ let $C^F_\blam$ be the corresponding cell module
and $D_\blam^F=C^F_\blam/\rad C^F_\blam$ be a graded simple $\RG(F)$-module (or
zero).  By \autoref{C:diagonalisable}, $\sigma_n^\bLam$ is $\varepsilon$-splittable, so
we can apply \autoref{T:ShiftedSimples} to deduce the following:

\begin{Lemma}\label{L:SimpleReduction}
  Suppose that $F$ is a field. Then
  $\Partspz=\set{(\blam,k)\in\Partsp|D^F_\blam\ne0}$.
\end{Lemma}

The simple $\RG(F)$-modules have been independently classified by
Bowman~\cite{Bowman:ManyCellular} and Kerschul~\cite{Kerschl:simples}.
To state their result we first need some new notation.

Let $\ltx$ be the total order on the set
of nodes $\Nodes$ where $A\ltx B$ if $\xcoord(A)<\xcoord(B)$.
Extending notation from \autoref{S:Webster}, for
$i\in\cI$ let
\[
  \Add_i(\blam)=\set{A\in\Add(\blam)|\res(A)=i} \quad\text{ and }\quad
  \Rem_i(\blam)=\set{A\in\Rem(\blam)|\res(A)=i}
\]
be the sets of addable and removable
$i$-nodes of $\blam\in\Parts$.
If $A\in\Add_i(\blam)\cup\Rem_i(\blam)$ define
\[
  d_A(\blam)=\#\set{B\in\Add_i(\blam)|A\ltx B} -\#\set{B\in\Rem_i(\blam)|A\ltx B}.
\]
Given a second node $C$ with $A\ltx C$ set
\[
  d^C_A(\blam)=\#\set{B\in\Add_i(\blam)|A\ltx B\ltx C} -\#\set{B\in\Rem_i(\blam)|A\ltx B\ltx C}.
\]
A \emph{good $i$-node} of~$\blam$ is a removable $i$-node
$A\in\Rem_i(\blam)$ that is minimal node with respect to~$\ltx$ such
that $d_A(\blam)\le0$ and~$d^C_A(\blam)<0$ whenever $A\ltx C$ and
$C\in\Rem_i(\blam)$.

\begin{Definition}
  If $n>0$ then the set of \emph{$\brho$-Uglov  $\ell$-partitions} of $n$ is
  \[
  \Uglov = \set{\blam\in\Parts|\blam=\bmu\cup\set{A}\text{ where $A$ is
  a good $i$-node of $\blam$ for some $i\in\cI$}},
  \]
  where $\Uglov[0]=\Parts[0]$.
\end{Definition}

For our particular choice of $\xcoord$-coordinate function, a special
case of the results of Bowman~\cite{Bowman:ManyCellular} and
Kerschl~\cite{Kerschl:simples} is the following:

\begin{Theorem}[{Bowman~\cite[Theorem~B]{Bowman:ManyCellular}
  and Kerschl~\cite[Main~Theorem]{Kerschl:simples}}]
  Suppose that $F$ is a field. Then
  $\set{D_\blam\<s\>|\blam\in\Uglov\text{ and } s\in\Z}$ is a complete set
  of pairwise non-isomorphic graded simple $\RG(F)$-modules.
\end{Theorem}

Hence, by \autoref{L:SimpleReduction} we obtain:

\begin{Corollary}
  Suppose that $F$ is a field. Then  $\Partspz = \set{(\blam,k) \in \Partsp | \blam \in \Uglov}$. That is,
  \(\set{D_{\blam}^{(k)}\<s\>|(\blam,k)\in\Partsp, \blam\in\Uglov\text{ and } s\in\Z}\)
  is a complete set of pairwise non-isomorphic graded simple $\Rlpn(F)$-modules.
\end{Corollary}

\section*{Index of notation}
\label{notation}
  \begin{longtable}{lp{0.7\textwidth}l}
    Symbol & Description & Page\\\hline
    \endhead
    \notation{D:ICompositions}{alpha}{An $\cI$-composition of $n$}
    \notation{E:WAA}{}[{[\alpha]}]{The $\sigma$-orbit of $\alpha$}
    \notation{D:RationalCherednik}{WA}{A rational Cherednik algebra}
    \notation{D:RationalCherednik}{WA}[[\alpha]]{A block of $\WA$}[(alpha)]
    \notation{E:WAbom}{WA}[(\bom)]{The $\bom$-weight space of $\WA$}
    \notation{E:WAA}{WAA}{$=\bigoplus_{\beta\in[\alpha]}\WA[\beta]$}
    \notation{L:BTreg}{Btreg}{A regular diagram indexed by $\t$}
    \notation{L:BTsing}{Btsing}{A singular diagram indexed by $\t$}
    \notation{D:RegSing}{}[B_\t]{The diagram $\Btreg \Btsing$}
    \notation{D:Cslam}{Cslam}{The diagram $\Cslamsing \Cslamreg$}
    \notation{E:def_Csigmakblam_sigmakt}{Ckt}{The diagram $\Cslamreg \Btreg \in \WA(\gbconf_\t^k,\gbconf_\t)$}
    \notation{D:ICompositions}{Comp}{The set of $\cI$-compositions of $n$}
    \notation{D:CST}{}[C_{\s\t}]{Cellular basis element of $\WA$}
    \notation{SS:CycHecke}{eps}{A fixed $p$th root of unity}
    \notation{E:WAbom}{Eom}{An idempotent in $\WA$}
    \notation{SS:CycHecke1}{Hln}{Hecke algebra of type $\Gln$}
    \notation{D:Hlpn}{Hlpn}{Hecke algebra of type $\Glpn$}
    \notation{D:ICompositions}{cI}{The vertex set $\set{\eps^j\q^i|i,j\in\Z}$ for the quiver $\Gamma$}
    \notation{D:multicharge}{charge}{The multicharge $\charge=(\rho_1,\dots,\rho_d)$}
    \notation{Ex:RussianDiagrams}{blam}{An $\ell$-partition of $n$}
    \notation{SS:CycHecke}{}[n,p,d=\frac\ell p]{The parameters for $\Rlpn$}
    \notation{nodes}{Nodes}{The set of nodes}
    \notation{E:so_sp}{}[\so_\lambda]{The size of the orbit of $\lambda \in \P$ under the action of $\sigma_\P$}
    \notation{S:skew_cellular_algebras}{P}{A poset}
    \notation{nodes}{Parts}{The set of $\ell$-partitions of $n$}
    \notation{C:maps_fI_gI}{Partss}{A set of $\sigma$-orbit representatives in $\Parts$}
    \notation{E:Partsp}{Partsp}{The set $\Partss \times \mathbb{Z}/\so_\blam \mathbb{Z}$}
    \notation{L:Efree}{Partspz}{The set $\set{(\blam,k)\in\Partsp|D_{\blam}^{(k),\Qcal}\ne0}$}
    \notation{D:Hlpn}{q}{The Hecke parameter}
    \notation{SS:CycHecke}{bvQ}{The cyclotomic parameters of $\Hln$}
    \notation{E:residues}{res}[\gamma]{The residue of $\gamma\in\Nodes$}
    \notation{D:QuiverHecke}{RG}{The quiver Hecke algebra of type $G(\ell,1,n)$}
    \notation{D:QuiverHecke}{RG}[{[\alpha]}]{A block of $\RG$}[{[alpha]}]
    \notation{D:Rlpn}{Rlpn}{The quiver Hecke algebra of type $G(\ell,p,n)$}
    \notation{E:sigma}{sigma}{An automorphism of $\Hln$ or $\RG$}
    \notation{D:sigma_cellular}{}[\sigma_\P]{A poset automorphism of $\P$}
    \notation{sigmaIso}{}[\sigma^\bLam_n]{An automorphism of $\RG$}
    \notation{D:sigmaP}{}[\sigma_\mathscr{P}]{An automorphism of $\Parts$}
    \notation{D:sigmaP}{}[\sigma_{\Std}]{An automorphism of $\Std(\Parts)$}
    \notation{D:standard}{Std}[(\blam)]{The set of standard tableaux of shape $\blam$}
    \notation{S:OrbitDiagrams}{Stds}{A set of $\sigma$-orbit representatives of standard tableaux of shape $\blam$}
    \notation{D:standard}{}[\s,\t]{Standard tableaux}
    \notation{E:gen-box-config}{gbconf}{A generalised partition}
    \notation{D:gbconf_Ct}{}[\gbconf_\t]{A generalised partition constructed from $\t$}
    \notation{E:WebsterDiagrams}{Web}{The set of (isotopy classes of) Webster diagrams}
    \notation{E:omega}{bom}{The $\ell$-partition $(1^n|\dots|0)\in\Parts$}
    \notation{E:xcoord}{xcoord}{The loading function $\xcoord\map\Nodes\Q$}
    \notation{D:DominanceOrder}{tdom}{The $\brho$-dominance order on $\ell$-partitions}
  \end{longtable}


\bibliographystyle{andrew}

\end{document}